\documentclass[a4paper,11pt]{article}
\usepackage[utf8]{inputenc}
\usepackage[T1]{fontenc}
\usepackage{amsthm}
\usepackage[font={footnotesize}]{caption}
\usepackage{amssymb}
\usepackage{mathrsfs}
\usepackage{amsmath}
\usepackage{amsfonts}
\usepackage{mathtools}
\usepackage{graphicx}
\usepackage{float}
\usepackage{dsfont}
\usepackage{enumerate}
\usepackage[a4paper]{geometry}
\usepackage{url}
\usepackage{csquotes}
\usepackage{IEEEtrantools}
\usepackage{appendix}
\usepackage{constants}
\usepackage{tikz-cd}
\usepackage{multicol}
\usepackage{parcolumns}
\usepackage{stackrel}
\usepackage{comment}

\usepackage{xparse}

\makeatletter
\RenewDocumentCommand\eqref{s m}{%
  \IfBooleanTF#1%
  {\textup{\tagform@{\ref*{#2}}}}
  {\textup{\tagform@{\ref{#2}}}}
}
\makeatother
\usepackage{hyperref}
\hypersetup{linktocpage,
  colorlinks   = true,
  urlcolor     = blue,
  linkcolor    = blue,
  citecolor   = blue
}

\newtheorem{theorem}{Theorem}[section]
\newtheorem{lemma}[theorem]{Lemma}
\newtheorem{proposition}[theorem]{Proposition}
\newtheorem{corollary}[theorem]{Corollary}
\theoremstyle{definition}
\newtheorem{definition}[theorem]{Definition}

\theoremstyle{remark}
\newtheorem{remark}[theorem]{Remark}

\title{\large \textbf{PERCOLATION FOR THE GAUSSIAN FREE FIELD ON THE CABLE SYSTEM: COUNTEREXAMPLES}}
\author{}
\date{}
\makeatletter
\newcommand{\Ed}{-}
\newcommand{\E}{\mathbb{E}}

\newcommand{\R}{\mathbb{R}}
\newcommand{\Z}{\mathbb{Z}}
\newcommand{\N}{\mathbb{N}}
\newcommand{\SV}{\mathcal{S}}
\newcommand{\SK}{\mathcal{SK}}
\newcommand{\KS}{\mathcal{KS}}

\newcommand{\F}{\mathcal{F}}
\newcommand{\G}{\mathcal{G}}
\newcommand{\K}{\mathcal{K}}

\newcommand{\W}{\mathcal{W}}
\renewcommand{\P}{\mathbb{P}}
\newcommand{\T}{\mathbb{T}}
\newcommand{\eps}{\varepsilon}
\newcommand{\1}{\mathds{1}}

\newcommand{\I}{{\cal I}}

\renewcommand{\phi}{\varphi}
\renewcommand{\tilde}{\widetilde}
\renewcommand{\hat}{\widehat}
\renewcommand{\epsilon}{\varepsilon}

\@addtoreset{equation}{section}

\newcommand{\tend}[2]{\displaystyle\mathop{\longrightarrow}_{#1\rightarrow#2}}
\newconstantfamily{c}{symbol=c}
\newcommand{\diff}{\,\mathrm{d}}
\newcommand{\h}{\mathit{\mathbf{h}}}
\makeatother

\usepackage{titletoc}

\dottedcontents{section}[4em]{}{2.9em}{0.7pc}
\dottedcontents{subsection}[0em]{}{3.3em}{1pc}
\begin{document}

\maketitle
\thispagestyle{empty}
\begin{center}
\vspace{-1.3cm}
Alexis Pr\'evost$^1$
\end{center}
\vspace{-0.1cm}
\begin{abstract}
\centering
\begin{minipage}{0.9\textwidth}
\vspace{0.3cm}
    For massless vertex-transitive transient graphs, the percolation phase transition for the level sets of the Gaussian free field on the associated continuous cable system is particularly well understood, and in particular the associated critical parameter $\tilde{h}_*$ is always equal to zero. On general transient graphs, two weak conditions on the graph $\G$ are given in \cite{DrePreRod3}, each of which implies one of the two inequalities $\tilde{h}_*\leq0$ and $\tilde{h}_*\geq0.$ In this article, we give two counterexamples to show that none of these two conditions are necessary, prove that the strict inequality $\tilde{h}_*<0$ is typical on massive graphs with bounded weights, and provide an example of a graph on which $\tilde{h}_*=\infty.$ On the way, we obtain another characterization of random interlacements on massive graphs, as well as an isomorphism between the Gaussian free field and the Doob $\h$-transform of random interlacements, and between the two-dimensional pinned free field and random interlacements.
\end{minipage}
\end{abstract}

\vspace{0.8cm}
\begin{minipage}{0.9\textwidth}
{\small
\tableofcontents
}
\end{minipage}
\vspace{5cm}
\begin{flushleft}

\noindent\rule{5cm}{0.4pt} \hfill April 2023 \\
\bigskip

$^1$University of Geneva\\
Section of Mathematics \\
24, rue du Général Dufour\\
1211 Geneva 4, Switzerland\\
\url{alexis.prevost@unige.ch}\\[2em]

\end{flushleft}

\newpage

\section{Introduction}
\label{sec:intro}

We consider level set percolation for the Gaussian free field on the cable system of transient weighted graphs. This percolation model was first studied on the corresponding discrete graph in \cite{MR914444}, where it is proven that the associated critical parameter is non-negative, and in fact the proof can be extended to the cable system. On a large class of discrete graphs, the critical parameter is actually positive and finite, see \cite{MR3053773,MR3492939,MR3765885,DrePreRod,DrePreRod2}. However on the cable system, this model often undergoes a phase transition at level zero, which was first proved in \cite{MR3502602} on $\Z^d,$ $d\geq3,$ then on trees in \cite{MR3492939,MR3765885}, and later on a large class of transient graphs, possibly with positive killing measure, in \cite{DrePreRod3}. In many cases, this phase transition is also particularly well understood in the near-critical regime, see \cite{DiWiLu} on $\Z^d,$ $d\geq3,$ or \cite{DrePreRod5} for additional results on general graphs. The very precise description of critical exponents derived therein designate the Gaussian free field on the cable system as a central model for percolation with long range correlations due to its strong integrability. Knowing that phase transition happens at level zero seems to be the root of this integrability, and we study in this article the limitations of this result on general transient weighted graph.
 
In \cite{DrePreRod3}, two simple criteria on the underlying graph are introduced, which together imply that the critical parameter associated to this percolation problem is equal to zero: condition \eqref{capcondition}, which says that the capacity of any unbounded set is infinite, and $\h_{\text{kill}}<1,$ see \eqref{defh0}, which says that with positive probability the discrete random walk on the graph will not be killed (by the killing measure). Canonical examples of graphs satisfying \eqref{capcondition} and $\h_{\text{kill}}<1$ are massless (i.e.\ with zero killing measure) vertex-transitive graphs. In this paper, we are interested in understanding the limit of this result, and we present various examples of graphs which answer positively the following questions: If either \eqref{capcondition} or $\h_{\text{kill}}<1$ do not hold, is it possible to find a graph with strictly positive or negative critical parameter? Can one find a graph such that the critical parameter is equal to zero, but \eqref{capcondition} or $\h_{\text{kill}}<1$ do not hold?

While trying to answer these questions, we will prove several results which are interesting in their own right: another characterization of random interlacements when $\h_{\text{kill}}\equiv1,$ see \eqref{disinterh0=1}, or Corollary~\ref{killedinterdes} for a more general result, and an isomorphism between the Gaussian free field and the $\h$-transform of random interlacements for any harmonic function $\h,$ see Theorem~\ref{couplingintergffh}, which holds under the same conditions as the isomorphism in \cite[Theorem~\ref*{4couplingintergff}]{DrePreRod3}. When $\h=\h_{\text{kill}},$ this corresponds to an isomorphism between the Gaussian free field and killed random interlacements, that is the trajectories in the random interlacement process which are killed, see Corollary~\ref{h0transformiso}, and when $\h$ is the potential kernel on $\Z^2,$ this corresponds to an isomorphism between the pinned free field and two-dimensional random interlacements, see Theorem~\ref{couplingintergffdim2}.

We use the same setting and notation as in \cite{DrePreRod3} that we now describe briefly, and we refer to Section~\ref{sec:notation} for details. We consider a transient weighted graph $ \G= (\overline{G},\bar{\lambda},\bar{\kappa}),$ where $\overline{G}$ is a finite or countably infinite set, the weights $\bar{\lambda}_{x,y},$ $x,y\in{\overline{G}},$ describe the rate at which the canonical jump process on $\overline{G}$ jumps to a neighbour, and the killing measure $\overline{\kappa}_x,$ $x\in{V},$ describe the rate at which it is killed. We allow the killing measure $\overline{\kappa}$ to be infinite, and, using network equivalence, we define a triplet $(G,\lambda,\kappa)$ so that $\kappa$ is finite, $\{x\in{\overline{G}}:\bar{\kappa}_x<\infty\}\subset G,$ and the restriction of the jump process on $(G,\lambda,\kappa)$ to $\{x\in{\overline{G}}:\bar{\kappa}<\infty\}$ corresponds to the jump process on $(\overline{G},\bar{\lambda},\bar{\kappa}),$ see around \cite[\eqref*{4eq:defGfinite}]{DrePreRod3} for details. When $\bar{\kappa}$ is finite, which will be the case in most of the examples in this article, this simply corresponds to the choice $(\overline{G},\bar{\lambda},\bar{\kappa})=(G,\lambda,\kappa).$ Note that our percolation results are only interesting when ${G}$ is infinite (otherwise percolation questions are always trivial), but allowing for $G$ to be finite does not add any complexity and is relevant in other results, for instance in Section~\ref{sec:Doob}.

Unless explicitly mentioned otherwise, we will assume that the jump process on $\G$ is transient. One can naturally associate to $\G$ the cable system, or metric graph, $\tilde{\G},$ corresponding to a continuous version of the graph where each edge $e=\{x,y\}$ is replaced by an open interval $I_e$ linking $x$ to $y,$ and where we add a half-open interval $I_x$ starting in each vertex $x.$ \phantomsection
\label{deftildeged}We denote by $\tilde{\G}^{\Ed}$ the subset of $\tilde{\G}$ consisting only of the edges $I_e$ for $e\in{E},$ that is removing the edges $I_x$ starting from $x$ for all $x\in{G}.$

One can then define a diffusion $X=(X_t)_{t\geq0}$ on the cable system $\tilde{\G},$ starting in $x$ under the probability $P_x^{\tilde{\G}},$ see for instance \cite{MR3152724}. It behaves like a Brownian motion inside the continuous edges and like the jump process on $\G$ on the vertices. The diffusion $X$ then stays in $\tilde{\G}$ until a time $\zeta\in{[0,\infty]},$ after which it remains in some cemetery state $\Delta,$ and, as $t\nearrow\zeta,$ either $X_t$ reaches the open end of the cable $I_x$ for some $x\in{G},$ and we say that $X$ has been killed, or $X_t$ exits every bounded and connected sets, and we say that $X$ blows up or survives. Note that, when starting in $x\in{G},$ the event that $X$ is killed corresponds to the event that the jump process on $\G$ is killed (i.e.\ by the killing measure). We define $\h_{\text{kill}}$ as the probability to be killed and $\h_{\text{surv}}$ as the probability to blow up: for all $x\in{\tilde{\G}},$
 \begin{equation}
\label{defh0}
    \h_{\text{kill}}(x)\stackrel{\mathrm{def.}}{=}P^{\tilde{\G}}_x\big((X_t)_{t\geq0}\text{ is killed}\big)\text{ and }\h_{\text{surv}}(x)\stackrel{\mathrm{def.}}{=}P^{\tilde{\G}}_x\big((X_t)_{t\geq0}\text{ blows up}\big)=1-\h_{\text{kill}}(x).
\end{equation}

The Gaussian free field on $\tilde{\G}$ is then defined under some probability $\P^G$ as the centered Gaussian field $(\phi_x)_{x\in{\tilde{\G}}}$ with covariance function given by
\begin{equation}
    \label{defGFF}
    \E^G[\phi_x\phi_y]=g(x,y),
\end{equation}
where $g(x,y)$ is the average time spent in $y$ by the diffusion $X$ under $P_x^{\tilde{\G}},$ see \eqref{Greendef} for a precise definition. We are interested in the level sets of the Gaussian free field, defined as
\begin{equation}
    \label{deflevelsets}
    E^{\geq h}\stackrel{\mathrm{def.}}{=}\{x\in{\tilde{\G}}:\phi_x\geq h\}.
\end{equation}
The critical parameter associated with the percolation of the level sets of the Gaussian free field on the cable system is defined as 
\begin{equation}
\label{defh*}
    \tilde{h}_*\stackrel{\mathrm{def.}}{=}\inf\big\{h\in\R:\,\P^{G}(E^{\geq h}\text{ contains an unbounded connected component})=0\big\},
\end{equation}
where we say that a connected set $F\subset\tilde{\G}$ is unbounded if and only if $F\cap G$ is infinite. Note that in the case of general transient graph $\G,$ it is not clear whether $E^{\geq h}$ contains an unbounded connected components with probability one for all $h<\tilde{h}_*$ or not. In \cite{DrePreRod3}, two main results are proved about the percolation of the Gaussian free field on the cable system. First,
\begin{equation}
\label{eq:h0<1thenh*>0}
    \text{if $\h_{\text{kill}}<1,$ then $\tilde{h}_*\geq0,$}
\end{equation} 
see  \cite[\eqref*{4ifhkill<1thenh_*>0}]{DrePreRod3}, where we write $\h_{\text{kill}}<1$ when $\h_{\text{kill}}(x)<1$ for all $x\in{\tilde{\G}},$ or equivalently $\h_{\text{kill}}(x)<1$ for some $x\in{\tilde{\G}}.$ This first result is actually an easy consequence of an extension of the isomorphism between random interlacements and the Gaussian free field on the cable system,  \cite[Proposition~6.3]{MR3502602}, to massive weighted graphs, see \eqref*{4couplingusualiso} and below in \cite{DrePreRod3} for details. Let us now introduce the following conditions
\begin{equation}
    \label{capcondition}
    \tag{Cap}
\mathrm{cap}(A)=\infty\text{ for all unbounded, closed, connected sets }A\subset \widetilde{\mathcal{G}},
\end{equation}
where $\mathrm{cap}(A)$ is the capacity of the set $A,$ which one can interpret as the size of the set $A$ from the point of view of the diffusion $X,$ see \eqref{defcap} and below for a precise definition, and 
\begin{equation}
\label{0bounded}
\tag{Sign}
    E^{\geq0}\text{ contains }\P^G\text{-a.s.\ only bounded connected components}.
\end{equation}
Note that if \eqref{0bounded} holds, then $\tilde{h}_*\leq 0.$ The second main result of interest from \cite{DrePreRod3} is that 
\begin{gather}
\label{eq:capfinite}
\P^G\big(\mathrm{cap}(E^{\geq0}(x_0))<\infty\big)=1\text{ for all } x_0\in{\tilde{\G}}, 
\\\text{ and therefore \eqref{capcondition} $\Longrightarrow$ \eqref{0bounded},}\label{eq:capimplies0bounded}
\end{gather}
see \cite[Theorem~\ref*{4T:main},1)]{DrePreRod3}. The implication \eqref{eq:capimplies0bounded} follows from the fact that, under condition \eqref{capcondition}, every closed and connected sets with finite capacity are bounded, and generalizes results on $\Z^d,$ $d\geq3$ from \cite{MR3502602} and on trees from \cite{MR3492939} and \cite{MR3765885}. It moreover follows from \eqref{eq:h0<1thenh*>0} and \eqref{eq:capimplies0bounded} that if \eqref{capcondition} is satisfied and $\h_{\text{kill}}<1,$ then $\tilde{h}_*=0.$

In this article, we present some (non-trivial) examples related to the implications in \eqref{eq:h0<1thenh*>0} and \eqref{eq:capimplies0bounded}, in order to understand better if the two conditions \eqref{capcondition} and $\h_{\text{kill}}<1$ are optimal or not. We sum up these examples in the following theorem. 

\begin{theorem}
\label{mainth}
\begin{enumerate}[1)]
    \item There exist graphs with $\h_{\text{kill}}\equiv1$ and $\tilde{h}_*\geq0,$ see Corollary~\ref{Cor:h_0=1andh_*=0}, and so \eqref{eq:h0<1thenh*>0} is not an equivalence.
    \item On $(d+1)$-regular trees with large enough killing measure, or on a large class of graphs with sub-exponential volume growth such that $\kappa\geq c$ and $\lambda\leq C$ we have $\tilde{h}_*<0,$ see Theorem~\ref{The:h_*<0}, and so the implication \eqref{eq:h0<1thenh*>0} is not trivial.
    \item There exists a graph for which \eqref{0bounded} holds, but for which this property cannot be directly deduced from \eqref{eq:capfinite}, see Proposition~\ref{Z20counterexample} and the discussion below \eqref{capGfinite}. In particular \eqref{capcondition} does not hold for this graph, and so the implication \eqref{eq:capimplies0bounded} is not an equivalence.
    \item There exist graphs for which \eqref{0bounded} does not hold, see Proposition~\ref{h*infinity}, and so the implication \eqref{eq:capimplies0bounded} is not trivial.
\end{enumerate}  
\end{theorem}

In order to obtain most of the examples in Theorem~\ref{mainth} a major role will be played by random interlacements, which were initially introduced on $\Z^d,$ $d\geq3,$ in \cite{MR2680403}, then on any transient massless graphs in \cite{MR2525105}, and extended to their cable system in \cite{MR3502602}. We explain in details how to extend the definition of the random interlacement process to the cable system of any transient weighted massive graphs in Section~\ref{sec:notation}, see in particular Theorem~\ref{nuexists} and its proof in Appendix~\ref{app:inter}, as a Poisson point process of doubly non-compact trajectories modulo time-shift, and we denote by $\I^u\subset\tilde{\G}$ the set of points visited by at least one of these trajectories, see Section~\ref{sec:notation} for details. Random interlacements are linked to the Gaussian free field via an isomorphism theorem, first derived in \cite{MR2892408} on discrete graphs and in \cite{MR3502602} on the cable system, and then strengthened in \cite{MR3492939} and \cite{DrePreRod3}. In order to illustrate the role of random interlacements to study percolation for the level sets of the Gaussian free field on $\tilde{\G},$ we recall the following consequence of the isomorphism theorem, which follows from \cite[Theorem~\ref*{4T:main},2)]{DrePreRod3}:
\begin{equation}
    \label{levelsetsvsIu}
    \begin{gathered}
    \text{ if \eqref{0bounded} holds, then }E^{\geq-\sqrt{2u}}\text{ has the same law as }\mathcal{C}_u\cup E^{\geq 0}\text{ under }\P^G\otimes\P^I,
    \\\text{ where }\mathcal{C}_u\text{ denotes the closure of the union of the connected components of}
    \\\text{the sign clusters $\{x\in{\tilde{\G}}:|\phi_x|>0\}$ intersecting the interlacement set $\I^u.$}
    \end{gathered}
\end{equation}

Let us now comment on the proofs and the four class of examples in Theorem~\ref{mainth} in more details. The first example consists of a class of graphs with $\h_{\text{kill}}\equiv1,$ killing measure diverging to infinity, and total weight from a vertex at generation $n$ to all its children also diverging to infinity, see \eqref{condtreeh_*=0} or \eqref{condlineh_*=0} for precise conditions. We prove that $\tilde{h}_*\geq0$ on these graphs using the coupling with random interlacements described in \eqref{levelsetsvsIu}. When $\h_{\text{kill}}\equiv1,$ one can describe random interlacements on $\tilde{\G}^{\Ed}$ as follows:
\begin{equation}
    \label{disinterh0=1}
    \begin{gathered}
    \text{if $\h_{\text{kill}}\equiv1,$ then the trace on $\tilde{\G}^{\Ed}$ of the random interlacement process has the same}
    \\\text{law as a Poisson point process with intensity }u\sum_{x\in{G}}\kappa_xP^{\tilde{\G}^{\Ed}}_x\text{ modulo time-shift,}
    \end{gathered}
\end{equation}
where $P^{\tilde{\G}^{\Ed}}_x$ is the law of the trace $X^{\tilde{\G}^{\Ed}}$ of $X$ on $\tilde{\G}^{\Ed},$ see \cite[Section~\ref*{4S:I_x}]{DrePreRod3} for a more precise description of this law and \eqref*{4traceonG} and above in \cite{DrePreRod3} for a definition of the trace of a process. The description \eqref{disinterh0=1} of random interlacements is a direct consequence of our construction of random interlacements, see Theorem~\ref{nuexists}, and its proof can be found below Corollary~\ref{killedinterdes}. In Section~\ref{sec:h0=1h_*=0}, we will use the description \eqref{disinterh0=1} to prove that, on the class of graphs that we consider, the random interlacement set either always contains a supercritical Galton-Watson tree or percolates along some fixed infinite path, and thus always contains an unbounded connected component. Using \eqref{levelsetsvsIu}, will let us, in turn, deduce that $\tilde{h}_*\geq0,$ and in fact $\tilde{h}_*=0$ under some additional conditions, see Corollary~\ref{Cor:h_0=1andh_*=0}. Note that one can find examples of such graphs with either sub-exponential volume growth or exponential volume growth, see the examples below Corollary~\ref{Cor:h_0=1andh_*=0}, and that one can also find examples of graphs with $\tilde{h}_*>0$ and $\h_{\text{kill}}\equiv1,$ see Remark~\ref{h0=1h_*=infinity}.

The second class of examples in Theorem~\ref{mainth} consist of any graphs with $\lambda\leq C,$ $\kappa\geq c,$ or just exponential decay of the Green function, see \eqref{condgvscap}, and which have either sub-exponential volume growth, see \eqref{condsizeball} for a more precise condition, or are $(d+1)$-regular tree with large enough killing measure. The proof of the inequality $\tilde{h}_*<0$ relies on a suitable renormalization scheme, see \eqref{defLk} and below. We first prove a quantitative bound on the probability that a cluster of $E^{\geq-h}$ has diameter $L$ for small, but positive, $h$ (depending on $L$ in an explicit form) by combining \eqref{levelsetsvsIu} and a result about the two-points function for the sign clusters of the Gaussian free field,  \cite[Proposition~5.2]{MR3502602}, see Lemma~\ref{iniren} for details. We then iterate this bound using our renormalization scheme to remove the dependency of $h$ on $L,$ relying on decoupling inequalities for the Gaussian free field from \cite{MR3325312}, see Lemma~\ref{lemmainduction}. One can then easily check that this last bound indeed imply that $\tilde{h}_*<0,$ and in fact we show that the probability that a component of $E^{\geq h}$ has diameter at least $L$ decays exponentially fast for some $h<0,$ see \eqref{connectingtoballdecayexpo}. Combining this with \eqref{levelsetsvsIu} leads, in turn, to similar results for the random interlacement set, see Corollary~\ref{cor:u_*>0}.

The third example in Theorem~\ref{mainth} is more challenging to find than a graph simply satisfying \eqref{0bounded} but not \eqref{capcondition}. Indeed, as noted in \cite[Remark~\ref*{4R:mainresults1},\ref*{4signwithoutcap}]{DrePreRod3}, one can easily obtain such graphs by simply adding vertices to a graph satisfying \eqref{capcondition}, and so \eqref{0bounded} by \eqref{eq:capimplies0bounded}, so that the new graph does not satisfy \eqref{capcondition} but still \eqref{0bounded}, see the beginning of Section~\ref{sec:Z20} for details. For these examples, one can however still almost directly deduce \eqref{0bounded} from \eqref{eq:capfinite}. In order to avoid this reasoning, and, in a sense, find a real counterexample to the implication \eqref{eq:capimplies0bounded}, we introduce a condition \eqref{capGfinite} on the graph $\G,$ under which no information about the percolation of $E^{\geq0}$ can be obtained from \eqref{eq:capfinite}, see the discussion below \eqref{capGfinite} for details. One can interpret this condition \eqref{capGfinite} as a stronger form of the complement of the condition \eqref{capcondition}. An example of a graph satisfying \eqref{0bounded} and \eqref{capGfinite} is the graph $\Z^{2,0},$ which correspond to the two-dimensional square lattice with infinite killing measure at the origin and zero killing measure everywhere else.  

It is easy to prove that \eqref{capGfinite} holds for $\Z^{2,0},$ and so \eqref{capcondition} does not hold, but in order to prove that \eqref{0bounded} holds, we introduce the notion of Doob-transform $\G_{\h}$ of the graph $\G,$ where $\h$ is an harmonic function on $\tilde{\G},$ see Definition \ref{harmo}. The diffusion $X$ on $\tilde{\G}_{\h}$ corresponds to a time-changed version of the usual $\h$-transform of the diffusion $X$ on $\tilde{\G},$ see for example \cite[Chapter~11]{MR2152573}, and we refer to \eqref{semigrouph} for a more precise statement. The sign clusters of the Gaussian free field on $\tilde{\G}_{\h}$ correspond to the sign clusters of the Gaussian free field on $\tilde{\G},$ see \eqref{eq:GFFh}, and so \eqref{0bounded} is equivalent on $\G$ or $\G_{\h}.$ In particular, if \eqref{capcondition} holds on $\G_{\h}$ for some harmonic function $\h,$ then \eqref{0bounded} holds on $\G$ by \eqref{eq:capimplies0bounded}, see Corollary~\ref{capimplieseverythingh}. 

\label{pagedefa}In the case of $\Z^{2,0},$ let $(\textbf{a}(x))_{x\in{\tilde{\Z}}^2}$ be the continuous potential kernel associated to the diffusion $X$ on $\tilde{\Z}^2,$ defined as in \cite[(4.15)]{MR2677157} for $x\in{\Z^2},$ with $\textbf{a}$ constant on $I_x$ for each $x\in{\Z^2},$ and linear on $I_e$ for each edge $e$ of $\Z^2.$ It is a classical result that the potential kernel $\textbf{a}$ is harmonic, see \cite[Proposition~4.4.2]{MR2677157}. We check that \eqref{capcondition} holds for $\Z^{2,0}_\textbf{a},$ and this proves that $\Z^{2,0}$ indeed satisfies \eqref{0bounded}, see Proposition~\ref{Z20counterexample}.

Finally, let us comment on the last class of examples in Theorem~\ref{mainth}. They consist of $(d+1)$-regular trees such that the length of the edge between a vertex at generation $n$ and one of its children is $1/(2\alpha^n),$ $\alpha<1.$ One can show under some conditions on $\alpha$ and $d,$ see \eqref{condondalpha}, that, for some vertex $x_0,$ $E^{\geq h}(x_0)$ is included in a supercritical Galton-Watson tree for all $h\in\R,$ and so $\tilde{h}_*=\infty,$ see Proposition~\ref{h*infinity}. One can find such trees with zero killing measure, or with $\h_{\text{kill}}\equiv1,$ depending on the choice of $\alpha,$ see Remark~\ref{endremark},\ref{alphadexist}). Note that on a large class of graphs, \eqref{0bounded} being not satisfied implies $\tilde{h}_*=\infty,$ see \cite[Corollary~\ref*{4dichotomy}]{DrePreRod3}. We actually expect this implication to be true on any transient graph, that is $\tilde{h}_*\in{(0,\infty)}$ never happens, which would explain why $\tilde{h}_*=\infty$ in Proposition~\ref{h*infinity}.

We finish this section by mentioning some interesting results that we obtain along the way. In Section~\ref{sec:Doob}, we use the notion of $\h$-transform of the graph $\G$ to prove that one can define a notion of $\h$-transform of random interlacements, see Definition \ref{defhtransforminter}, and an isomorphism between the Gaussian free field on $\G$ and the $\h$-transform of random interlacements similar to \cite[Theorem~2.4]{MR3492939}, under the same condition as in  \cite[Theorem~\ref*{4T:main},2)]{DrePreRod3}, see Theorem~\ref{couplingintergffh}. In particular, it implies a coupling similar to \eqref{levelsetsvsIu} but for the $\h$-transform of random interlacements, see \eqref{levelsetsvsIuh}. At the end of Section~\ref{sec:Doob}, we gather some interesting consequences of this isomorphism when choosing $\h=\h_{\text{kill}}$ or $\h=\h_{\text{surv}},$ see Corollaries~\ref{h0transformiso} and \ref{couplingintergffK}. 

For the moment, in order to illustrate more precisely the kind of results one can obtain using the $\h$-transform of random interlacements and the Gaussian free field, we present a theorem concerning the two-dimensional random interlacements and pinned Gaussian free field. We consider the lattice $\Z^2$ with its usual edge set $E_2,$ weights $\frac14$ between every two neighbours of $\Z^2$ and $0$ otherwise, and killing measure equal to $0,$ and $\tilde{\Z}^2$ the associated cable system. Let us also define $\Z^2_n$ the graph with same vertices and weights as $\Z^2,$ but with killing measure equal to infinity on $B_n^c$ and zero on $B_n,$ where $B_n$ is the discrete ball of radius $n$ centered at the origin. Even if the graph $\Z^2$ is not transient, one can define a pinned version of the Gaussian free field $(\phi^p_x)_{x\in{\tilde{\Z}^2}}$ under some probability $\P^{G,p}$ with covariance function given by
\begin{equation}
    \label{def2dgff}
    \begin{split}
    \E^{G,p}[\phi_x^p\phi_y^p]&=\lim_{n\rightarrow\infty}\E^G_{\tilde{\Z}^2_n}[(\phi_x-\phi_0)(\phi_y-\phi_0)]
    \\&=\lim\limits_{n\rightarrow\infty}g_{\tilde{\Z}^2_{n}}(x,y)-g_{\tilde{\Z}^2_{n}}(x,0)-g_{\tilde{\Z}^2_{n}}(y,0)+g_{\tilde{\Z}^2_{n}}(0,0)
    \end{split}
\end{equation}
for all $x,y\in{\tilde{\Z}^2}.$ The limit in \eqref{def2dgff} exists for each $x,y\in{\Z^2},$ see \cite[(2.27)]{MR3936156}, and, since one can obtain the value of $g_{\tilde{\Z}^2_{n}}(x,y),$ $x,y\in{\tilde{\Z}_n^2},$ from its value on $\Z_n^2$ by interpolation, see \cite[(2.1)]{MR3502602}, it is easy to show that the limit in \eqref{def2dgff} exists in fact for each $x,y\in{\tilde{\Z}^2}.$ In fact, one can show that $\phi^p$ corresponds to the Gaussian free field for the graph $\Z^{2,0},$ see Lemma~\ref{le:phiZ2otherdef}. 

We are interested in the percolation of the level sets $E_\textbf{a}^{p,\geq h}=\{x\in{\tilde{\Z}^2}:\,\phi_x^p\geq h\times \textbf{a}(x)\},$ $h\in\R,$ of the pinned free field on the cable system, and we denote by $E_\textbf{a}^{p,\geq h}(x_0)$ the connected component of $x_0\in{\tilde{\Z}^2}$ in $E_\textbf{a}^{p,\geq h}.$ Note that one could also consider percolation for the usual level sets of the pinned field $\{x\in{\tilde{\Z}^2}:\,\phi_x^p\geq h\}$ but the phase transition is then trivial, see Remark~\ref{endrkdim2},\ref{normalpinnedlevelsets}), and the level sets $E_\textbf{a}^{p,\geq h}$ appear more naturally in the context of the isomorphism with random interlacements, see \cite[Theorem~5.5]{MR3936156}.

In \cite{MR3475663}, a definition of random interlacements on $\Z^2$ was given using the $\textbf{a}$-transform of the random walk on $\Z^2,$ and corresponds to a Poisson soup of trajectories conditioned on never hitting the origin. We extend this definition here to a point process $\omega^{(2)}$ of trajectories on the cable system $\tilde{\Z}^2$ under some probability $\P^{I,2},$ see \eqref{defRIdim2}. We denote by $\ell_{x,u}^{(2)}$ the continuous field of local times associated with the trajectories in $\omega^{(2)}$ with label at most $u,$ and by $\I^u_2$ the associated interlacement set. For all closed sets $F\subset\tilde{\Z}^2$ such that $0\in{F}$ we have
\begin{equation}
\label{defIudim2}
    \P^{I,2}(\I^u_2\cap F=\varnothing)=\exp\big(-u\mathrm{cap}_{\tilde{\Z}^2}(F)\big),
\end{equation}
where $\mathrm{cap}_{\tilde{\Z}^2}(A),$ see \eqref{capdim2}, is an extension of the usual definition of two-dimensional capacity to the cable system, as defined in \cite[Section~6.6]{MR2677157} for instance. Note that since $\Z^2$ is recurrent, this definition of capacity differs from the usual definition of capacity in the rest of the paper, see \eqref{defequicap}. Moreover, the normalization in \eqref{defIudim2} corresponds to the one from \cite{MR3936156}, and differs by a constant from \cite{MR3475663}, but is more natural in our context.

As a consequence of our results about the $\h$-transform and \cite{DrePreRod3}, we obtain that the critical parameter associated with the percolation of $E^{p,\geq h}_\textbf{a}$ is equal to zero, as well as the law of the capacity of the level sets of the pinned Gaussian free field. We also obtain an isomorphism between the pinned Gaussian free field and two-dimensional random interlacements, which corresponds to a signed version of \cite[Theorem~5.5]{MR3936156} on the cable system. The proof of Theorem~\ref{couplingintergffdim2} appears at the end of Section~\ref{sec:Z20}.

\begin{theorem}
\label{couplingintergffdim2}
The sign clusters $E^{p,\geq 0}_\textbf{a}$ of the pinned free field in dimension two are $\P^{G,p}$-a.s.\ bounded and the level sets $E^{p,\geq h}_\textbf{a}$ contain with $\P^{G,p}$ positive probability an unbounded connected component for all $h<0.$ Moreover, for all $h,u\geq0$ and $x_0\in{\tilde{\Z}^2}$
\begin{equation}
\label{eq:laplacecapkilleddim2}
   \E^{G,p}\left[\exp\left(-u\mathrm{cap}_{\tilde{\Z}^2}\big({E}^{p,\geq h}_\textbf{a}(x_0)\cup\{0\}\big)\right)\1_{\phi_{x_0}^p\geq h\times \textbf{a}(x_0)}\right]
   =\P^{G,p}\big(\phi_{x_0}^p\geq \textbf{a}(x_0)\sqrt{2u+h^2}\big)
\end{equation}
and
\begin{equation}
\label{eqcouplingintergffdim2}
	\begin{gathered}
	\big(\phi_x^p\1_{x\notin{\mathcal{C}_u^{2}}}+\sqrt{2\ell_{x,u}^{(2)}+(\phi_x^p)^2}\1_{x\in{\mathcal{C}_u^{2}}}\big)_{x\in{\tilde{\Z}^2}}\text{ has the same law under }\P^{I,2}\otimes\P^{G,p}
	\\\text{ as }\big(\phi_x^p+\sqrt{2u}\textbf{a}(x)\big)_{x\in{\tilde{\Z}^2}}\text{ under }\P^{G,p}\text{ for all }u\geq0,
	\end{gathered}
\end{equation}
where $\mathcal{C}_u^2$ denotes the closure of the union of the connected components of the sign clusters $\{x\in{\tilde{\Z}^2}:|\phi_x^p|>0\}$ intersecting the random interlacement set $\I^u_2.$
\end{theorem}

Section~2 introduces the setup and the various objects studied in this article, and contains in particular the definition of killed and surviving random interlacements. Section~\ref{sec:h0=1h_*=0}, Section~\ref{sec:h_*<0}, Sections~\ref{sec:Doob}+\ref{sec:Z20}, and Section~\ref{sec:example} correspond to the proofs of the respective items in Theorem~\ref{mainth}, and can essentially be read independently, unless exceptionally explicitly mentioned otherwise. Section~\ref{sec:Doob} also presents results of independent interest about isomorphism theorems between the Gaussian free field and the $\h$-transform of random interlacements, see Theorem~\ref{couplingintergffh} and Corollaries~\ref{h0transformiso} and~\ref{couplingintergffK}. Appendix~\ref{app:inter} contains the proof of Theorem~\ref{nuexists} about the construction of random interlacements on the cable system of massive graphs, and Appendix~\ref{subsec:capandkappa=0} the proof of Proposition~\ref{corh} about the $\h$-transform of the diffusion $X$ on the cable system and the Gaussian free field.

\vspace{2mm}
\noindent{\bf Acknowledgments.} The author thanks Alexander Drewitz and Pierre-François Rodriguez for several useful discussions about the various problems solved in this article, as well as an anonymous referee for useful suggestions. This research is supported by the Engineering and Physical Sciences Research
Council (EPSRC) grant EP/R022615/1 and the Isaac Newton Trust (INT) grant G101121.

\section{Notation, definition and interlacements on massive graphs}
\label{sec:notation}
In this section, we introduce the setup and first recall the definition of the diffusion $X$ on the cable system \eqref{Dirichlet}, as well as its associated equilibrium measure \eqref{defequicap} and capacity \eqref{defcap}. We then describe the construction of random interlacements on the cable system of massive graphs in a slightly more general setup than usual, see Theorem~\ref{nuexists}. Finally, we define killed and surviving interlacements, and give a simple description of their law in Corollary~\ref{killedinterdes}.

We consider a massive weighted graph $\G=(\overline{G},\bar{\lambda},\bar{\kappa}),$ where $\bar{\kappa}$ is possibly infinite, to which we associate an equivalent triplet $(G,\lambda,\kappa),$ as defined in \cite[\eqref*{4eq:defGfinite}]{DrePreRod3}, for which $\kappa$ is finite. $G$ is then a finite or countable set of vertices, $\lambda=(\lambda_{x,y})_{x,y\in{V}}\in{[0,\infty)^{V\times V}}$ are called weights, and $\kappa=(\kappa_x)_{x\in{V}}\in{[0,\infty)^G}$ is called the killing measure. We assume that the associated graph with vertex set $G$ and edge set $E =\{ \{x,y \}\in{G^2} : \lambda_{x,y}>0\}$ is connected and locally finite. One can associate to the weighted graph $\G$ its canonical jump process, that is the continuous-time Markov chain on $\overline{G}$ which jumps from $x \in \overline{G}$ to $y \in \overline{G}$ at rate $\bar{\lambda}_{x,y}$ and is killed at rate $\bar{\kappa}_x.$ We define $\lambda_x=\kappa_x+\sum_{y\sim x}\lambda_{x,y}$ the total weight of a vertex $x\in{G},$ where $y\sim x$ means that $\{x,y\}\in{E}.$

The cable system $\tilde{\G}$ associated to $\G$ is defined by glueing together open segments $I_e$ with length $\rho_{e}=1/(2\lambda_{x,y}),$ $e=\{x,y\}\in{E},$ through their endpoints, and glueing the closed endpoint of half-open intervals $I_x$ with length $\rho_x=1/(2\kappa_x)$ to $x,$ $x\in{G}.$ For all $e=\{x,y\}\in{E}$ and $t\in{[0,\rho_{e}]}$ we denote by $x+t\cdot I_e$ the point of $I_e$ at distance $t$ from $x,$ that is $x=x+0\cdot I_e=y+\rho_e\cdot I_e,$ and similarly for all $x\in{G}$ and $t\in{[0,\rho_x)},$ we denote by $x+t\cdot I_x$ the point of $I_x$ at distance $t$ from $x.$ One can endow $\tilde{\G}$ with a distance $d_{\tilde{\G}},$ or simply $d$ when there is no ambiguity about the choice of the graph $\G,$ such that $d_{\tilde{\G}}(x,y)$ is the length of the shortest path between $x$ and $y$ when replacing the length of $I_e$ by $1$ for each $e\in{E\cup G},$ through some given increasing bijection $[0,\infty)\rightarrow[0,1)$ for $I_x$ when $\kappa_x=0.$ The associated metric space $\tilde{\G}$ is a Polish space, and a connected set $K$ is compact for this topology if and only $K\cap G$ is finite and $K\cap \overline{I_e}$ is a connected compact of $\overline{I_e}$ for all $e\in{E},$ and $K\cap I_x$ is a connected compact of $I_x$ for all $x\in{G}.$ A connected set $F$ is unbounded if and only if $F\cap G$ is infinite.

Let $m$ be the Lebesgue measure on $\tilde{\G},$ that is the sum of the Lebesgue measure on each $I_e,$ $e\in{E\cup G},$ with the normalization $m(I_e)=\rho_e$ for all $e\in{E\cap G},$ and $W_{\tilde{\G}}^+,$ be the set of continuous functions from $[0,\infty)$ to $\tilde{\G}\cup\{\Delta\},$ where $\Delta$ is some cemetery point, that is for each $w\in{W_{\tilde{\G}}^+}$ there exists a time $\zeta\in{[0,\infty]}$ such that $w_{[0,\zeta)}$ is continuous on $\tilde{\G}$ and $w(t)=\Delta$ for all $t\geq\zeta.$ We also define $W_{\tilde{\G}}^{\K,+}$ the set of forwards trajectories in $W_{\tilde{\G}}^+$ which are killed, that is escape $\tilde{\G}$ through some $I_x,$ $x\in{G},$ and $W_{\tilde{\G}}^{\SV,+}$ the set of forwards trajectories in $W_{\tilde{\G}}^+$ which blow up, that is exit every bounded and connected sets before time $\zeta.$ Let $X_t$ be the projection function at time $t$ for all $t\geq0,$ and $\mathcal{W}_{\tilde{\G}}^+$ the sigma-algebra generated by $X_t,$ $t\geq0.$ For all measures $\tilde{m}$ on $\tilde{\G},$ that is $\tilde{m}_{|I_e}$ is a measure on $(I_e,\mathcal{B}(I_e))$ for all $e\in{E\cup{{G}}},$ and measurable function $f:\tilde{\G}\rightarrow\R,$ we define
\begin{equation*}
    (f,f)_{\tilde{m}}\stackrel{\mathrm{def.}}{=}\sum_{e\in{E\cup G}}\int_{I_e}f^2\diff\tilde{m}_{|I_e},
\end{equation*}
$L^2(\tilde{\G},\tilde{m})=\{f:\,(f,f)_{\tilde{m}}<\infty\},$ and $(f,g)_{\tilde{m}}$ the associated Dirichlet form on $L^2(\tilde{\G},\tilde{m}).$ Let also $D(\tilde{\G},\tilde{m})\subset L^2(\tilde{\G},\tilde{m})$ be the space of function $f\in{C_0(\tilde{\G})},$ the closure of the space of functions with compact support with respect to the $\|\cdot\|_{\infty}$ norm, such that $f_{|I_e}\in{W^{1,2}(I_e,\tilde{m}_{|I_e}})$ for all $e\in{E\cup{G}}$ and 
\begin{equation*}
    \sum_{e\in E\cup{G}}\|f_{|I_e}\|_{W^{1,2}(I_e,\tilde{m}_{|I_e})}^2<\infty.
\end{equation*} 
Following \cite{MR3152724}, the canonical Brownian motion on $\tilde{\G}$ is then defined by taking probabilities $P_x^{\tilde{\G}},$ or simply $P_x$ when there is no ambiguity about the choice of the graph $\G,$ $x\in{\tilde{\G}},$ under which the process $X$ is an $m$-symmetric diffusion on $\tilde{\G}$ starting in $x$ and with associated Dirichlet form on $L^2(\tilde{\G},m)$
\begin{equation}
\label{Dirichlet}
    \mathcal{E}_{\tilde{\G}}(f,g)\stackrel{\text{def.}}{=}\frac12(f',g')_m\text{ for all }f,g\in{D({\tilde{\G}},m)}.
\end{equation}
We refer to \cite[Section~\ref*{4s:usefulresults}]{DrePreRod3} for more details and properties of the cable system $\tilde{\G}$ and its associated diffusion $X.$ If $F$ is either a subset of $G$ or a union of edges $I_e$ for $e\in{A\subset G\cup E},$ we denote by $X^F$ the trace of $X$ on $F,$ that is the time changed process with respect to the positive continuous additive
functional corresponding to the time spent by $X$ in $F,$ see above \cite[\eqref*{4traceonG}]{DrePreRod3} for details. The diffusion $X$ then behaves locally like a Brownian motion on each $I_e,$ $e\in{G\cup E},$ and
\begin{equation}
    \label{printonG}
    \begin{gathered}
    \text{the trace }Z\stackrel{\text{def.}}{=}X^G\text{ of }X\text{ on }G\text{ has the same law}\\\text{ as the canonical jump process on }(G,\lambda,\kappa).
    \end{gathered}
\end{equation}
One can also see $\{x\in{\overline{G}}:\,\bar{\kappa}_x<\infty\}$ as a subset of $G,$ and then the law of the trace of $X$ on $\{x\in{\overline{G}}:\,\bar{\kappa}_x<\infty\}$ is the same as the law of the canonical jump process on $\G,$ which justifies our choice of $(G,\lambda,\kappa).$ We also denote by $(\hat{Z}_n)_{n\in\N}$ the discrete-time skeleton of $Z;$ i.e. the sequence of elements of $G$ visited by the process $Z,$ with the convention that $\hat{Z}_n =\Delta$ for all large enough $n$ if $Z$ gets killed. 

Unless explicitly mentioned otherwise, we will from now on assume that the graph $\G$ is transient, that is that $\ell_{y}(\zeta)$ is $P_x^{\tilde{\G}}$-a.s.\ finite for all $x,y\in{\tilde{\G}}$, where $(\ell_y(t))_{y\in{\tilde{\G}},t\geq0}$ is the continuous field of local times with respect to $m$ associated with $X.$ We can now define the Green function on $\tilde{\G}$ by taking 
\begin{equation}
\label{Greendef}
g_{\tilde{\G}}(x,y)=E_x^{\tilde{\G}}[\ell_y(\zeta)]\text{ for all } x,y\in{\tilde{\G}},
\end{equation}
where $E_x^{\tilde{\G}}$ denotes expectation with respect to $P^{\tilde{\G}}_x.$ Denoting by $\phi_x,$ $x\in{\tilde{\G}},$ the coordinate maps on the space $C(\tilde{\G},\R),$ endowed with $\sigma$-algebra they generate, we define a probability $\P^G_{\tilde{\G}}$ on $C(\tilde{\G},\R)$ so that \eqref{defGFF} holds, and we call $\phi$ under $\P^G_{\tilde{\G}}$ the Gaussian free field on $\tilde{\G}.$ We often write $g$ instead of $g_{\tilde{\G}}$ and $\P^G$ instead of $\P^G_{\tilde{\G}}$ when there is no ambiguity about the choice of the graph $\G.$ We refer to \cite{MR3502602} or \cite[Section~\ref*{4subsec:GFF}]{DrePreRod3} for a description of the main properties of the Gaussian free field, and in particular on how to construct the Gaussian free field on the cable system from the discrete Gaussian free field on $G.$ 

If $A\subset\tilde{\G}$ is a set without accumulation point in $\tilde{\G},$ it follows from  \cite[Lemma~\ref*{4GA}]{DrePreRod3} that 
\begin{equation}
\label{eq:enhancements}
\begin{gathered}
    \text{there exists a unique graph }\G^{A}\text{ with vertex set }G^{A}\stackrel{\mathrm{def.}}{=}A\cup G
    \\\text{ such that the trace of }X\text{ on }G^{A}\text{ has the same law under }P_x^{\tilde{\G}^{A}}\text{ and }P_x^{\tilde{\G}}.
\end{gathered}
\end{equation}
The graph $\G^{A}$ informally corresponds to the graph on which for each $e\in{E\cup G},$ we used network equivalence to add $A\cap I_e$ as new vertices, by adapting the weights and killing measure so that each edge $I_e,$ $e\in{E\cup G},$ of $\tilde{\G}$ corresponds to a union of edges of $\tilde{\G}^A$ with total length $\rho_e.$ Since $\tilde{\G}^A$ has additional half-open intervals $I_x,$ $x\in{G^A},$ we will often consider $\tilde{\G}$ as a subset of $\tilde{\G}^{A}.$

We now define the equilibrium measure of a set $F\subset G$ by
\begin{equation}
    \label{defequicap}
    e_{F,\tilde{\G}}(x)\stackrel{\mathrm{def.}}{=}\lambda_xP_x^{\tilde{\G}}(\tilde{H}_F=\infty)\text{ for all }x\in{G},
\end{equation}
where $\tilde{H}_F=\inf\{n\geq1:\hat{Z}_n\in{F}\}$ is the first time the discrete time Markov chain $\hat{Z}$ on $G$ return in $F$ after time one, which is equal to $\infty$ if $\hat{Z}$ never returns in $F.$ When $F\subset\tilde{\G}$ is a closed set, we define the hitting time $H_F$ of $F$ by $H_F=\inf\{t\geq0:\,X_t\in{F}\},$ with $\inf\varnothing=\zeta,$ the last exit time $L_F$ of $F$ by $L_F=\sup\{t\geq0:X_t\in{F}\},$ with $\sup\varnothing=0$ and let 
\begin{equation}
\label{defpartialext}
    \hat{\partial}F=\left\{x\in{ F}:\,P_x\left(X_{L_F}=x,L_F\in{(0,\zeta)}\right)>0\right\}.
\end{equation} 
Note that $\hat{\partial}F$ does not have any accumulation point since $I_e$ contains at most two points of $\hat{\partial}F$ for all $e\in{E\cup G}.$ 
One can then extend the definition \eqref{defequicap} of the equilibrium measure to any closed set $F\subset\tilde{\G}$ by considering the graph $\G^{\hat{\partial}F}$ similarly as in \cite[\eqref*{4defequilibriumcable}]{DrePreRod3}. We can now define the capacity via
\begin{equation}
\label{defcap}
    \mathrm{cap}_{\tilde{\G}}(K)=\sum_{x\in{\hat{\partial}K}}e_{K,\tilde{\G}}(x)\text{ for all compacts }K\subset\tilde{\G}.
\end{equation}
One can then extend this definition of the capacity to any closed sets $F\subset\tilde{\G}$ by approximating $F$ by an increasing sequence of compacts, see \cite[\eqref*{4defcapinfinity}]{DrePreRod3}. Note that the capacity of a non-compact set is not necessarily equal to the sum of its equilibrium measure: for instance on $\Z^3$ with unit weights and zero killing measure, the capacity of $\Z^3$ is infinite, but $e_{\Z^3,\tilde{\Z}^3}(x)=0$ for all $x\in{\Z^3}.$ We simply write $e_{F}(x)$ instead of $e_{F,\tilde{\G}}$ and $\mathrm{cap}(F)$ instead of $\mathrm{cap}_{\tilde{\G}}(F)$ when there is no ambiguity about the choice of the graph $\G.$

We turn to the definition of random interlacements on $\tilde{\G}.$ The random interlacement measure was first defined on $\Z^d,$ $d\geq3,$ in \cite{MR2680403}, and then on any discrete transient graph with $\kappa\equiv0$ in \cite{MR2525105}. It was then extended to the cable system of $\Z^d$ in \cite{MR3502602} using the fact that one can obtain the diffusion $X$ by adding Brownian excursions on the edges to a discrete random walk on $\Z^d,$ and this proof can easily be extended to any transient graphs on which discrete random interlacements exist. Let us now recall this definition. For any $x\in{\hat{\partial} F},$ we define
\begin{equation}
    \label{defPxF}
    \text{$P^{F,\tilde{\G}}_x$ the law of $(X_{t+L_F})_{t\geq0}$ under $P_x(\cdot\,|\,X_{L_F}=x,L_F\in{(0,\zeta)}),$}
\end{equation}
or simply $P^{F}_x$ when there is no ambiguity about the choice of the graph $\G.$ Note that $X_{L_F}\in{\hat{\partial}F}$ if $L_F\in{(0,\zeta)}$ a.s. since $\hat{\partial}F$ has no accumulation points. Using similar ideas as in the proof of \cite[(1.56)]{MR2932978}, one can use the Markov property to prove that when considering only events which depend on the trace $Z$ of $X$ on $G,$ the probability $P_x^{F,\tilde{\G}}$ can be rewritten for all $F\subset\tilde{\G}$ with $\hat{\partial}F\subset G$ as follows
\begin{equation}
\label{PxFforZ}
    P^{F,\tilde{\G}}_x(Z\in{\cdot})=P^{\tilde{\G}}_x(Z\in{\cdot}\,|\,\tilde{H}_F=\infty)\text{ for all }x\in{\hat{\partial}F}.
\end{equation}
We now define the set of doubly non-compact trajectories $W_{\tilde{\G}}$ as the set of continuous functions from $\R$ to $\tilde{\G}\cup{\Delta},$ which take values in $\tilde{\G}$ between times $\zeta^-\in{[-\infty,\infty)}$ and $\zeta^+\in{(-\infty,\infty]},$ and is equal to $\Delta$ on $(\zeta^-,\zeta^+)^c.$ We denote by
$p_{\tilde{\G}}^*(w)$ the equivalence class of $w$ modulo time-shift for each $w\in{W_{\tilde{\G}}},$ and $W_{\tilde{\G}}^*=\{p_{\tilde{\G}}^*(w),\,w\in{W_{\tilde{\G}}}\}.$ We define $\mathcal{W}_{\tilde{\G}}$ the $\sigma$-algebra on $W_{\tilde{\G}}$ generated by the coordinate functions, and $\mathcal{W}_{\tilde{\G}}^*=\{A\subset W_{\tilde{\G}}^*:(p_{\tilde{\G}}^*)^{-1}(A)\in{\mathcal{W}_{\tilde{\G}}}\}.$ For each closed set $F\subset\tilde{\G}$ and $w\in{W_{\tilde{\G}}},$ we denote by $H_F(w)=\inf\{t\in\R:\,w(t)\in{F}\},$ the first hitting time of $F,$ with the convention $\inf\varnothing=\zeta^+.$ Let 
\begin{equation}
\label{defWFG}
    W_{F,\tilde{\G}}^0=\{\zeta^-<H_F=0<\zeta^+\}\text{ and }W_{F,\tilde{\G}}^*=p_{\tilde{\G}}^*(W_{F,\tilde{\G}}^0),
\end{equation}
that is $W_{F,\tilde{\G}}^*$ is the set of trajectories in $W_{\tilde{\G}}^*$ not starting in $F$ but which hit $F$ in finite time. If $F$ is compact, then $W_{F,\tilde{\G}}^*$ is simply the set of trajectories hitting $F.$ The forwards part of a trajectory $w\in{W_{\tilde{\G}}}$ is $(w(t))_{t\geq0}$ and its backwards part $(w(-t))_{t\geq0},$ and we denote by $\mathcal{W}_{K,\tilde{\G}}^0$ the set of events $B\in{\mathcal{W}_{\tilde{\G}}},$ $B\subset W_{F,\tilde{\G}}^0,$ such that $B$ is equal to the set of trajectories $w$ with forwards part in $B^+=\{(w(t))_{t\geq0}:\,w\in{B}\}$ and backwards part in $B^-=\{(w(-t))_{t\geq0}:\,w\in{B}\}.$ We define a measure $Q_{F,\tilde{\G}}$ on $\mathcal{W}_{\tilde{\G}},$ whose restriction to $\mathcal{W}_{F,\tilde{\G}}^0$ is given by
\begin{equation}
\label{defQK}
Q_{F,\tilde{\G}}\stackrel{\mathrm{def.}}{=}\sum_{x\in{\hat{\partial} F}}e_{F,\tilde{\G}}(x)P_x^{\tilde{\G}}({\cdot^+})P^{F,\tilde{\G}}_x({\cdot^-}),
\end{equation}
and such that $Q_{F,\tilde{\G}}(A)=0$ for all $A\in{\mathcal{W}_{\tilde{\G}}}$ with $A\cap W_{F,\tilde{\G}}^0=\varnothing.$ 

\begin{theorem}
\label{nuexists}
There exists a unique measure $\nu$ on $(W^*_{\tilde{\G}},\W^*_{\tilde{\G}})$ such that for all closed sets $F\subset\tilde{\G}$
\begin{equation}
    \label{definter}
    \nu(A)=Q_{F,\tilde{\G}}\big((p_{\tilde{\G}}^*)^{-1}(A)\big)\text{ for all }A\in{\mathcal{W}_{\tilde{\G}}^*},A\subset W_{F,\tilde{\G}}^*.
\end{equation}
\end{theorem}

We call $\nu$ the random interlacement measure. Theorem~\ref{nuexists} is classical when $K$ is compact and $\kappa\equiv0,$ see \cite{MR2525105} and \cite{MR3502602}, and can actually be extended to massive graphs, see \cite[Remark~\ref*{4R:nomassconversion}]{DrePreRod3}. The main interest of Theorem~\ref{nuexists} is that, contrary to the usual definition of the intensity measure $\nu,$ see for instance \cite[Theorem~1.1]{MR2680403}, it includes any closed sets $F,$ and not only compact, or finite, sets. This extension will later be useful, see Corollary~\ref{killedinterdes} and the proof of Proposition~\ref{Prop:condtreeh_*=0}. Since we could not find this statement in the literature, we will give in Appendix~\ref{app:inter} a direct proof of Theorem~\ref{nuexists}, which in particular provides a direct proof of the existence of the interlacement measure on the cable system of massive graphs (instead of using discrete interlacements and then adding Brownian excursions as in \cite{MR3502602}).

If $F$ is a closed set such that $P_x(L_F=\zeta)=0$ for all $x\in{\tilde{\G}},$ then for all $A\in{\W_{\tilde{\G}}^*}$ such that each $w\in{A}$ hits $F,$ we have by \eqref{defWFG}, \eqref{defQK} and \eqref{definter} that $\nu(A\cap W_{F,\tilde{\G}}^*)=\nu(A),$ which leads to the following description of the trajectories in the random interlacement process hitting $F.$ Start independently for each $x\in{\hat{\partial}F}$ a number $N_x$ of doubly-infinite trajectories with law $\text{Poi}(ue_F(x)),$ each trajectory hitting $x$ at time zero, with forwards trajectory having law $P_x$ and backwards trajectory law $P_x^{F}.$ Then by \eqref{definter} the point process consisting of all these trajectories modulo time-shift has the same law as the trajectories in $\omega_u$ hitting $F.$ Examples of such sets $F$ are $F=G$ if $\h_{\text{kill}}\equiv1,$ or $F=\Z^k\times\{0\}^{d-k}$ if $\G=\Z^d$ with $\kappa\equiv0$ and $k\in{\{1,\dots,d-3\}},$ which follows easily from Wiener’s test, see \cite[Theorem~2.2.5]{MR1117680}.

\begin{remark}
\begin{enumerate}[1)]
    \item Similarly as in \cite[(1.40)]{MR2680403} or \cite[(2.16)]{MR3167123}, it is easy to show that random interlacements on the cable system are invariant under time reversal. Indeed for all connected compacts $K$ of $\tilde{\G}$ we have by \eqref{defQK}, \eqref{lastexitdec} and \eqref{symbridge} that for all $A,A''\in{W_{\tilde{\G}}^+}$ and $A'\in{W_{\tilde{\G}}^{+,f}},$
\begin{align*}
    &Q_{K,\tilde{\G}}\big((X_{-t})_{t\geq0}\in{A},(X_t)_{t\in{[0,L_K]}}\in{A'},(X_t)_{t\geq L_K}\in{A''}\big)
    \\&=\sum_{x,y\in{\hat{\partial} K}}e_{K}(x)e_{K}(y)P^{K}_y(A'')P^{K}_x(A)\int_{0}^{\infty}P_{x,y,s}(\pi_s^{-1}(A'))p_s(x,y)\diff s
    \\&=Q_{K,\tilde{\G}}\big((X_{-t})_{t\geq0}\in{A''},(X_{L_K-t})_{t\in{[0,L_K]}}\in{A'},(X_t)_{t\geq L_K}\in{A}\big).
\end{align*}
 Denoting by $\check{\nu}$ the image of $\nu$ under time reversal, taking a sequence of compacts increasing to $\tilde{\G},$ we thus directly obtain by \eqref{definter} that
\begin{equation}
\label{nuinvrev}
    \nu=\check{\nu}.
\end{equation}
    \item One can find a result similar to \eqref{definter} but with
    \begin{equation}
    \label{defWFGprime}
        W_{F,\tilde{\G}}^{'0}=\{\zeta^-<0=L_F<\zeta^+\}\text{ and }W_{F,\tilde{\G}}^{'*}=p_{\tilde{\G}}^*(W_{F,\tilde{\G}}^{'0})
    \end{equation}
    instead of $W_{F,\tilde{\G}}^{0}$ and $W_{F,\tilde{\G}}^{*},$ see \eqref{defWFG}, where $L_F(w)=\sup\{t\in\R:\,w(t)\in{F}\}$ for all $w\in{W_{\tilde{\G}}},$ with the convention $\sup\varnothing=\zeta^-.$ Indeed defining $\W_{F,\tilde{\G}}^{'0}$ as the set of events $B\in{\W_{\tilde{\G}}},$ $B\subset W_{F,\tilde{\G}}^{'0},$ which are a product of $B^+$ and $B^-,$ and taking
    \begin{equation}
\label{defQKprime}
Q_{F,\tilde{\G}}'\stackrel{\mathrm{def.}}{=}\sum_{x\in{\hat{\partial} F}}e_{F,\tilde{\G}}(x)P_x^{\tilde{\G}}({\cdot^-})P^{F,\tilde{\G}}_x({\cdot^+}),
\end{equation}
we have by \eqref{nuinvrev} that
\begin{equation}
    \label{definterprime}
    \nu(A)=Q_{F,\tilde{\G}}'\big((p_{\tilde{\G}}^*)^{-1}(A)\big)\text{ for all }A\in{\mathcal{W}_{\tilde{\G}}^*},A\subset W_{F,\tilde{\G}}^{'*}.
\end{equation}
\end{enumerate}
\end{remark}

We now define the random interlacement process $\omega$ under some probability $\P^I_{\tilde{\G}},$ or $\P^I$ when there is no ambiguity about the choice of the graph $\G,$ as a Poisson point process with intensity measure $\nu\otimes\lambda,$ where $\lambda$ is the Lebesgue measure on $(0,\infty).$ We also denote by $\omega_u$ the point process, which consists of the trajectories in $\omega$ with label less than $u,$ by $(\ell_{x,u})_{x\in{\tilde{\G}}}$ the continuous total field of local times with respect to $m$ on $\tilde{\G}$ of $\omega_u$ and by $\I^u=\{x\in{\tilde{\G}}:\ \ell_{x,u}>0\}$ the interlacement set at level $u.$ The trace $\hat{\omega}_u\stackrel{\mathrm{def.}}{=}\omega_u^{G}$ of $\omega_u$ on $G$ corresponds to a random interlacement process on the discrete graph $\G,$ and the random interlacement set $\I^u$ is characterized by the following relation
\begin{equation}
    \label{defIu}
    \P^I(\I^u\cap F=\varnothing)=\exp(-u\mathrm{cap}(F))\text{ for all closed sets $F$}.
\end{equation}

When $0<\h_{\text{kill}}<1,$ see \eqref{defh0}, there are four type of trajectories in the interlacement process: either their forwards and backwards parts are killed, or their forwards and backwards parts blow up, or their backwards parts are killed and forwards parts blow up, or their forwards parts are killed and backwards parts blow up, and we denote the respective sets of trajectories by $W_{\tilde{\G}}^{\K,*},$ $W_{\tilde{\G}}^{\SV,*},$ $W_{\tilde{\G}}^{\KS,*}$ and $W_{\tilde{\G}}^{\SK,*}.$ We call respectively killed random interlacements, surviving random interlacements, killed-surviving random interlacements and surviving-killed random interlacements the point processes corresponding to each type of trajectory, which are Poisson point processes with respective intensity measure $\nu^{\K},$ $\nu^{\SV},$ $\nu^{\KS}$ and $\nu^{\SK}$ (that is $\nu^{\K}(A)=\nu(A\cap W_{\tilde{\G}}^{\K,*})$ for instance). 

These notions of killed and surviving random interlacements are not really necessary here to prove Theorem~\ref{mainth}, but they seem to be a natural object to study on massive graphs and let us obtain a larger class of massive graphs with $\tilde{h}_*\geq0,$   see the proof of Proposition~\ref{Prop:condtreeh_*=0}. Moreover, in Section~\ref{sec:Doob}, we are going to prove general results about random interlacements, which in turn will be useful to find examples as in Theorem~\ref{mainth}, but are also interesting in their own rights, see for instance Theorem~\ref{couplingintergffh}. We will then use killed or surviving interlacements to give a first illustration of the independent interest of these results, see Corollaries \ref{h0transformiso} and \ref{couplingintergffK}. We also refer to   \cite[Section~6]{DrePreRod5} for an example of a proof where killed-surviving interlacements play an essential role. The main interest of these objects is that their intensity measure can be described in a much easier fashion than the usual intensity measure of random interlacements, see \eqref{definter}. Indeed,  defining $G_{\kappa}=\{x\in{G}:\,\kappa_x>0\}$ and recalling the definition of $\tilde{\G}^{\Ed}$ from page~\pageref{deftildeged}, Theorem~\ref{nuexists} implies the following.

\begin{corollary}
\label{killedinterdes}
Let $\tilde{\nu}^\K,$ $\tilde{\nu}^{\KS}$ and $\tilde{\nu}^{\SK}$ be the probabilities on $(W_{\tilde{\G}},\mathcal{W}_{\tilde{\G}})$ given by 
\begin{align*}
    \nonumber\tilde{\nu}^{\K}\stackrel{\mathrm{def.}}{=}\sum_{x\in{G_{\kappa}}}\kappa_x\h_{\text{kill}}(x)P^{\tilde{\G}}_x\big(\cdot^+\,|\,\mathcal{W}_{\tilde{\G}}^{\K,+}\big)P^{\overline{I_x^c},\tilde{\G}}_x\big(\cdot^-\big)\text{ on }\mathcal{W}_{\tilde{\G}^{\Ed},\tilde{\G}}^{0},
    \\\tilde{\nu}^{\KS}\stackrel{\mathrm{def.}}{=}\sum_{x\in{G_{\kappa}}}\kappa_x\h_{\text{surv}}(x)P^{\tilde{\G}}_x\big(\cdot^+\,|\,\mathcal{W}_{\tilde{\G}}^{\SV,+}\big)P^{\overline{I_x^c},\tilde{\G}}_x\big(\cdot^-\big)\text{ on }\mathcal{W}_{\tilde{\G}^{\Ed},\tilde{\G}}^{0},
    \\\nonumber\tilde{\nu}^{\SK}\stackrel{\mathrm{def.}}{=}\sum_{x\in{G_{\kappa}}}\kappa_x\h_{\text{surv}}(x)P^{\overline{I_x^c},\tilde{\G}}_x\big(\cdot^+\big)P^{\tilde{\G}}_x\big(\cdot^-\,|\,\mathcal{W}_{\tilde{\G}}^{\SV,+}\big)\text{ on }\mathcal{W}_{\tilde{\G}^{\Ed},\tilde{\G}}^{'0},
\end{align*}
and such that $\tilde{\nu}^{\K}=\tilde{\nu}^\KS_{\tilde{\G}}=0$ for all $A\in\mathcal{W}_{\tilde{\G}}$ with $A\cap W_{\tilde{\G}^{\Ed},\tilde{\G}}^0=\varnothing,$ and $\tilde{\nu}^\SK_{\tilde{\G}}=0$ for all $A\in\mathcal{W}_{\tilde{\G}}$ with $A\cap W_{\tilde{\G}^{\Ed},\tilde{\G}}^{'0}=\varnothing.$ Then 
\begin{equation}
\label{newdescriptionnuK}
    \nu^{\K}(A)=\tilde{\nu}^{\K}\big((p_{\tilde{\G}}^*)^{-1}(A)\big)\text{ for all }A\in{\mathcal{W}_{\tilde{\G}}^{*}},A\subset W^{\K,*}_{\tilde{\G}}\cap W^{*}_{\tilde{\G}^{\Ed},\tilde{\G}},
\end{equation}
and similarly for killed-surviving random interlacements if $A\subset W^{\KS,*}_{\tilde{\G}}\cap W^{*}_{\tilde{\G}^{\Ed},\tilde{\G}},$ and for surviving-killed random interlacements if $A\subset W^{\SK,*}_{\tilde{\G}}\cap W^{'*}_{\tilde{\G}^{\Ed},\tilde{\G}}$.
\end{corollary}
\begin{proof}
Let us first consider killed random interlacements. We have $\hat{\partial}\tilde{\G}^{\Ed}=G_{\kappa},$ $e_{\tilde{\G}^{\Ed}}(x)=e_{G}(x)=\lambda_xP_x(\tilde{H}_G=\infty)=\kappa_x$ and $P^{\tilde{\G}^{\Ed}}_x=P^{\overline{I_x^c}}_x$ since $L_{\overline{I_x^c}}=L_{\tilde{\G}^{\Ed}}$ on the event $\{L_{\tilde{\G}^{\Ed}}\in{(0,\zeta)},X_{L_{\tilde{\G}^{\Ed}}}=x\}=\{L_{\overline{I_x^c}}\in{(0,\zeta)},X_{L_{\overline{I_x^c}}}=x\}$ for all $x\in{G_{\kappa}}.$ Therefore $\tilde{\nu}^{\K}(A)=Q_{\tilde{\G}^{\Ed},\tilde{\G}}(A)$ for all sets $A\subset W_{\tilde{\G}}$ with $A^+,A^-\subset W_{\tilde{\G}}^{\K,+},$ and so \eqref{newdescriptionnuK} follows readily from \eqref{definter}. The proof is similar for killed-surviving interlacements, as well as for surviving-killed interlacements using \eqref{defQKprime} and \eqref{definterprime} with $F=\tilde{\G}^{\Ed}.$ 
\end{proof}

Corollary~\ref{killedinterdes} provides us with the following description of the point process consisting of the trajectories of killed interlacements hitting $\tilde{\G}^{\Ed}$: for each $x\in{G},$ take a Poisson number of trajectories with parameter $u\kappa_x\h_{\text{kill}}(x),$ each independent, with law $P_x(\cdot\,|\,W^{\K,+}_{\tilde{\G}})$ for their forwards part and $P_x^{\overline{I_x^c}}$ for their backwards part. Then the point process which consist of all these trajectories modulo time-shift has the same law as the point process consisting of the trajectories of killed interlacements hitting $\tilde{\G}^{\Ed}.$ Note that $P_x^{\overline{I_x^c}}$ can also be seen as the law as a $\text{BES}^3(0)$ process on $I_x$ starting in $x$ and stopped when reaching the open end of $I_x,$ see for instance \cite[Theorem~4.5, Chapter~XII]{MR1725357}. We will see another description of killed random interlacements as an $\h_{\text{kill}}$-transform of random interlacements in Corollary~\ref{h0transformiso}.

When $\h_{\text{kill}}\equiv1$ killed random interlacements and random interlacements coincide, and one can thus describe the trace on $\tilde{\G}^{\Ed}$ of the random interlacement process, that is the point process consisting of the trace on $\tilde{\G}^{\Ed}$ of each trajectory in the random interlacement process $\omega_u$ hitting $\tilde{\G}^{\Ed},$ as in \eqref{disinterh0=1}. One can describe similarly the discrete killed random interlacement process, that is the point process consisting of the trajectories in $\hat{\omega}_u$ whose forwards and backwards parts are killed, as well as killed-surviving and surviving-killed random interlacements. Note that finitary interlacements, as introduced in \cite{MR3962876}, are a special case of killed random interlacements, and \eqref{disinterh0=1} can be seen as generalization of \cite[Proposition~4.1]{MR3962876}, and we refer to Remark~\ref{linkwithfinitary} for additional results on finitary random interlacements.

\begin{remark}
\label{rqkilled}
\begin{enumerate}[1)]
    \item \label{desgiveall}If one applies \eqref{newdescriptionnuK} to a new graph $\G'$ which is the same graph as $\G,$ plus an additional vertex $x+t_x\cdot I_x$ on each $I_x,$ $x\in{G}$ and $t_x\in{(0,\rho_x)},$ then \eqref{newdescriptionnuK} describes the law of $\nu^\K$ on $(\tilde{\G}')^{\Ed},$ that is on $\tilde{\G}^{\Ed}$ and on $[x,x+t_x\cdot I_x](\subset I_x),$ $x\in{G}.$ We can approximate the whole cable system $\tilde{\G}$ in that way by letting $t_x\rightarrow\rho_x$ for all $x\in{G},$ and thus \eqref{newdescriptionnuK} is enough to obtain the complete law of $\nu^{\K}.$ One cannot however find a direct description similar to \eqref{newdescriptionnuK} for the complete law of $\nu^{\K}$ since for all $x\in{G}$ with $\kappa_x>0,$ 
    \begin{equation*}
        \nu^{\K}([x,x+t\cdot I_x])\geq e_{[x,x+t\cdot I_x]}(x+t\cdot I_x)\h_{\text{kill}}(x+t\cdot I_x)\rightarrow{\infty}
    \end{equation*}
    as $t\nearrow\rho_x$ by a similar argument as in the proof of \cite[\eqref*{4capIx}]{DrePreRod3}, and so there is an infinite number of trajectories in the killed random interlacement process hitting $I_x$ (or in fact $[x+t\cdot I_x,x+\rho_x\cdot I_x)$ for any $t<\rho_x$).

    \item\label{rk:killedonARI} For a set $A\subset G,$ we say that a trajectory on $\tilde{\G}$ is killed on $A$ if it reaches the open end of the cable $I_x$ for some $x\in{A}.$ One could also want to study "killed on $A$" random interlacements, that is the point process consisting of the trajectories in the random interlacement process $\omega$ whose forwards and backwards trajectories have been killed on $A,$ or "surviving on $A^c$" random interlacements, that is the point process consisting of the trajectories in the random interlacement process $\omega$ whose forwards and backwards trajectories have not been killed on $A.$ 

Let $A_{\infty}=\{x+\rho_x(1-2^{-n})\cdot I_x:\,x\in{A^c\cap G},n\in\N\},$ and $\G^{A_{\infty}}$ the graph defined in \eqref{eq:enhancements}. Since $\kappa^{A_{\infty}}_x=0$ for all $x\in{A_{\infty}\cup (G\setminus A)},$ killed trajectories in $W_{\tilde{\G}^{A_{\infty}}}^+$ are always killed on $I_x$ for some $x\in{A}.$ Therefore, the trace on $\tilde{\G}(\subset\tilde{\G}^{A_{\infty}})$ of the killed random interlacement process under $\P^{I}_{\tilde{\G}^{A_{\infty}}}$ has the same law as the killed on $A$ random interlacement process under $\P^I_{\tilde{\G}},$ and the trace on $\tilde{\G}(\subset\tilde{\G}^{A_{\infty}})$ of the surviving random interlacement process under $\P^{I}_{\tilde{\G}^{A_{\infty}}}$ has the same law as the surviving on $A^c$ random interlacement process under $\P^I_{\tilde{\G}}.$ Therefore in the sequel we will mainly focus on killed and surviving random interlacements, since all our definitions and results could be extended to killed on $A,$ or surviving on $A^c,$ random interlacements by considering the graph $\G^{A_{\infty}}$ instead of $\G.$ 
    
If now $A\subset \overline{G},$ that is we allow $\bar{\kappa}_x=\infty$ for some $x\in{A},$ we say that a trajectory on $\tilde{\G}$ is killed on $A$ if this trajectory is killed on $A'$ for the graph $(G,\lambda,\kappa),$ where $A'\subset G$ is the union of $A\cap G$ and all the vertices $z\in{G}$ for which there exist $x\in{A\cap G^c}$ and $y\in{G}$ with $\bar{\lambda}_{x,z}>0$ and $\bar{\lambda}_{z,y}>0.$ We adapt the definition of killed on $A$ and surviving on $A^c$ random interlacements accordingly, and similarly as before, it is enough to study killed, or surviving, random interlacements to obtain results on killed on $A,$ or surviving on $A^c,$ random interlacements, even if $A\subset\overline{G},$ which we will from now on do.

    \item \label{killedonA} If $A\subset G,$ Corollary~\ref{killedinterdes} still holds for killed on $A$ random interlacements instead of killed random interlacements when replacing $G_{\kappa}$ by $G_{\kappa}^{A_{\infty}}=\{x\in{A}:\,\kappa_x>0\},$ killed trajectories $W_{\tilde{\G}}^{\K,+}$ by killed on $A$ trajectories, and $\h_{\text{kill}}$ by the probability to be killed on $A,$ which follows from considering the graph $\G^{A_{\infty}}$ instead of $\G$.
\end{enumerate}
\end{remark}

\section{Some trees with \texorpdfstring{$\h_{\text{kill}}\equiv1$}{hkill=1} and \texorpdfstring{$\tilde{h}_*=0$}{h*=0}}
\label{sec:h0=1h_*=0}
We present a class of weighted graphs for which $\h_{\text{kill}}\equiv1$ and the critical parameter $\tilde{h}_*=0,$ thus proving that the implication \eqref{eq:h0<1thenh*>0} is not an equivalence. First, we use the description \eqref{disinterh0=1} of random interlacements when $\h_{\text{kill}}\equiv1$ (or more precisely Corollary~\ref{killedinterdes}) and the link between the Gaussian free field and random interlacements from \eqref{eq:capimplies0bounded} to prove that $\tilde{h}_*\geq0$ under certain conditions \eqref{condtreeh_*=0} or \eqref{condlineh_*=0} on the growth of the weights, killing measure and number of neighbours in a subset of our graph in Proposition~\ref{Prop:condtreeh_*=0}, which combined with a simple condition which implies $\h_{\text{kill}}=1$ provides us with the desired class of graphs, see Corollary~\ref{Cor:h_0=1andh_*=0}. We say that $\G$ contains a tree $(T_i)_{i\geq0}$ if $T_i,$ $i\geq0,$ is a sequence of disjoints subsets of $G$ such that  $|T_0|=1$ and each vertex of $T_{i+1}$ has a unique neighbour in $T_i$ for all $i\geq0,$ and we say that $(T_i)_{i\geq0}$ is binary if $|\{y\in{T_{n+1}:\,y\sim x}\}|=2$ for all $x\in{T_n}.$ 

\begin{proposition}
\label{Prop:condtreeh_*=0}
Let $\G$ be a transient graph such that either
\begin{gather}
\label{condtreeh_*=0}
    \text{$\G$ contains a tree }(T_i)_{i\geq0}\text{ with }\inf_{x\in{T_n}}\sum_{\substack{y\in{T_{n+1}}\\y\sim x}}1\wedge\frac{\lambda_{x,y}\kappa_x}{\lambda_x}\tend{n}{\infty}\infty,
    \\\label{condlineh_*=0}
    \text{or there exists an infinite connected path } (x_1,x_2,\dots)\subset G^\N\text{ with }\frac{\lambda_{x_n,x_{n+1}}\kappa_{x_n}}{\lambda_{x_n}\log(n)}\tend{n}{\infty}\infty,
\end{gather}
then $\tilde{h}_*\geq0.$
\end{proposition}
\begin{proof}

Let us first assume that \eqref{condtreeh_*=0} holds. For each $x\in{T_n}$ and $y\in{T_{n+1}},$ the probability that a discrete time trajectory starting at $x$ first jump to $y$ is $\lambda_{x,y}/\lambda_x.$ Starting an independent Poisson number of discrete trajectories at $x$ with parameter $u\kappa_x$ for all $x\in{T_i},$ $i\geq0,$ let us denote by $A_u\subset G$ the union of the set of vertices at generation $n+1$ visited at time one by a trajectory starting at generation $n,$ $n\in\N.$ The average number of neighbours $y\in{T_{n+1}\cap A_u}$ of a vertex $x\in{T_n}$ is then
\begin{equation*}
    \sum_{\substack{y\in{T_{n+1}}\\y\sim x}}1-\exp\Big(-\frac{u\lambda_{x,y}\kappa_x}{\lambda_x}\Big)\geq \frac12\sum_{\substack{y\in{T_{n+1}}\\y\sim x}}1\wedge \frac{u\lambda_{x,y}\kappa_x}{\lambda_x}\geq\frac{u\wedge1}{2}\sum_{\substack{y\in{T_{n+1}}\\y\sim x}}1\wedge \frac{\lambda_{x,y}\kappa_x}{\lambda_x}.
\end{equation*}
Therefore, if $n$ is large enough, the intersection of $A_u$ and all the descendants of a vertex at generation $n$ stochastically dominates a supercritical Galton-Watson tree, and thus contains an infinite connected component with positive probability. If we denote by $\tilde{A}_u\subset\tilde{\G}$ the set obtained by adding the cable between each $x\in{A_u}$ and its first ancestor, we thus have that $\tilde{A}_u$ contains an unbounded connected component with positive probability. It moreover follows from Corollary~\ref{killedinterdes} that the trace on $\tilde{\G}^-$ of the sum of the killed random interlacements process and the killed-surviving random interlacements process has the same law as a Poisson point process with intensity $\sum_{x\in{G}}u\kappa_xP_x^{\tilde{\G}^-}$, and so  $\I^u\cap\tilde{\G}^{\Ed}$ stochastically dominates $\tilde{A}_u.$ Therefore, $\I^u$ contains an unbounded connected component with positive probability for all $u>0.$ If \eqref{0bounded} holds, it moreover follows from \eqref{levelsetsvsIu} that $\I^u$ is stochastically dominated by $E^{\geq-\sqrt{2u}},$ and so $E^{\geq-\sqrt{2u}}$ contains an unbounded connected component with positive probability for all $u>0,$ that is $\tilde{h}_*\geq0.$ If \eqref{0bounded} does not hold, it is clear that $\tilde{h}_*\geq0$ by monotonicity.

Let us now assume that \eqref{condlineh_*=0} holds. For all $u>0$ we have
\begin{equation*}
    \sum_{n\in\N}\exp\Big(-\frac{u\lambda_{x_n,x_{n+1}}\kappa_{x_n}}{\lambda_{x_n}}\Big)=\sum_{n\in\N}\left(\frac{1}{n}\right)^{\frac{u\lambda_{x_n,x_{n+1}}\kappa_{x_n}}{\lambda_{x_n}\log(n)}}<\infty.
\end{equation*}
Therefore by Corollary~\ref{killedinterdes} and Borel-Cantelli lemma, there exists a.s.\ $N\in\N$ such that the random interlacements process at level $u$ on $\G$ contains for all $n\geq N$ a trajectory starting in $x_n$ and visiting $x_{n+1}$ at time $1,$ and so $\I^u$ contains an unbounded connected component, and we can conclude similarly as above.
\end{proof}

It is moreover easy to find a condition which implies that $\h_{\text{kill}}\equiv1$: it is enough that the probability to be instantly killed for a discrete trajectory is uniformly larger than a constant. 

\begin{corollary}
\label{Cor:h_0=1andh_*=0}
Let $\G$ be a graph satisfying either \eqref{condtreeh_*=0} or \eqref{condlineh_*=0} such that $\kappa_x\geq c$ and $\frac{\kappa_x}{\lambda_x}\geq c'$ for all $x\in{G}.$ Then $\tilde{h}_*=0$ and $\h_{\text{kill}}\equiv1.$
\end{corollary}
\begin{proof}
Since $\frac{\kappa_x}{\lambda_x}\geq c',$ the canonical random walk on $\G$ has probability at least $c'>0$ to be killed at each step, and so  $\h_{\text{kill}}\equiv1.$ Moreover, since for all compacts $K\subset G,$ we have $e_K(x)=\lambda_xP_x(\tilde{H}_K=\infty)\geq\kappa_x\geq c$ for all $x\in{K},$ we deduce that $\mathrm{cap}(A)\geq c|A|$ for all $A\subset G.$ Thus by \cite[\eqref*{4capconditiondis}]{DrePreRod3} we have that \eqref{capcondition} holds, and so $\tilde{h}_*\leq0$ by \eqref{eq:capimplies0bounded}. We can now conclude by Proposition~\ref{Prop:condtreeh_*=0}. 
\end{proof}

Let us give three canonical examples of graphs satisfying the conditions of Corollary~\ref{Cor:h_0=1andh_*=0}.
\begin{enumerate}
    \item any trees with constant weights such that the number $N_n$ of children of each vertex at generation $n$ and the killing measure $\kappa_n$ of each vertex at generation $n$ satisfy $\kappa_n\geq cN_n$ and $N_n\tend{n}{\infty}\infty,$
    \item any trees with bounded degree such that the weight $\lambda_n$ of each edge between a vertex at generation $n$ and a vertex at generation $n+1$ and the killing measure $\kappa_n$ of each vertex at generation $n$ satisfy $\kappa_n\geq\lambda_n\wedge\lambda_{n-1}$ and $\lambda_n\tend{n}{\infty}\infty,$
    \item any graphs $\G$ with $\kappa_x\geq c$ and $\kappa_x/\lambda_x\geq c'$ for all $x\in{G}$ containing an infinite connected path $(x_1,x_2,\dots)\subset G^\N$ with  $\lambda_{x_n,x_{n+1}}/\log(n)\tend{n}{\infty}\infty.$
\end{enumerate}
Note that all these examples have an unbounded killing measure, and in fact one can easily show that this is always the case for any graphs satisfying either \eqref{condtreeh_*=0} or \eqref{condlineh_*=0}. Example~(c) shows that one can find graphs with sub-exponential volume growth satisfying $\tilde{h}_*=0$ and $\h_{\text{kill}}\equiv1$  if the killing measure and weights increase sufficiently fast along some path, and example~(b) shows that it is in fact enough that they diverge at the same speed to infinity on trees with bounded degree. We refer to Remark~\ref{rkboundedlambda} for an explanation of why these are important examples in view of Theorem~\ref{The:h_*<0}. Finally, example~(a) shows that it is also possible to find  graphs satisfying $\tilde{h}_*=0$ and $\h_{\text{kill}}\equiv1$ with bounded weights $\lambda_{x,y},$ $x,y\in{G},$ or equivalently bounded lengths $\rho_{e},$ $e\in{E},$ of the cables, but with unbounded degree.

\begin{remark}
\label{h0=1h_*=infinity}
In Proposition~\ref{h*infinity}, we also give an example of a graph satisfying $\h_{\text{kill}}\equiv1$ for which $\tilde{h}_*=\infty$ (this is the case when $\alpha\leq\frac1d$), which also show that the implication \eqref{eq:h0<1thenh*>0} is not an equivalence. The advantage of Corollary~\ref{Cor:h_0=1andh_*=0} is that it shows that even under condition \eqref{capcondition}, under which $\tilde{h}_*\leq0,$ the implication $\h_{\text{kill}}<1\ \Rightarrow\ \tilde{h}_*=0$ is still not an equivalence. 
\end{remark}

\section{Proof of \texorpdfstring{$\tilde{h}_*<0$}{h*<0} on massive graphs}
\label{sec:h_*<0}
In this section we prove that $\tilde{h}_*<0$ on almost any massive graph with sub-exponential volume growth, bounded weights and $\kappa\geq c,$ and on a class of trees including the $(d+1)$-regular tree with unit weights and large enough constant killing measure, see Theorem~\ref{The:h_*<0} and Corollary~\ref{Cor:h_*<0} for exact statements. This shows in particular that the implication \eqref{eq:h0<1thenh*>0} is not trivial, and in fact that one has $\tilde{h}_*<0$ on most typical examples of graphs with $\h_{\text{kill}}\equiv1.$ Finally, we use the isomorphism between random interlacements and the Gaussian free field to derive similar results for the interlacement set on the cable system in Corollary~\ref{cor:u_*>0}.

Throughout this section, we will use $c$ and $C$ for constants changing from place to place, and numbered constants $c_0,$ $C_0$... for fixed constants, which appear in increasing numerical order. For each $L\geq0$ and $x\in{G},$ let us define the discrete ball $B(x,L)=\{y\in{G}:d(x,y)< L\}$ with internal boundary $\partial B(x,L)=\{y\in{B(x,L)}:\,\exists\,z\sim y,z\notin{B(x,L)}\},$
\begin{equation}
\label{defgb}
    g(L)=\sup_{x\in{G},y\in{B(x,L)^c}}g(x,y)\text{ and }b(L)=\sup_{x\in{G}}|\partial B(x,L)|.
\end{equation}

\begin{theorem}
\label{The:h_*<0}
Assume that $\lambda_x\leq C$ for all $x\in{G},$ and that there exist $C<\infty$ and $c_0>0$ such that
\begin{equation}
    \label{condgvscap}
    g(L)\leq C\exp(-c_0L)\text{ for all }L\geq0.
\end{equation}
Assume moreover that either there exist $\alpha>0$ and $c>0$ such that
\begin{equation}
\label{condsizeball}
    b(L)\leq \exp\left(\frac{cL}{\log(L)^{1+\alpha}}\right)\text{ for all }L\geq2,
\end{equation}
or $\G$ is a tree and there exist $C<\infty$ and $c_1\in{(0,c_0\cdot c_2)},$ where $c_2\in{(0,1]}$ is some absolute constant independent of the choice of $\G,$ such that
\begin{equation}
    \label{condsizeballtree}
    b(L)\leq C\exp\left(c_1L\right)\text{ for all }L\geq1.
\end{equation}
Then there exist $h<0,$ $C<\infty$ and $c>0$ such that
\begin{equation}
\label{connectingtoballdecayexpo}
    \P^G(x\leftrightarrow \partial B(x,L)\text{ in }E^{\geq h})\leq C\exp(-cL)\text{ for all }x\in{G}\text{ and }L\geq1,
\end{equation}
where $\{x\leftrightarrow \partial B(x,L)\text{ in }E^{\geq h}\}$ is the event that there exists a continuous path $\pi\subset E^{\geq h}$ between $x$ and $\partial B(x,L).$ In particular, $\tilde{h}_*(\G)<0.$
\end{theorem}

Before proving Theorem~\ref{The:h_*<0}, let us give a few examples of graphs on which it is easy to check that the conditions of Theorem~\ref{The:h_*<0} are fulfilled, and thus $\tilde{h}_*<0.$
\begin{corollary}
\label{Cor:h_*<0}
If $\G$ is either a graph such that \eqref{condsizeball} holds and there exist $c>0$ and $C<\infty$ with $\kappa_x\geq c$ and $\lambda_x\leq C$ for all $x\in{G},$ or if $\G$ is a $(d+1)$-regular tree endowed with unit weights and constant killing measure $\kappa\equiv \hat{\kappa}\in{[C_0,\infty)},$ where $C_0\in{(0,\infty)}$ is a fixed constant depending on $d,$ then \eqref{connectingtoballdecayexpo} holds and $\tilde{h}_*(\G)<0.$
\end{corollary}
\begin{proof}
Let us first assume that \eqref{condsizeball} holds and there exist $c>0$ and $C<\infty$ with $\kappa_x\geq c$ and $\lambda_x\leq C$ for all $x\in{G}.$ By Theorem~\ref{The:h_*<0}, one only needs to prove that \eqref{condgvscap} holds. Under these conditions, the probability that the discrete random walk $\hat{Z}$ on $G$ is killed at time one is uniformly bounded from below by $\inf_{x\in{G}}\kappa_x/\lambda_x.$ Therefore, if the graph distance between $x$ and $y$ is $L,$ we have that $P_x(H_y<\infty)\leq (1-\inf_{x\in{G}}\kappa_x/\lambda_x)^L,$ and that the number of times a random walk starting in $y$ return in $y$ is smaller than a geometric random variable with parameter $\kappa_y/\lambda_y.$ Therefore  \eqref{condgvscap} holds since
\begin{equation}
\label{boundonGL}
    g(L)\leq \sup_{x\in{G}}\frac{1}{\kappa_x}\exp\Big(\log\Big(1-\inf_{x\in{G}}\frac{\kappa_x}{\lambda_x}\Big)L\Big).
\end{equation}

Let us now assume that $\G$ is a $(d+1)$-regular tree endowed with unit weights and constant killing measure $\kappa\equiv \hat{\kappa}\in{[C_0,\infty)}.$ It is clear that \eqref{condsizeballtree} holds for $c_1=\log(d+1)$ and that $\lambda_x=(d+1)+\hat{\kappa}$ for all $x\in{G},$ and thus in view of Theorem~\ref{The:h_*<0}, we only need to prove that \eqref{condgvscap} holds for some $c_0$ such that $\log(d+1)\leq c_0\cdot c_2.$ Let us call $p=p(\hat{\kappa},d)=P_x(H_y<\zeta)$ for some neighbours $x$ and $y,$ which does not depend on the choice of these neighbours by transitivity. Since for all $x,y\in{G},$ the random walk starting in $x$ visits $y$ if and only if it visits the first vertex on the geodesic between $x$ and $y,$ then the second and so on, one can easily show by the Markov property that $P_x(H_y<\zeta)=p^{d(x,y)}.$ By transitivity, we immediately obtain that \eqref{condgvscap} holds for $c_0=\log(1/p).$ But since $p(\hat{\kappa},d)\rightarrow0$ when $\hat{\kappa}\rightarrow\infty,$ one can choose the constant $C_0$ large enough so that $\log(d+1)\leq \log(1/p(\hat{\kappa},d))\cdot c_2$ if $\hat{\kappa}\geq C_0,$ and we can conclude.
\end{proof}

Corollary~\ref{Cor:h_*<0} indicates that when the weights and killing measure are uniformly bounded away from $0$ and $\infty,$ we typically have $\tilde{h}_*(\G)<0.$ This holds for almost any graph with sub-exponential volume growth, for instance for the typical example of the massive $d$ dimensional lattice $\Z^d,$ $d\geq 3,$ or for $(d+1)$-regular trees, which have exponential volume growth for $d\geq2,$ when the killing measure is large enough. Note that on the contrary $\tilde{h}_*(\G)\geq0$ in the massless case $\kappa\equiv0$ by \eqref{eq:h0<1thenh*>0}.

We now turn to the proof of Theorem~\ref{The:h_*<0}, and we are first going to show that the probability in \eqref{connectingtoballdecayexpo} can be made arbitrarily small for an adapted choice of $h,$ depending on $L,$ which will follow from the coupling between random interlacements and the Gaussian free field \eqref{levelsetsvsIu} and a result from \cite{MR3502602} about the probability to connect two vertices in the level sets at level $0.$ Under the assumption \eqref{condgvscap}, to simplify notation we define
\begin{equation}
    \label{defc2}
    c_3=\frac{c_0}{2},
\end{equation}
for all $L>0$
\begin{equation}
\label{defpL}
    p(L)=\left\{\begin{array}{ll}
         \sup_{x\in{G}}\mathrm{cap}(B(x,L))&\text{if }\eqref{condsizeball}\text{ holds,}\\
          \sup_{x\in{G},y\in{\partial B(x,L)}}\mathrm{cap}([x,y])&\text{if }\G\text{ is a tree,}
    \end{array}\right.
\end{equation}
where $[x,y]\subset\tilde{\G}$ denotes the geodesic path between $x$ and $y$ in the cable system, and for all $h>0$
\begin{equation}
    a(L,h)=\left\{\begin{array}{ll}
         \sup_{x\in{G}}\P^G\left(x\leftrightarrow \partial B(x,L)\text{ in }E^{\geq -h}\right)&\text{if }\eqref{condsizeball}\text{ holds,}\\
          \sup_{x\in{G},y\in{\partial B(x,L)}}\P^G\left(x\leftrightarrow y\text{ in }E^{\geq -h}\right)&\text{if }\G\text{ is a tree.}
          \end{array}\right.
\end{equation}
Note that there are trees such that \eqref{condsizeball} holds, and one can then choose for instance the first definition for $p(L)$ and $a(L,h).$

\begin{lemma}
\label{iniren}
Let $\G$ be a graph such that $\lambda_x\leq C$ for all $x\in{G},$ \eqref{condgvscap} holds and either \eqref{condsizeball} holds or $\G$ is a tree, then there exists $C<\infty$ such that
\begin{equation}
\label{eq:iniren}
    a\left(L,\frac{t}{\sqrt{p(L)}}\right)\leq t^2+e^{-c_3L}\text{ for all }L\geq C\text{ and }t\geq0.
\end{equation}
\end{lemma}
\begin{proof}
First note that 
\begin{equation}
    \label{Greenimpliescap}
    \text{ if \eqref{condgvscap} is satisfied, then the condition \eqref{capcondition} holds.}
\end{equation}
Indeed, let $A\subset G$ be a finite connected set with diameter $n+1,$ and for each $i\in\{0,\dots,n\}$ let $x_i\in{A}$ be such that $d(x_0,x_i)=i.$ Then by \cite[\eqref*{4variational}]{DrePreRod3} and \eqref{condgvscap} we have 
\begin{equation*}
    \mathrm{cap}(A)\geq\left(\frac{1}{(n+1)^2}\sum_{i,j=0}^ng(x_i,x_j)\right)^{-1}\geq\left(\frac{2}{n+1}\sum_{k=0}^{n}C\exp(-ck)\right)^{-1}\geq cn,
\end{equation*}
and \eqref{Greenimpliescap} then follows from \cite[\eqref*{4capconditiondis}]{DrePreRod3}. 

Let us take $u=t^2/(2p(L)),$ and first assume that \eqref{condsizeball} holds. Recalling the definition of $\mathcal{C}_u$ from \eqref{levelsetsvsIu}, if $x\leftrightarrow \partial B(x,L)$ in $\mathcal{C}_u\cup \{x\in{\tilde{\G}}:\,\phi_x>0\},$ then either $\I^u\cap B(x,L)\neq\varnothing$ or $x\leftrightarrow \partial B(x,L)$ in $\{x\in{\tilde{\G}}:\,|\phi_x|>0\}.$ By symmetry of the Gaussian free field, \eqref{Greenimpliescap}, \eqref{eq:capimplies0bounded} and \eqref{levelsetsvsIu}, we obtain
\begin{equation*}
    \P^G\Big(x\leftrightarrow \partial B(x,L)\text{ in }E^{\geq -\frac{t}{\sqrt{p(L)}}}\Big)\leq \P^I(\I^u\cap B(x,L)\neq\varnothing)+2\P^G(x\leftrightarrow\partial B(x,L)\text{ in }E^{\geq0}).
\end{equation*}
By \eqref{defIu} and \eqref{defgb} we moreover have
\begin{equation*}
    \P^I(\I^u\cap B(x,L)\neq\varnothing)=1-\exp\big(-u\mathrm{cap}(B(x,L))\big)\leq 1-\exp\left(-\frac{t^2}{2}\right)\leq t^2.
\end{equation*}
Moreover, by \cite[Propositions~2.1 and 5.2]{MR3502602}, since $g(y,y)\geq \lambda_y^{-1}\geq c$ for all $y\in{G},$ we have by a union bound that for all $L$ large enough
\begin{equation*}
    \P^G(x\leftrightarrow\partial B(x,L)\text{ in }E^{\geq0})\leq Cb(L)\arcsin(Cg(L))\leq \exp(-c_3L),
\end{equation*}
where we used \eqref{condgvscap} and \eqref{condsizeball} in the last inequality, as well as the fact that $\arcsin(x)\leq Cx$ for all $x\leq1.$

Let us now assume that $\G$ is a tree and fix some $x\in{G},$ $L\geq1,$ $y\in{\partial B(x,L)}$ and $t\geq0.$ We can prove similarly as before that since $\G$ is a tree 
\begin{equation}
\label{eq:whytrees}
\begin{split}
    \P^G\Big(x\leftrightarrow y\text{ in }E^{\geq -\frac{t}{\sqrt{p(L)}}}\Big)&=\P^G\Big([x,y]\subset E^{\geq -\frac{t}{\sqrt{p(L)}}}\Big)
\\&\leq \P^I(\I^u\cap [x,y]\neq\varnothing)+2\P^G( x\leftrightarrow y\text{ in } E^{\geq0}).
\end{split}
\end{equation}
Using \eqref{defIu} and \cite[Propositions~2.1 and 5.2]{MR3502602}, we can conclude.
\end{proof}

Note that in \eqref{eq:whytrees}, the first equality only holds when $\G$ is a tree, which is the reason why the assumption \eqref{condsizeball} of sub-exponential volume growth can be replaced by the assumption \eqref{condsizeballtree} of exponential volume growth in the case of trees. Lemma~\ref{iniren} implies that the probability in \eqref{connectingtoballdecayexpo} can be made arbitrarily small by taking $h=\frac{t}{\sqrt{p(L)}},$ $t$ small enough and $L$ large enough. This will serve as the base of a renormalization scheme, similar to the one presented in \cite[Section~7]{MR3420516}, that we now explain. For some $L_0>0$ we define recursively
\begin{equation}
\label{defLk}
    L_{k+1}=2L_k\Big(1+\frac{1}{(k+1)^{1+\alpha}}\Big)\text{ for all }k\geq0,
\end{equation}
where $\alpha$ is the same constant as in \eqref{condsizeball} if \eqref{condsizeball} holds, and $\alpha=1$ otherwise. Then there exists a constant $C_1<\infty$ only depending on $\alpha$ such that
\begin{equation}
\label{ren:boundonL_k}
    2^kL_0\leq L_k\leq C_12^kL_0\text{ for all }k\in\N.
\end{equation}
Let us also define for all $t\geq0$ and $k\in\N_0$
\begin{equation}
    h_k(t)=\frac{t}{\sqrt{p(L_0)}}\left(1-\sum_{i=1}^{k}\frac{c}{i^{1+\alpha}}\right),
\end{equation}
where the constant $c=c(\alpha)$ is chosen small enough so that $h_{\infty}(t)>0$ for all $t>0,$ where $h_{\infty}(t)$ is the limit of $h_k(t)$  as $t\rightarrow\infty.$ For any $h\in\R$ we have
\begin{equation}
\label{Lk+1inclusLk}
    \{x\leftrightarrow \partial B(x,L_{k+1})\text{ in }E^{\geq h}\}\subset\hspace{-4mm}\bigcup_{y\in{\partial B(x,L_{k+1})}}\hspace{-3mm}\{y\leftrightarrow \partial B(y,L_k)\text{ in }E^{\geq h}\}\cap\{x\leftrightarrow \partial B(x,L_k)\text{ in }E^{\geq h}\},
\end{equation}
and if $\G$ is a tree, for any $y\in{\partial B(x,L_{k+1})},$ letting $z$ and $z'$ be the unique points on the geodesic between $x$ and $y$ such that $z\in{\partial B(x,L_k)}$ and $z'\in{\partial B(y,L_k)},$ 
\begin{equation}
\label{Lk+1inclusLktree}
    \{x\leftrightarrow y\text{ in }E^{\geq h}\}\subset\{y\leftrightarrow z'\text{ in }E^{\geq h}\}\cap\{x\leftrightarrow z\text{ in }E^{\geq h}\}.
\end{equation}
The two events on the right-hand side of \eqref{Lk+1inclusLk} and \eqref{Lk+1inclusLktree} are measurable with respect to the field on distant sets since for all $y\in{\partial{B(x,L_{k+1})}}$ 
\begin{equation}
\label{xfarfromy}
    d\big(B(y,L_k),B(x,L_k)\big)\geq \frac{L_k}{(k+1)^{1+\alpha}},
\end{equation}
upon choosing $L_0$ large enough.

Let us now recall the decoupling inequality from \cite{MR3325312}, which implies that the events on the right-hand side of \eqref{Lk+1inclusLk} and \eqref{Lk+1inclusLktree}  are almost uncorrelated, up to some sprinkling parameter.

\begin{lemma}
Let $\G$ be a transient graph. For all $L\geq1,$ $\delta>0,$ $x_1,x_2\in{G}$ with $s\stackrel{\mathrm{def.}}{=}d(B(x_1,L),B(x_2,L))>0,$ increasing events $A_1$ and $A_2$ such that $A_i\subset C(\tilde{\G},\R)$ is measurable with respect to the $\sigma$-algebra generated by the coordinate functions and depend only on the value of the function on the closure of the edges $I_{\{x,y\}},$ $x,y\in{B(x_i,L)}$ for each $i\in{\{1,2\}},$ 
\begin{equation}
    \label{decoupling}
    \begin{split}
    \P^G(\phi\in{A_1\cap A_2})&\leq\P^G(\phi+\delta\in{A_1})\P^G(\phi+\delta\in{A_2})
    +2b(L)\exp\left(-\frac{\delta^2}{8g(s)}\right).
    \end{split}
\end{equation}
\end{lemma}
\begin{proof}
As noted in \cite[Theorem~6.2]{DrePreRod2}, one can easily adapt the proof of \cite[Corollary~1.3]{MR3325312} to the cable system, and we obtain 
\begin{equation*}
    \P^G(\phi\in{A_1\cap A_2})\leq\P^G(\phi+\delta\in{A_1})\P^G(\phi+\delta\in{A_2})
    +2\P^G\left(H_{\delta}^c\right),
\end{equation*}
where, denoting by $h$ the harmonic average of the field $\phi$ in $B(x_1,L),$
\begin{equation*}
H_{\delta}=\left\{\sup_{x,y\in{B(x_2,L)}}\sup_{z\in{\overline{I_{\{x,y\}}}}}|h_z|\leq\frac\delta2\right\}
\end{equation*}
Note that $h_z$ is linear on the closure of $I_{\{x,y\}},$ and so $H_{\delta}$ is actually equal to the event $G_{\delta}$ from \cite[(1.6)]{MR3325312} with $K_i=B(x_i,L).$ 
We can now conclude by \cite[Proposition~1.4]{MR3325312}, whose proof can easily be adapted to any transient graph.
\end{proof}

Combining \eqref{Lk+1inclusLk}, \eqref{xfarfromy} and \eqref{decoupling} with Lemma~\ref{iniren}, we can now derive a bound on $a(L_k,h_k(t)).$

\begin{lemma}
\label{lemmainduction}
Let $\G$ be a graph such that $\lambda_x\leq C$ for all $x\in{G}$ and \eqref{condgvscap} holds. If \eqref{condsizeball} holds, then there exist a constant $C_2<\infty$ depending only on $\alpha,$ as well as $C<\infty$ such that for all $t\in{(0,1/2]}$ and $L_0\geq C$ with 
\begin{equation}
\label{L0mustbelarge}
    t^2\geq CL_0\exp\left(-\frac{c_3L_0}{C_2}\right),
\end{equation}
and for all $k\in\N_0,$ we have
\begin{equation}
\label{toproverecursively}
    a\big(L_k,h_k(t)\big)\leq \frac12\exp\left(C2^kL_0\sum_{i=0}^k\frac{1}{\log(L_i)^{1+\alpha}}\right)(t^2+e^{-c_3L_0})^{2^k}.
\end{equation}
If $\G$ is a tree such that \eqref{condsizeballtree} holds, then there exists $C<\infty$ such that for all $t\in{(0,1/2]}$ and $L_0\geq C$ satisfying \eqref{L0mustbelarge}, and for all $k\in\N_0,$ we have
\begin{equation}
\label{toproverecursivelytree}
    a\big(L_k,h_k(t)\big)\leq \frac12(2(t^2+e^{-c_3L_0}))^{2^k}.
\end{equation}
\end{lemma}
\begin{proof}
We fix some $t\in{(0,1/2]}$ and first prove \eqref{toproverecursively} by induction on $k$ under the assumptions \eqref{L0mustbelarge} and \eqref{condsizeball}. The statement for $k=0$ follows directly from \eqref{eq:iniren} if $L_0$ is large enough. Let us now assume that \eqref{toproverecursively} holds for some $k\in{\N_0}.$ Combining \eqref{Lk+1inclusLk}, \eqref{xfarfromy} and \eqref{decoupling} we obtain by a union bound that
\begin{equation}
\label{eq:recursion1}
\begin{split}
    a&\big(L_{k+1},h_{k+1}(t)\big)
    \\&\leq b(L_{k+1})\left(a\big(L_k,h_k(t)\big)^2+2b(L_k)\exp\left(-\Big(\frac{ct}{\sqrt{p(L_0)}(k+1)^{1+\alpha}}\Big)^2\frac{1}{8g(s_k)}\right)\right),
  \end{split}
\end{equation}
where $s_k=\frac{L_k}{(k+1)^{1+\alpha}}.$ By \eqref{toproverecursively}, we moreover have that
\begin{equation}
\label{eq:recursion2}
\begin{split}
    b(L_{k+1})a\big(L_k,h_k(t)\big)^2&\leq\frac14b(L_{k+1})\exp\left(C2^{k+1}L_0\sum_{i=0}^k\frac{1}{\log(L_i)^{1+\alpha}}\right)(t^2+e^{-c_3L_0})^{2^{k+1}}\\&\leq\frac14\exp\left(C2^{k+1}L_0\sum_{i=0}^{k+1}\frac{1}{\log(L_i)^{1+\alpha}}\right)(t^2+e^{-c_3L_0})^{2^{k+1}},
\end{split}
\end{equation}
where the last inequality holds by \eqref{condsizeball} and \eqref{ren:boundonL_k} when choosing the constant $C$ large enough, independently of $t$ and $k.$ Moreover, by \eqref{condgvscap}, \eqref{condsizeball} and \eqref{ren:boundonL_k}, noting that 
\begin{equation*}
    p(L_0)\leq C\sum_{k=0}^{L_0}b(k)\leq CL_0\exp\Big(\frac{cL_0}{\log(L_0)^{1+\alpha}}\Big)\leq g(s_k)^{-1/2}
\end{equation*}
for $L_0$ large enough, independently of $k,$ we have
\begin{align*}
   \left(\frac{ct}{\sqrt{p(L_0)}(k+1)^{1+\alpha}}\right)^2\frac{1}{8g(s_k)}\geq \frac{ct^2g(s_k)^{-\frac12}}{(k+1)^{2(1+\alpha)}}&\geq \frac{ct^2}{(k+1)^{2(1+\alpha)}}\exp\left(\frac{c_02^kL_0}{2(k+1)^{1+\alpha}}\right)
    \\&\geq ct^22^{k+1}\exp(c_3L_0/C_2),
\end{align*}
for $L_0$ large enough, where $C_2$ is a constant depending only on $\alpha.$ Therefore, upon choosing $C$ large enough, independently of $k,$ if \eqref{L0mustbelarge} holds then by \eqref{condsizeball} and \eqref{ren:boundonL_k}
\begin{equation}
\label{eq:recursion3}
    2b(L_k)\exp\left(-\Big(\frac{ct}{\sqrt{p(L_0)}(k+1)^{1+\alpha}}\Big)^2\frac{1}{8g(s_k)}\right)\leq \frac14e^{-c_3L_02^{k+1}}\leq \frac14(t^2+e^{-c_3L_0})^{2^{k+1}},
\end{equation}
and by \eqref{condsizeball} and \eqref{ren:boundonL_k}, upon choosing $C$ large enough, 
\begin{equation}
\label{eq:recursion4}
   b(L_{k+1})\leq\exp\left(C2^{k+1}L_0\sum_{i=0}^{k+1}\frac{1}{\log(L_i)^{1+\alpha}}\right).
\end{equation}
We can easily conclude that \eqref{toproverecursively} holds for $k+1$ by combining \eqref{eq:recursion1}, \eqref{eq:recursion2}, \eqref{eq:recursion3} and \eqref{eq:recursion4}.

The proof of \eqref{toproverecursivelytree} is similar when $\G$ is a tree and \eqref{condsizeballtree} holds. Indeed, when $k=0,$ \eqref{toproverecursivelytree} is just \eqref{eq:iniren}. Let us now assume that \eqref{toproverecursivelytree} holds for some $k\in\N_0.$ Combining \eqref{Lk+1inclusLktree}, \eqref{xfarfromy} and \eqref{decoupling} we have
\begin{equation}
\label{eq:recursion1tree}
    a\big(L_{k+1},h_{k+1}(t)\big)\leq a\big(L_k,h_k(t)\big)^2+2b(L_k)\exp\left(-\left(\frac{ct}{\sqrt{p(L_0)}(k+1)^{1+\alpha}}\right)^2\frac{1}{8g(s_k)}\right).
\end{equation}
Moreover, one can easily prove that \eqref{eq:recursion3} still holds under condition \eqref{condsizeballtree} since $p(L_0)\leq CL_0\leq g(s_k)^{-1/2}$ if $L_0$ is large enough, and we can easily deduce from \eqref{eq:recursion1tree} that \eqref{toproverecursivelytree} also holds for $k+1.$
\end{proof}

In view of Lemma~\ref{lemmainduction}, we are now ready to prove Theorem~\ref{The:h_*<0}.

\begin{proof}[Proof of Theorem~\ref{The:h_*<0}]
We take $t\equiv t(L_0)=e^{-c_3L_0/(4C_2)},$ then \eqref{L0mustbelarge} holds if $L_0$ is large enough. Let us first assume that \eqref{condsizeball} holds. By Lemma~\ref{lemmainduction} and \eqref{ren:boundonL_k} one can fix $L_0$ large enough, so that for all $k\in\N_0$
\begin{align*}
    a\big(L_k,h_k(t)\big)&\leq\frac12\exp\left(C2^kL_0\sum_{i=0}^{\infty}\frac{1}{\big(i+\log(L_0)\big)^{1+\alpha}}\right)\left(2\exp(-cL_0)\right)^{2^k}
    \\&\leq2^{2^k}\exp\left(C2^k\frac{L_0}{\log(L_0)^{\alpha}}\right)\exp\big(-c2^{k}L_0\big)
    \\&\leq\exp\big(-c2^{k}L_0\big)\leq\exp(-cL_k).
\end{align*}
For any $L\geq1,$ let us take $k\in\N_0$ such that $L_{k}\leq L\leq L_{k+1},$ then for all $x\in{G}$ by \eqref{ren:boundonL_k}
\begin{equation*}
     \P^G\big(x\leftrightarrow \partial B(x,L)\text{ in }E^{\geq -h_{\infty}(t)}\big)\leq a\big(L_k,h_k(t)\big)\leq\exp(-cL_k)\leq\exp(-cL).
\end{equation*}
Taking the limit as $L\rightarrow\infty,$ we get that the component of $x$ in $E^{\geq -h_{\infty}(t)}$ is $\P^G$-a.s.\ bounded for all $x\in{G},$ and by a union bound $E^{\geq -h_{\infty}(t)}$ contains $\P^G$-a.s.\ only bounded components, that is $\tilde{h}_*(\G)\leq -h_{\infty}(t)<0.$

Let us now assume that $\G$ is a tree such that \eqref{condsizeballtree} holds. We then have by a union bound that
\begin{align*}
    \P^G\big(x\leftrightarrow \partial B(x,L_k)\text{ in }E^{\geq -h_{k}(t)}\big)&\leq b(L_k)a\big(L_k,h_k(t)\big)
    \\&\leq C\exp(c_1C_1L_02^k)4^{2^k}\exp\Big(-\frac{c_3}{2C_2\vee1}L_02^k\Big),
\end{align*}
where we used \eqref{condsizeballtree}, \eqref{ren:boundonL_k} and \eqref{toproverecursivelytree} in the last inequality. In view of \eqref{defc2}, since $C_1$ and $C_2$ depend only on $\alpha$ and we can choose $\alpha=1$ when $\G$ is a tree, we can define the absolute constant 
\begin{equation*}
    c_2\stackrel{\text{def.}}{=}\frac{1}{4C_1(2C_2\vee1)},
\end{equation*}
and, if $c_1\leq c_0\cdot c_2,$ then, taking $L_0$ large enough, we can easily conclude that \eqref{connectingtoballdecayexpo} holds and $\tilde{h}_*(\G)<0$ similarly as before.
\end{proof}

\begin{remark}
\label{rkboundedlambda}
    Consider any graphs satisfying \eqref{condsizeball} and the conditions of example~(c) below  Corollary~\ref{Cor:h_0=1andh_*=0}, then \eqref{condgvscap} holds by \eqref{boundonGL}. Similarly, if $\G$ is the $(d+1)$-regular tree, $d\geq2,$ with $\kappa_n= t\lambda_n$ for some $t>0$ and sequence $\lambda_n$ increasing to infinity, which is a graph as in example~(b) below  Corollary~\ref{Cor:h_0=1andh_*=0}, then it satisfies \eqref{condgvscap} with $c_0=\log(1+t/(d+1))$ by \eqref{boundonGL}, as well as \eqref{condsizeballtree} with $c_1=\log(d+1),$ and in particular $c_1< c_0\cdot c_2$ if $t$ is chosen large enough. By Corollary~\ref{Cor:h_0=1andh_*=0}, we have $\tilde{h}_*=0$ on these graphs, and so the condition $\lambda_x\leq C$ in Theorem~\ref{The:h_*<0} is necessary. It is however not clear whether one could replace the condition \eqref{condgvscap} by $\h_{\text{kill}}\equiv1,$ which is necessary in view of \eqref{eq:h0<1thenh*>0}, and allow graphs with exponential growth instead of either assuming \eqref{condsizeball} or that the graph is a tree with fast enough decay of the Green function.
\end{remark}

One can easily derive from Theorem~\ref{The:h_*<0} and \eqref{levelsetsvsIu} a similar result but for the random interlacement set $\I^u$ on the cable system.

\begin{corollary}
\label{cor:u_*>0}
Let $\G$ be any graph satisfying the assumptions of either Theorem~\ref{The:h_*<0} or Corollary~\ref{Cor:h_*<0}. There exist $u>0$ and $c>0$ such that
\begin{equation}
\label{connectingtoballdecayexpoIu}
    \P^I(x\leftrightarrow \partial B(x,L)\text{ in }\I^{u})\leq\exp(-cL)\text{ for all }x\in{G}\text{ and }L\geq1,
\end{equation}
and in particular the critical parameter associated with the percolation of $\I^u$ on the cable system is positive. 
\end{corollary}
\begin{proof}
It follows from \eqref{Greenimpliescap}, \eqref{eq:capimplies0bounded} and \eqref{levelsetsvsIu} that $\I^u$ is stochastically dominated by $E^{\geq -\sqrt{2u}}.$ The inequality \eqref{connectingtoballdecayexpoIu} then follows from \eqref{connectingtoballdecayexpo} for $u=h^2/2.$ Moreover taking the limit as $L\rightarrow\infty,$ we obtain that $\I^u$ contains $\P^I$-a.s.\ only bounded components for any $u$ satisfying \eqref{connectingtoballdecayexpoIu}, and we can conclude. 
\end{proof}

\begin{remark}
\begin{enumerate}[1)]
\label{linkwithfinitary}
\item In \cite{cai2021rigorous}, a result similar to Corollary~\ref{cor:u_*>0} is proven. They consider finitary random interlacements, which by \cite[Proposition~4.1]{MR3962876}, see also \eqref{disinterh0=1}, correspond to random interlacements on the graph $\G^T=(\Z^d,\lambda^T,\kappa^T),$ $d\geq3,$ where $\lambda^T_{x,y}=\frac{T}{T+1}$ and $\kappa^T_x=\frac{2d}{T+1}.$ If $\I^u(\subset\tilde{\G})$ contains only bounded components, then the set of edges crossed by a trajectory in the random interlacement process contains only finite components, which happens a.s.\ for $u$ small enough by Corollary~\ref{cor:u_*>0}. This is thus an alternative proof of the bound $u_c>0$ from \cite[Theorem~4]{cai2021rigorous}.
\item \label{h_*u_*finite}One can easily show that $\tilde{h}_*>-\infty$ and that the critical parameter associated with the percolation of $\I^u,$ $u>0,$ as a subset of the cable system or as a set of edges crossed by trajectories in the random interlacement process, is finite on any graph $\G$ with uniformly bounded weights and killing measure, that is $c\leq \kappa_x\leq C$ and $c\leq \lambda_{x,y}\leq C,$ and such that $p_c<1,$ where $p_c$ is the critical parameter for Bernoulli bond percolation on $\G.$ We refer to \cite{DcGoRaSeYa} for a review of the literature and a proof of the inequality $p_c<1$ under rather general conditions, which incidentally uses the property \eqref{eq:h0<1thenh*>0}. Finitary random interlacements on $\Z^d$ clearly fulfill all these hypotheses, and combining with the previous remark, we obtain  \cite[Theorem~4]{cai2021rigorous}, but on a more general class of graphs.

Indeed, let us start for each $x\in{G}$ a poissonian number of independent trajectories starting in $x$ on $\tilde{\G}^{\Ed}$ with parameter $u\kappa_x.$ Then by \eqref{disinterh0=1}, $\I^u\cap\tilde{\G}^{\Ed}$ has the same law as the set of points visited by one of these trajectories. For each edge $e=\{x,y\},$ the number of trajectories starting in either $x$ or $y$ and crossing first $e$ has law 
\begin{equation*}
    \text{Poi}\Big(\frac{u\kappa_x\lambda_{x,y}}{\lambda_x}+\frac{u\kappa_y\lambda_{x,y}}{\lambda_y}\Big).
\end{equation*}
Therefore for any $u$ large enough so that $1-\exp\big(-u\lambda_{x,y}(\kappa_x\lambda_x^{-1}+\kappa_y\lambda_y^{-1})\big)>p_c$ for all $x\in{G},$ there is an infinite connected component of edges crossed by the discrete killed random interlacement process, and thus $\I^u(\subset\tilde{\G})$ contains an unbounded connected component with positive probability. One can easily check that \eqref{capcondition} holds for a graph with uniformly bounded weights and killing measure since $e_A(x)\geq \kappa_x\geq c$ for all $x\in{G}$ and $A\subset G.$ In particular, using \eqref{eq:capimplies0bounded} and \eqref{levelsetsvsIu}, we know that $E^{\geq-\sqrt{2u}}$ contains an infinite connected component with positive probability if $\I^u$ does, and we can conclude.
\item There also exist graphs for which $\tilde{h}_*=-\infty,$ for instance finite graphs or $\mathbb{Z}$ with constant weights and constant positive killing measure, as for any $h\in{\mathbb{R}}$ there is a.s.\ infinitely many $n\in{\mathbb{Z}}$ with $\phi_n<h$ by ergodicity. We refer to Remark~\ref{rk:Z20counterexample},\ref{rk:h*-inftyZ20}) for a less trivial example of a graph for which $\tilde{h}_*=-\infty.$
\end{enumerate}
\end{remark}

\section{Doob \texorpdfstring{$\h$}{h}-transform}
\label{sec:Doob}
In this section, we introduce the notion of the Doob $\h$-transform $\G_{\h}$ of a graph $\G,$ when $\h:\tilde{\G}\rightarrow(0,\infty)$ is an harmonic function, so that the diffusion $X$ on the cable system ${\tilde{\G}_\h}$ of $\G_{\h}$ is related to the $\h$-transform of the diffusion $X$ on $\tilde{\G},$ see \eqref{semigrouph}. One can then also relate the law of the Gaussian free field on $\tilde{\G}_{\h}$ to the Gaussian free field on $\tilde{\G},$ see \eqref{eq:GFFh}, from which one can deduce an effective criterion for \eqref{0bounded} to hold in terms of $\h$-transform, see Corollary~\ref{capimplieseverythingh}, which will be useful in Section~\ref{sec:Z20}. Introducing the notion of $\h$-transform of random interlacements, see Definition \ref{defhtransforminter}, we finally use results from \cite{DrePreRod3} to obtain under condition \eqref{0bounded} the law of the $\h$-transform of the capacity of the level sets $E_{\h}^{\geq h}(x_0)$ and an isomorphism between the Gaussian free field and the $\h$-transform of random interlacements for various choices of the harmonic function $\h,$ see Theorem~\ref{couplingintergffh} and Corollaries~\ref{h0transformiso} and \ref{couplingintergffK}.

\begin{definition}
\label{harmo}
We say that a function $\h:\tilde{\G}\rightarrow(0,\infty)$ is harmonic on $\tilde{\G}$ if for all $x\in{G}$
\begin{enumerate}[1)]
    \item $\h(\partial I_x)\stackrel{\text{def.}}{=}\lim\limits_{t\nearrow\rho_x}\h(x+t\cdot I_x)$ exists and is finite,
    \item if $e=\{x,y\}\in{E}$ or $e=x\in{G},$ then $t\mapsto \h(x+t\cdot I_e)\in{C^2([0,\rho_e),(0,\infty))}$ and 
    \begin{equation*}
        \frac{\mathrm{d}^2\h(x+t\cdot I_e)}{\mathrm{d}t^2}=0\text{ for all }t\in{[0,\rho_e)},
    \end{equation*}
    \item and for all $x\in{G}$ 
\begin{equation}
\label{defharmonic}
    \Big(\frac{\mathrm{d}\h(x+t\cdot I_{x})}{\mathrm{d}t}+\sum_{y\sim x}\frac{\mathrm{d}\h(x+t\cdot I_{\{x,y\}})}{\mathrm{d}t}\Big)\Big|_{t=0}=0.
\end{equation}
\end{enumerate} 
We define the $\h$-transform $\G_{\h}$ of the graph $\G$ as the graph with vertex set $G_{\h}=G,$ with weights $\lambda^{(\h)}_{x,y}=\h(x)\h(y)\lambda_{x,y},$ $x,y\in{G},$ and with killing measure $\kappa^{(\h)}_x=\kappa_x\h(x)\h(\partial I_x),$ $x\in{G}.$
\end{definition}

Note that this corresponds to the usual definition of harmonicity, that is the generator associated with the form $\mathcal{E}_{\tilde{\G}}$ from \eqref{Dirichlet} applied to $\h$ is equal to $0,$ see around \cite[(2.1)]{MR3152724} for a description of this generator. In order to prove that the graph $\G_{\h}$ is part of the setting described at the beginning of Section~\ref{sec:intro}, we only need to prove that the diffusion $X$ on $\tilde{\G}_{\h}$ is transient, which follows from \eqref{relationGreenfunctionh}. The conditions 2) and 3) of Definition \ref{harmo} can be respectively equivalently stated as follows: for all $e\in{E\cup G}$ and $t\in{[0,\rho_e)},$ 
\begin{equation}
\label{hformulaonedges}
    \h(x+t\cdot I_e)=\begin{cases}2t\lambda_{x,y}\h(y)+\big(1-2t\lambda_{x,y}\big)\h(x)&\text{if }e=\{x,y\}\in{E}
    \\2t\kappa_x\h(\partial Ix)+\big(1-2t\kappa_x\big)\h(x)&\text{if }e=x\in{G},
    \end{cases}
\end{equation}
and for all $x\in{G}$
\begin{equation}
\label{hformulaoutsideofedge}
    \kappa_x\h(\partial I_x)+\sum_{y\sim x}\lambda_{x,y}\h(y)=\lambda_x\h(x).
\end{equation}
In particular, the total weight of a vertex ${x}\in{{G}_{\h}}$ is $\lambda_{{x}}^{(\h)}=\h(x)^2\lambda_x.$ Note that since $\h>0,$ the edge set $E_{\h}$ of $\G_{\h}$ is equal to $E,$ and we will often identify the edges and vertices of $\G_{\h}$ to the edges and vertices of $\G.$ 

Let us define a function $\psi_\h:\tilde{\G}\rightarrow\tilde{\G}_\h$ such that for all $e=\{x,y\}\in{E}$ or $e=x\in{G}$
\begin{equation}
\label{defpsih}
    \psi_\h(x+t\cdot I_e)\stackrel{\mathrm{def.}}{=}x+\frac{t}{\h(x)\h(x+t\cdot I_e)}\cdot I_{e}\text{ for all }t\in{[0,\rho_e)},
\end{equation}
where with a slight abuse of notation, $I_e\subset\tilde{\G}$ on the left-hand side of \eqref{defpsih} and $I_e\subset\tilde{\G}_{\h}$ on the right-hand side of \eqref{defpsih}. We also take $\psi_{\h}(\Delta)=\Delta.$ Using \eqref{hformulaonedges}, one can easily check that this definition does not depend on the choice of the endpoint $x$ or $y$ of $I_e$ when $e=\{x,y\}\in{E},$ and that $\psi_\h$ is bijective. For any forwards trajectories $w^+\in{W^+_{\tilde{\G}_\h}}$ on $\tilde{\G}_\h,$ we define the time change
\begin{equation}
\label{defthetah}
    \theta_{\h}^{w^+}(t)\stackrel{\mathrm{def.}}{=}\inf\Big\{s\geq0:\,\int_0^s\h\big(\psi_\h^{-1}(w^+(u))\big)^{4}\,\mathrm{d}u>t\Big\}\text{ for all }t\in{[0,\infty)},
\end{equation}
with the conventions $\h(\Delta)=0$ and $\inf\varnothing=\zeta,$ and
\begin{equation}
\label{defxi}
    (\xi_\h(w^+))(t)\stackrel{\mathrm{def.}}{=}\psi_\h^{-1}\big(w^+({\theta_{\h}^{w^+}(t)})\big)\text{ for all }t\in{[0,\infty)}.
\end{equation}
The process $\xi_{\h}(X)$ is thus a stochastic process on $\tilde{\G}$ under $P^{\tilde{\G}_{\h}}_{\psi_{\h}(x)},$ $x\in{\tilde{\G}},$ and we call it the $\h$-transform of $X.$ Indeed, we prove in \eqref{semigrouph} that if $T_t^{\tilde{\G}},$ $t\geq0,$ denotes the semigroup on $L^2(\tilde{\G},m)$ associated with $X$ under $P_{x}^{\tilde{\G}},$ $x\in{\tilde{\G}},$ one can relate the semigroup associated with $\xi_{\h}(X)$ to $T_t^{\tilde{\G}},$ in a way which corresponds to the usual definition of the $\h$-transform, see for instance \cite[Chapter~11]{MR2152573}. Moreover, one can also relate the local times associated with $\xi_{\h}(X)$ to the local times associated with $X$ under $P_{\cdot}^{\tilde{\G}_{\h}},$ and the Gaussian free field on $\tilde{\G}_{\h}$ and on $\tilde{\G}.$ 

\begin{proposition}
\label{corh}
Let $\h$ be an harmonic function on $\tilde{\G}.$ Under $P^{\tilde{\G}_{\h}}_{\psi_{\h}(x)},$ $x\in{\tilde{\G}},$
\begin{equation}
    \label{semigrouph}
    \text{the semigroup on $L^2(\tilde{\G},\h^2\cdot m)$ associated with $\xi_{\h}(X)$ is $f\mapsto\frac{1}{\h}T_t^{\tilde{\G}}(f\h),$}
\end{equation}
with respect to the measure $m,$
\begin{equation}
    \label{eq:localtimesh}
    \text{the field of local times associated with $\xi_{\h}(X)$ is } \big(\h(x)^2\ell_{\psi_{\h}(x)}(\theta_{\h}^{X}(t))\big)_{t\geq0,x\in{\tilde{\G}}},
\end{equation}
and the Gaussian field
\begin{equation}
\label{eq:GFFh}
  \big(\h(x)\phi_{\psi_{\h}(x)}\big)_{x\in{\tilde{\G}}}\text{ has the same law under }\P_{\tilde{\G}_{\h}}^G\text{ as }(\phi_x)_{x\in{\tilde{\G}}}\text{ under }\P_{\tilde{\G}}^G.
\end{equation}
\end{proposition}

Similar links between the graph $\G_{\h}$ and the $\h$-transform have already been noticed for the discrete graph in specific contexts, see the proof of \cite[Proposition~4.6]{LuSaTa} or the Appendix of \cite{MR4091511}, and Proposition~\ref{corh} can be seen as a generalization of these results, and its proof is presented in the Appendix.

In view of \eqref{eq:GFFh}, one can transfer the results \eqref{eq:h0<1thenh*>0} and \eqref{eq:capimplies0bounded} about the level sets for the Gaussian free field on $\tilde{\G}_{\h}$, defined as in \eqref{deflevelsets}, to similar results for the sets
\begin{equation}
    \label{defEf}
        E^{\geq h}_{\h}=\{x\in{\tilde{\G}}:\phi_x\geq h\times\h(x)\}\text{ and }E^{\geq h}_{\h}(x_0)=\text{ the connected component of }x_0\text{ in }E^{\geq h}_{\h}.
\end{equation}
It then follows directly from \eqref{eq:GFFh} that for all $h\in\R$ and $x_0\in{\tilde{\G}}$
\begin{equation}
    \label{killedlevelsetsvsh0transform}
     E^{\geq h}(\psi_{\h}(x_0))\text{ has the same law under }\P_{\tilde{\G}_{\h}}^G\text{ as }\psi_{\h}\big(E_{\h}^{\geq h}(x_0)\big)\text{ under }\P_{\tilde{\G}}^G,
\end{equation}
where $E^{\geq h}(x_0)=E^{\geq h}_\textbf{1}(x_0)$ is the connected component of $x_0$ in $E^{\geq h}.$

\begin{corollary}
\label{capimplieseverythingh}
\begin{enumerate}[1)]
    \item If there exists an harmonic function $\h$ on $\tilde{\G}$ such that \eqref{capcondition} is satisfied for $\G_{\h},$
then \eqref{0bounded} holds.
    \item If $\h$ is an harmonic function on $\tilde{\G}$ such that $\h_{\text{kill}}<1$ on $\G_{\h},$ then $E^{\geq h}_{\h}$ contains an unbounded connected component with $\P^G_{\tilde{\G}}$ positive probability for all $h<0.$
\end{enumerate}
\end{corollary}
\begin{proof}
Noting that $E_{\h}^{\geq 0}=E^{\geq 0}$ for any harmonic function $\h$ on $\tilde{\G},$ 1) follows directly from \eqref{eq:capimplies0bounded} and \eqref{killedlevelsetsvsh0transform}, whereas 2) follows directly from \eqref{eq:h0<1thenh*>0} and \eqref{killedlevelsetsvsh0transform}.
\end{proof}

In the rest of this section, we will deduce from Proposition~\ref{corh} and \cite[Theorem~\ref*{4T:main},2)]{DrePreRod3} an explicit formula for the law of the level sets $E^{\geq h}_{\h}$ from \eqref{defEf}, as well as various isomorphisms between the Gaussian free field and either the $\h$-transform of random interlacements, see Theorem~\ref{couplingintergffh}, or killed or surviving random interlacements, see Corollary~\ref{h0transformiso}, or the trajectories in the random interlacement process never hitting a compact $K,$ see Corollary~\ref{couplingintergffK}, which will lead to a proof of Theorem~\ref{couplingintergffdim2} in Section~\ref{sec:Z20}. These results are interesting and widely generalize similar theorems in dimension 2, see \cite[Theorems~5.3 and 5.5]{MR3936156}, or on finite graphs, see \cite[Proposition~2.4]{MR4091511}. They are however not needed to prove our main result Theorem~\ref{mainth}, and the hastened reader can directly skip to Section~\ref{sec:Z20}, see in particular Proposition~\ref{Z20counterexample} therein, to understand how to to use Corollary~\ref{capimplieseverythingh} to find a graph as in Theorem~\ref{mainth},3). 

Let us now define the $\h$-transform of random interlacements. If $w^*\in{W_{\tilde{\G}_\h}^*},$ we denote by $\xi_{\h}^*(w^*)$ the trajectory in $W_{\tilde{\G}}^*$ which corresponds to taking the image modulo time-shift of a trajectory with backwards part $\xi_{\h}((w(-t))_{t\geq0})$ and forwards part $\xi_{\h}((w(t)_{t\geq0}),$ for some $w\in{(p_{\tilde{\G}_\h}^*)^{-1}(w^*)},$ and one can easily check that this definition does not depend on the choice of $w.$ To simplify notation, we also denote by $\xi_{\h}^*:W_{\tilde{\G}_{\h}}^*\times[0,\infty)\rightarrow W_{\tilde{\G}}^*\times[0,\infty)$ the application which associates $(\xi_{\h}^*(w^*),u)$ to $(w^*,u).$

\begin{definition}
\label{defhtransforminter}
If $\h$ is an harmonic function on $\tilde{\G},$ under $\P_{\tilde{\G}_{\h}}^I,$ let us define $\omega^{\h}=\omega\circ(\xi_{\h}^*)^{-1}$ the $\h$-transform of the random interlacement process, and for all $u>0$ we denote by $(\ell^{\h}_{x,u})_{x\in{\tilde{\G}}}$ the family of local times with respect to $m$ associated with $\omega_u^{\h},$ the point process of trajectories in $\omega$ with label at most $u,$ and by $\I^u_{\h}=\{x\in{\tilde{\G}}:\ell_{x,u}^{\h}>0\}$ the $\h$-transform of the interlacement set. 
\end{definition}

In \cite{DrePreRod3}, an explicit formula for the law of the capacity of level sets of the Gaussian free field was given, as well as a signed version of the isomorphism between the Gaussian free field and local times of random interlacements on the cable system under the condition \eqref{0bounded}, generalizing the isomorphism from \cite{MR3492939}. Applying these results to the $\h$-transform $\G_{\h}$ of the graph $\G,$ we thus obtain the law of the capacity on $\tilde{\G}_{\h}$ of $\psi_{\h}(E_{\h}^{\geq h}(x_0)),$ $h\geq0,$ as well as an isomorphism between the Gaussian free field and the $\h$-transform of random interlacements on the cable system.

\begin{theorem}
\label{couplingintergffh}
Let $\G$ be a transient graph and $\h$ an harmonic function on $\tilde{\G}.$ If \eqref{0bounded} is satisfied, then for all $h\geq0,$
\begin{equation}
\label{eq:laplacecaph}
\begin{split}
   &\E^G_{\tilde{\G}}\left[\exp\left(-u\mathrm{cap}_{\tilde{\G}_{\h}}\big(\psi_{\h}({E}^{\geq h}_{\h}(x_0))\big)\right)\1_{\phi_{x_0}\geq h\times\h(x_0)}\right]\\&=\P^G_{\tilde{\G}}\big(\phi_{x_0}\geq \h(x_0)\sqrt{2u+h^2}\big)\text{ for all }u\geq0\text{ and }x_0\in{\tilde{\G}}.
   \end{split}
\end{equation}
Moreover, \eqref{eq:laplacecaph} for $h=0$ is equivalent to
\begin{equation}
\label{eqcouplingintergffh}
	\begin{gathered}
	\big(\phi_x\1_{x\notin{\mathcal{C}_u^{\h}}}+\sqrt{2\ell_{x,u}^{\h}+\phi_x^2}\1_{x\in{\mathcal{C}_u^{\h}}}\big)_{x\in{\tilde{\mathcal{G}}}}\text{ has the same law under }\P^{I}_{\tilde{\mathcal{G}}_{\h}}\otimes\P^{G}_{\tilde{\G}}
	\\\text{ as }\big(\phi_x+\sqrt{2u}{\h}(x)\big)_{x\in{\tilde{\mathcal{G}}}}\text{ under }\P^G_{\tilde{\mathcal{G}}}\text{ for all }u\geq0,
	\end{gathered}
\end{equation}
where $\mathcal{C}_u^{\h}$ denotes the closure of the union of the connected components of the sign clusters
$\{x\in{\tilde{\G}}:|\phi_x|>0\}$ intersecting $\I_{\h}^u.$
\end{theorem}

\begin{proof}
By \eqref{killedlevelsetsvsh0transform} for $h=0,$ if \eqref{0bounded} holds for $\G,$ then it also holds for $\G_{\h}.$ Applying \cite[Theorem~\ref*{4mainresultcap}]{DrePreRod3} to the graph $\G_{\h}$ and using \eqref{eq:GFFh} and \eqref{killedlevelsetsvsh0transform}, we obtain that \eqref{0bounded} implies \eqref{eq:laplacecaph} for all $h\geq0.$ Moreover, it follows easily from \eqref{eq:localtimesh} that for all $u>0$
\begin{equation}
\label{eq:localtimeRIh}
     (\ell_{x,u}^{\h})_{x\in{\tilde{\G}}}\text{ has the same law under }\P_{\tilde{\G}_{\h}}^{I}\text{ as }(\h(x)^2\ell_{\psi_{\h}(x),u})_{x\in{\tilde{\G}}}\text{ under }\P_{\tilde{\G}_{\h}}^{I}.
\end{equation}
Applying \cite[Theorem~\ref*{4couplingintergff}]{DrePreRod3} to the graph $\G_{\h},$ and using \eqref{eq:GFFh} and \eqref{eq:localtimeRIh}, we obtain that \eqref{eq:laplacecaph} for $h=0$ is equivalent to \eqref{eqcouplingintergffh}. 
\end{proof}

\begin{remark}
\begin{enumerate}[1)]
    \item The isomorphism \eqref{eqcouplingintergffh} generalizes the coupling presented in \eqref{levelsetsvsIu} since it directly implies that
    \begin{equation}
    \label{levelsetsvsIuh}
    \text{ if \eqref{0bounded} holds, then }E^{\geq-\sqrt{2u}}_{\h}\text{ has the same law as }\mathcal{C}_u^{\h}\cup E^{\geq 0}\text{ under }\P^G\otimes\P^I,
\end{equation}
where $\mathcal{C}_u^{\h}$ is defined similarly as $\mathcal{C}_u$ but for the $\h$-transform of random interlacements, see below \eqref{eqcouplingintergffh}.
    \item Following \cite{DrePreRod3}, one could obtain several other results on $E_{\h}^{\geq h},$ $h\in\R,$ when $\h$ is an harmonic function on $\tilde{\G}$: $E^{\geq{h}}_{\h}(x_0)$ is non-compact with $\P^G_{\tilde{\G}}$ positive probability for all $h<0,$ $\mathrm{cap}_{\tilde{\G}_{\h}}\big(\psi_{\h}({E}^{\geq 0}(x_0))\big)<\infty$ $\P^G_{\tilde{\G}}$-a.s.\ for all $x_0\in{\tilde{\G}}$ by \cite[Theorem~\ref*{4mainresult}]{DrePreRod3}, formulas similar to \eqref{eq:laplacecaph} for $h<0$ under condition \eqref{capcondition} for $\G_{\h},$ see \cite[\eqref*{4lawforhnegative} and \eqref*{4eq:capinfinity}]{DrePreRod3}, an equivalence between \eqref{eq:laplacecaph} for $h=0$ and \eqref{eq:laplacecaph} for all $h>0,$ see \cite[\eqref*{4equivisom}]{DrePreRod3} or another formulation of the isomorphism \eqref{eqcouplingintergffh}, see \cite[\eqref*{4eqcouplingintergff}]{DrePreRod3}. Finally, for \emph{any} transient graph $\G,$ one could also obtain a "squared" version of the isomorphism \eqref{eqcouplingintergffh}, both on the discrete graph $G$ as in \cite{MR2892408}, or on the cable system $\tilde{\G}$ as in \cite{MR3502602}. These results could also be extended to all the consequences of Theorem~\ref{couplingintergffh} for $\h=\h_{\text{kill}},$ $\h=\h_{\text{surv}}$ or $\h=a$ gathered in Corollaries~\ref{h0transformiso} and \ref{couplingintergffK} and Theorem~\ref{couplingintergffdim2}.
    \item Let us describe the analogue of the $\h$-transform but for the discrete graph $G.$ We identify $G_{\h}$ and $G$ and define for all continuous-time trajectories $\overline{w}^+$ on $G$ and $t\geq0$
\begin{align*}
    \overline{\theta}_{\h}^{\overline{w}^+}(t)&\stackrel{\mathrm{def.}}{=}\inf\Big\{s\geq0:\,\int_0^s\h\big(\overline{w}^+(u)\big)^{2}\,\mathrm{d}u>t\Big\}
    \\&=\inf\Big\{s\geq0:\,\sum_{x\in{G}}\ell_{x}(s)\h(x)^2>t\Big\},
\end{align*}
with the conventions $\h(\Delta)=0$ and $\inf\varnothing={\zeta},$ and
\begin{equation*}
    (\overline{\xi}_\h(\overline{w}^+))(t)\stackrel{\mathrm{def.}}{=}\overline{w}^+\big({\overline{\theta}_{\h}^{\overline{w}^+}(t)}\big)\text{ for all }t\in{[0,\infty)}.
\end{equation*}
Then the results from Proposition~\ref{corh} still hold when replacing $\xi_{\h}$ by $\overline{\xi}_{\h},$ $\psi_{\h}$ by the identity, the diffusion $X$ by the jump process $Z,$ and the Gaussian free field $(\phi_x)_{x\in{\tilde{\G}}}$ on $\tilde{\G}$ by the Gaussian free field $(\phi_x)_{x\in{G}}$ on $G.$ One can deduce this statement from Proposition~\ref{corh} by using the fact that $Z$ is the trace of $X$ on $G,$ or prove it directly, see the proof of \cite[Proposition~4.6]{LuSaTa} for a proof of a similar statement. We can then also define the $\h$-transform of discrete random interlacements directly with $\overline{\xi}_{\h}$ similarly as in Definition \ref{defhtransforminter}, and, if \eqref{eq:laplacecaph} holds for $h=0,$ obtain a version of the isomorphism \eqref{eqcouplingintergffh} on the discrete graph $G$ between the Gaussian free field on $G$ and the $\h$-transform of discrete random interlacements similar to \cite[\eqref*{4eqcouplingintergffdis}]{DrePreRod3}. 
\end{enumerate}
\end{remark}

Let us now give some some applications of Theorem~\ref{couplingintergffh} for some particular choices of the harmonic function $\h.$ The result \eqref{semigrouph} implies that $\xi_{\h}(X)$ corresponds under $\P^G_{\tilde{\G}_{\h}}$ to the $\h$-transform of $X$ under $\P^G_{\tilde{\G}},$ see for instance \cite[Chapter~11]{MR2152573}, and when $\h=\h_{\text{kill}},$ see \eqref{defh0}, one can then classically relate the law of the diffusion $X$ on $\tilde{\G}$ conditioned on being killed with the $\h_{\text{kill}}$ transform of $X,$ see \cite[Theorem~11.26]{MR2152573}. Therefore, the law of $X$ on $\tilde{\G}$ conditioned on being killed can be related to the diffusion $X$ on the $\h_{\text{kill}}$-transform $\tilde{\G}_{\h_{\text{kill}}}$ of $\tilde{\G},$ and since the proof of this result is short, we include it below for completeness. Similarly, the law of $X$ on $\tilde{\G}$ conditioned on blowing up can be related to the law of $X$ on $\tilde{\G}_{\h_{\text{surv}}}.$

\begin{lemma}
\label{Le:h0transformX}
If $\G$ is a graph with $\h_{\text{kill}}\neq0,$ then the function $\h_{\text{kill}}$ is harmonic on $\tilde{\G}.$ Moreover, for all $x\in{\tilde{\G}},$ the diffusion
\begin{equation}
    \label{h0transformX}
    \xi_{\h_{\text{kill}}}(X)\text{ has the same law under }P_{\psi_{\h_{\text{kill}}}(x)}^{\tilde{\G}_{\h_{\text{kill}}}}\text{ as }X\text{ under }P_x^{\tilde{\G}}(\cdot\,|\,\mathcal{W}_{\tilde{\G}}^{\K,+}).
\end{equation}
The same results also hold when replacing $\h_{\text{kill}}$ by $\h_{\text{surv}}$ and $W_{\tilde{\G}}^{\K,+}$ by $W_{\tilde{\G}}^{\SV,+}.$
\end{lemma}

\begin{proof}
We only do the proof for $\h_{\text{kill}},$ the proof for $\h_{\text{surv}}$ is similar. If $e=\{x,y\}\in{E}$ and $t\in{[0,\rho_e]},$ then the probability beginning in $x+t\cdot I_e$ that $X$ hits $y$ before $x$ is $t\rho_e^{-1},$ see for instance equation 3.0.4 (b), in Part II of \cite{MR1912205}, and by the Markov property $\h_{\text{kill}}(x+t\cdot I_e)=t\rho_{e}^{-1}\h_{\text{kill}}(y)+(1-\rho_{e}^{-1}t)\h_{\text{kill}}(x)=2t\lambda_{x,y}\h_{\text{kill}}(y)+(1-2t\lambda_{x,y})\h_{\text{kill}}(x).$ Similarly, if $x\in{G}$ and $t\in{[0,\rho_x)},$ $\h_{\text{kill}}(x+t\cdot I_x)=2t\kappa_x+(1-2t\kappa_x)\h_{\text{kill}}(x).$ Since $\h_{\text{kill}}(\partial I_x)=1$ if $\kappa_x\neq0,$ we deduce that \eqref{hformulaonedges} holds. Moreover for each $x\in{G}$ the formula \eqref{hformulaoutsideofedge} for $\h=\h_{\text{kill}}$ follows easily from the Markov property for $\hat{Z}$ at time one, and the function $\h_{\text{kill}}$ is thus harmonic on $\tilde\G.$ For all $x\in{\tilde{\G}},$ $t\in{[0,\infty)}$ and functions $f\in{L^2(\tilde{\G},\h^2\cdot m)}$ we have by the Markov property at time $t$
\begin{equation*}
    E^{\tilde{\G}}_x\big[f(X_t)\,|\,W_{\tilde{\G}}^{\K,+}\big]=\frac{1}{\h_{\text{kill}}(x)} E^{\tilde{\G}}_x\big[f(X_t)\1_{W_{\tilde{\G}}^{\K,+}}\big]=\frac{1}{\h_{\text{kill}}(x)}E^{\tilde{\G}}_x\left[f(X_t)\h_{\text{kill}}(X_t)\right],
\end{equation*}
and \eqref{h0transformX} follows from \eqref{semigrouph}.
\end{proof}
\begin{remark}
One can use \eqref{eq:GFFh} for $\h=\h_{\text{surv}},$ which is harmonic when $\h_{\text{kill}}<1$ in view of Lemma~\ref{Le:h0transformX}, to find an alternative proof of \eqref{eq:h0<1thenh*>0}, without using the isomorphism with random interlacements as in \cite{DrePreRod3}. Indeed, following the reasoning of \cite{MR914444}, see also see the Appendix of \cite{MR3765885}, one can directly prove using the Markov property for the Gaussian free field that $E^{\geq h}(x_0)$ contains an unbounded connected components with positive probability for all $h<0$ and $x_0\in{\tilde{\G}}$ on any graph $\G$ with $\kappa\equiv0.$ If $\G$ is a graph such that $\h_{\text{kill}}<1,$ then $\kappa^{(\h_{\text{surv}})}\equiv0,$ and we thus have by \eqref{killedlevelsetsvsh0transform} that $E_{\h_{\text{surv}}}^{\geq h}$ contains an unbounded connected component with $\P^G_{\tilde{\G}}$-positive probability for all $h<0.$ In particular, since $\h_{\text{surv}}\leq 1,$ we obtain that $E^{\geq h}(x_0)$ contains an unbounded connected component with $\P^G_{\tilde{\G}}$-positive probability for all $h<0,$ that is $\tilde{h}\geq0.$

\end{remark}

As a consequence, one can also relate killed random interlacements on $\tilde{\G},$ as defined below \eqref{defIu}, with random interlacements on $\tilde{\G}_{\h_{\text{kill}}},$ and surviving random interlacements on $\tilde{\G}$ with random interlacements on $\tilde{\G}_{\h_{\text{surv}}},$ and apply these results to Theorem~\ref{couplingintergffh}.

\begin{corollary}
\label{h0transformiso}
If $\G$ is a graph with $\h_{\text{kill}}\neq0,$ then the random interlacement process
\begin{equation}
\label{h0transforminter}
\begin{gathered}
    \omega^{\h_{\text{kill}}}\text{ has the same law under }\P^{I}_{\tilde{\G}_{\h_{\text{kill}}}}\text{ as the}
    \\\text{killed random interlacement process under }\P^{I}_{\tilde{\G}},
\end{gathered}
\end{equation}
and for all connected compacts $K\subset\tilde{\G},$
\begin{equation}
\label{capGhvscapkilled}
    \mathrm{cap}_{\tilde{\G}_{\h_{\text{kill}}}}\big(\psi_{\h_{\text{kill}}}(K)\big)=\sum_{x\in{\hat{\partial}}K}\lambda_x\h_{\text{kill}}(x)^2P_x^{\tilde{\G}^{\hat{\partial}K}}(\tilde{H}_{K}=\infty\,|\,W_{\tilde{\G}}^{\K,+}).
\end{equation}
Similar results hold when replacing $\h_{\text{kill}}$ by $\h_{\text{surv}},$ killed random interlacements by surviving random interlacements, and $W_{\tilde{\G}}^{\K,+}$ by $W_{\tilde{\G}}^{\SV,+}.$

In particular, if \eqref{0bounded} holds, or simply \eqref{capcondition} for $\G,$ then \eqref{eq:laplacecaph} for $\h=\h_{\text{kill}}$ provides us with the law of the capacity, given by \eqref{capGhvscapkilled}, of the $\h_{\text{kill}}$ level sets of the Gaussian free field, and \eqref{eqcouplingintergffh} for $\h=\h_{\text{kill}}$ with an isomorphism between the Gaussian free field and killed random interlacements, and similarly for the $\h_{\text{surv}}$ level sets of the Gaussian free field and surviving random interlacements when $\h=\h_{\text{surv}}.$
\end{corollary}

\begin{proof}
By \eqref{defequicap} and \eqref{h0transformX} we have for all finite sets $K\subset G$ and all $x\in{G}$ that
\begin{equation}
\label{eKh0}
\begin{split}
    e_{{\psi_{\h_{\text{kill}}}(K)},\tilde{\G}_{\h_{\text{kill}}}}(\psi_{\h_{\text{kill}}}(x))&=\lambda_{\psi_{\h_{\text{kill}}}(x)}^{(\h_{\text{kill}})}P_{\psi_{\h_{\text{kill}}}(x)}^{{\tilde{\G}_{\h_{\text{kill}}}}}(\tilde{H}_{\psi_{\h_{\text{kill}}}(K)}=\infty)
    \\&=\lambda_x\h_{\text{kill}}(x)^2P_x^{\tilde{\G}}(\tilde{H}_{K}=\infty\,|\,W_{\tilde{\G}}^{\K,+}).
\end{split}
\end{equation}
We thus obtain \eqref{capGhvscapkilled} when $K\subset G,$ is finite, and one can extend it to any compact $K\subset\tilde{\G}$ by considering the graph $\G^{\hat{\partial} K}$ from \eqref{eq:enhancements}.

We now turn to the proof of the identity \eqref{h0transforminter} for random interlacements. By definition of random interlacements, see for instance \eqref{definter}, it is enough to prove that
\begin{equation}
\label{QKh}
    Q_{\psi_{\h_{\text{kill}}}(K),\tilde{\G}_{\h_{\text{kill}}}}\circ(\xi_{\h_{\text{kill}}}^*\circ p_{\tilde{\G}_{\h_{\text{kill}}}}^*)^{-1}(A)=Q_{K,\tilde{\G}}\circ (p_{\tilde{\G}}^*)^{-1}(A\cap W_{\tilde{\G}}^{\K,*})
\end{equation} 
for all finite sets $K\subset G$ and measurable sets $A\subset W_{K,\tilde{\G}}^{*}.$ Using \eqref{h0transformX}, one can easily prove that
\begin{align*}
    P^{\psi_{\h_{\text{kill}}}(K),\tilde{\G}_{\h_{\text{kill}}}}_{\psi_{\h_{\text{kill}}}(x)}(\xi_{\h_{\text{kill}}}(X)\in{\cdot})&=\frac{P_x^{\tilde{\G}}\big(X_{L_K}=x,L_K\in{(0,\zeta)\big)}}{P_x^{\tilde{\G}}\big(X_{L_K}=x,L_K\in{(0,\zeta)},W_{\tilde{\G}}^{\K,+}\big)}P^{K,\tilde{\G}}_x(\cdot,W_{\tilde{\G}}^{\K,+}),
    \\&=\Big(P_x^{K,\tilde{\G}}\big(W_{\tilde{\G}}^{\K,+}\big)\Big)^{-1}P^{K,\tilde{\G}}_x(\cdot,W_{\tilde{\G}}^{\K,+})
    \\&=\frac{P_x^{\tilde{\G}}(\tilde{H}_K=\infty)}{P_x^{\tilde{\G}}\big(\tilde{H}_K=\infty,W_{\tilde{\G}}^{\K,+}\big)}P^{K,\tilde{\G}}_x(\cdot,W_{\tilde{\G}}^{\K,+}),
\end{align*}
where we used \eqref{PxFforZ} in the last equality. Combining with \eqref{eKh0} and Lemma~\ref{Le:h0transformX}, we obtain \eqref{QKh}, and thus \eqref{h0transforminter}. The proof is similar for surviving random interlacements. Using Theorem~\ref{couplingintergffh} and Corollary~\ref{capimplieseverythingh} with $\h=\h_{\text{kill}}$ and $\h=\h_{\text{surv}},$ we can conclude.
\end{proof}

\begin{remark}
\label{rkisokilled}
\begin{enumerate}[1)]
    \item One can characterize the killed random interlacement set similarly as the random interlacement set in \eqref{defIu}, that is the probability that no trajectory in the killed random interlacement process hit a closed set $F$ is given by $\exp(-u\mathrm{cap}_{\tilde{\G}_{\h_{\text{kill}}}}(\psi_{\h_{\text{kill}}}(F))),$ see \eqref{capGhvscapkilled} for an explicit formula, and similarly for surviving random interlacements when replacing $\h_{\text{kill}}$ by $\h_{\text{surv}}.$
    \item When $\kappa\not\equiv0$ and $\{x\in{G}:\,\kappa_x>0\}$ is finite, \eqref{eqcouplingintergffh} for $\h=\h_{\text{kill}}$ can be seen as a reformulation of a signed version of the second Ray-Knight theorem on the cable system. Indeed, one can then define the graph $\G^*$ which corresponds to $\G,$ but replacing the open end of each $I_x,$ $x\in{G}$ with $\kappa_x>0,$ by a common vertex $x_*,$ and using \eqref{h0transforminter}, one can show that the law of the excursions on $G$ of $(X_t)_{t<\tau_u^{x_*}}$ under $P^{\tilde{\G}^*}_{x_*}(\cdot\,|\,\tau_u^{x_*}<\zeta)$ is the same as the law of the trace of the killed random interlacement process on $G$ under $\P^{KI}_{\tilde{\G}},$ where $\tau_u^{x_*}=\inf\{s>0:\ell_{x_*}(s)>u\},$ see \cite[\eqref*{4RidisonG*}]{DrePreRod3} for a proof of a similar statement. One can then easily replace $\ell_{\cdot,u}^{\h_{\text{kill}}}$ by $\ell_{\cdot}(\tau_u^{x_*})$ in \eqref{eqcouplingintergffh}, which corresponds to a version of \cite[Theorem~8]{LuSaTa} on the cable system. In particular, following the proof of \cite[Theorem~8]{LuSaTa}, on any transient graph such that $\kappa\not\equiv0$ and $\{x\in{G}:\,\kappa_x>0\}$ is finite, we obtain that \eqref{eqcouplingintergffh} for $\h_{\text{kill}}$ holds, and thus \eqref{eq:laplacecaph} for $\h_{\text{kill}}$ as well.
    \item One can prove Theorem~\ref{eqcouplingintergffh} for $\h=\h_{\text{kill}}$ directly, without using \cite{DrePreRod3} and Doob $\h$-transform. Indeed, let $K_n,$ $n\in\N,$ be a sequence of finite subsets of $G$ increasing to $G,$ $\kappa^{(n)}=\kappa\1_{K_n},$ and $\G_n$ be the same graph as $\G,$ but with killing measure $\kappa^{(n)}$ instead of $\kappa.$ Since $\{x\in{G}:\kappa^{(n)}_x>0\}$ is finite, as explained before, one can use a version of \cite[Theorem~8]{LuSaTa} on the cable system to obtain \eqref{eqcouplingintergffh} on $\G_n$ for $\h=\h_{\text{kill}}$ and $n\in\N.$ Using the description of killed random interlacements from \eqref{newdescriptionnuK} and Remark~\ref{rqkilled},\ref{desgiveall}), one can compare for each $n\in\N$ the killed interlacement measures on the whole cable system of $\tilde{\G}_n$ and $\tilde{\G},$ instead of their restriction to compacts as in \cite[Lemma~\ref*{4limitKn}]{DrePreRod3}. Proceeding similarly as in the proof of \cite[Lemma~\ref*{4couplingisalwaystrue}]{DrePreRod3}, one can then approximate killed random interlacements on $\tilde{\G}$ by killed random interlacements on the sequence $\tilde{\G}_n,$ decreasing to $\tilde{\G},$ to obtain that, when $\h=\h_{\text{kill}},$ \eqref{eqcouplingintergffh} holds for ${\G}$ if \eqref{0bounded} or \eqref{eq:laplacecaph} holds for $h=0.$ Following the proof of \cite[Proposition~\ref*{4couplingimplytheorem}]{DrePreRod3}, we can then also prove that \eqref{eqcouplingintergffh} implies \eqref{eq:laplacecaph} for all $h\geq0$ when $\h=\h_{\text{kill}}.$ It is less clear how to obtain a direct proof of Theorem~\ref{eqcouplingintergffh} for $\h=\h_{\text{surv}},$ without using $\h$-transforms.
    \item We have $E^{\geq0}=E^{\geq0}_{\h_{\text{kill}}}=E^{\geq0}_{\h_{\text{surv}}},$ and so Theorem~\ref{couplingintergffh} provides us not only with an explicit formula for the capacity of the sign clusters of the Gaussian free field on the cable system of any graph $\G$ satisfying \eqref{0bounded} when $\h\equiv1,$ but also for the capacity of the sign clusters on the graph $\tilde{\G}_{\h_{\text{kill}}}$ when $\kappa\not\equiv0,$ as given in \eqref{capGhvscapkilled}, or on the graph $\tilde{\G}_{\h_{\text{surv}}}$ when $\h_{\text{kill}}<1,$ given similarly as in \eqref{capGhvscapkilled} but with $\h_{\text{surv}}$ and $W_{\tilde{\G}}^{\SV,+}$ instead of $\h_{\text{kill}}$ and $W_{\tilde{\G}}^{\K,+}.$
    \item \label{couplingkilledonA} One also obtains for any $A\subset\overline{G}$ results similar to Corollary~\ref{h0transformiso} for $\h=P_x(X\text{ is killed on }A),$ by replacing $W_{\tilde{\G}}^{\K,+}$ by $\{X\text{ is killed on }A\}$ and killed interlacements by killed on $A$ random interlacements, as defined in Remark~\ref{rqkilled},\ref{rk:killedonARI}), or for $\h=P_x(X\text{ is not killed on }A),$ when replacing $W_{\tilde{\G}}^{\K,+}$ by $\{X\text{ is not killed on }A\}$ and killed interlacements by surviving on $A^c$ random interlacements. This follows from considering the graph $\G^{A_{\infty}}$ as in Remark~\ref{rqkilled},\ref{rk:killedonARI}).
\end{enumerate}
\end{remark}

Let us now give a consequence of Corollary~\ref{h0transformiso}, Corollary~\ref{couplingintergffK}, which is another example of interesting results one can obtain from our Doob $\h$-transform method. It states that, for any compact $K$ of $\tilde{\G},$ one can prove results similar to Theorem~\ref{couplingintergffh} but for the Gaussian free field conditioned on being equal to $0$ on $K$ and the trajectories in the random interlacement process $\omega_u$ avoiding $K.$ In particular, we obtain an isomorphism similar to \eqref{eqcouplingintergffh} between these two objects, which can be seen as a generalization of \cite[Theorem~5.3]{MR3936156}. Recalling the definition of the hitting time $H_K$ from above \eqref{defpartialext}, we define $\h_K(x)=P_x^{\tilde{\G}}(H_K=\zeta)$ for all $x\in{\tilde{\G}},$ $\omega^{(K)}$ the trajectories in the random interlacement process $\omega$ never hitting $K,$ $\ell^{(K)}_{\cdot,u}$ the total local times of the trajectories in $\omega^{(K)}$ with label at most $u$ and $\I^u_{(K)}=\{x\in{\tilde{\G}}:\ell^{(K)}_{x,u}>0\}.$ Let us also define for all compacts $K,K'$ such that $K'\subset K^c,$ $\hat{\partial} K'\subset \hat{\partial} G$ and $\hat{\partial} K\subset\hat{\partial} G,$
\begin{equation}
\label{defcapK}
\begin{gathered}
    e^{(K)}_{K',\tilde{\G}}(x)\stackrel{\mathrm{def.}}{=}\lambda_x\h_K(x)P_x^{\tilde{\G}}\big(\tilde{H}_{K'}={\infty},H_K=\zeta\big)\text{ for all }x\in{\hat{\partial} K'}\\\text{ and }\mathrm{cap}_{\tilde{\G}}^{(K)}(K')=\sum_{x\in\hat{\partial} K'}e_{K',\tilde{\G}}(x).
\end{gathered}
\end{equation}
Using \eqref{eq:enhancements}, one can extend this definition of capacity to any compact $K,K'$ of $\tilde{\G}$ with $K'\subset K^c.$ Finally, let ${E}^{\geq h}_{(K)}(x_0)$ be the cluster of $x_0\in{K^c}$ in $\{x\in{K^c}:\phi_x\geq h\times\h_K(x)\}.$

\begin{corollary}
\label{couplingintergffK}
Let $\G$ be a transient graph satisfying \eqref{capcondition} and $K$ a compact of $\tilde{\G}.$ The identities \eqref{eq:laplacecaph} and \eqref{eqcouplingintergffh} still hold when replacing $\P_{\tilde{\G}}^G$ by $\P_{\tilde{\G}}^G(\cdot\,|\,\phi_{|K}=0),$ $\mathrm{cap}_{\tilde{\G}_{\h}}\big(\psi_{\h}(E_{\h}^{\geq h}(x_0))\big)$ by $\mathrm{cap}_{\tilde{\G}}^{(K)}\big(E_{(K)}^{\geq h}(x_0)\big),$ $\h(x)$ by $\h_K(x),$ $\P_{\tilde{\G}_{\h}}^I$ by $\P_{\tilde{\G}}^I,$ $\ell_{x,u}^{\h}$ by $\ell_{x,u}^{(K)},$ and $\I_{\h}^u$ by $\I_{(K)}^u.$
\end{corollary}

\begin{proof}
By \eqref{eq:enhancements}, we can assume without loss of generality that $\hat{\partial} K\subset G.$ Up to considering each connected component of $K^c$ individually, we will assume that $K^c$ is connected. We also assume that $\h_K>0,$ otherwise the result is trivially true since otherwise $\mathrm{cap}^{(K)}(K')=0$ for all compacts $K'\subset K^c$ and $\I^u_{(K)}=\varnothing$ a.s. We call $\G_{K^c}=(\overline{G}_{K^c},\bar{\lambda}_{K^c},\bar{\kappa}_{K^c})$ the graph such that $(\overline{G}_{K^c},\bar{\lambda}_{K^c})=({G},{\lambda}),$ $\bar{\kappa}_{K^c}={\kappa}$ on $K^c\cap G,$ and $\bar{\kappa}_{K^c}=\infty$ on $K\cap G.$ Similarly as in the beginning of Section~\ref{sec:notation}, we associate to $\G_{K^c}$ an equivalent triplet $(G_{K^c},\lambda_{K^c},\kappa_{K^c})$ with $\kappa_{K^c}<\infty.$ We finally denote by $\G'=\G^{G_{K^c}}=(G',\lambda',\kappa')$ the enhancement of $\G$ containing $G_{K^c}$ in its vertex set, see \eqref{eq:enhancements}. One can then identify $\tilde{\G}_{K^c}$ with $K^c(\subset\tilde{\G})$ and, using \cite[Theorem~4.4.2]{MR2778606}, show that the law of $(X_t)_{t<H_K}$ under $P^{\tilde{\G}}_x$ is $P^{\tilde{\G}_{K^c}}_x$ for all $x\in{K^c},$ and that $(Z_t)_{t<H_K},$ see \eqref{printonG}, has the same law under $P^{\tilde{\G}'}_x$ as $Z$ under $P^{\tilde{\G}_{K^c}}_x$ for all $x\in{G_{K^c}}.$ Moreover, using the Markov property for the Gaussian free field, see \cite[(1.8)]{MR3492939} for instance, one can easily see that
\begin{equation}
    \label{phionGvsphionGAinfinity}
    (\phi_x)_{x\in{K^c}}\text{ has the same law under }\P_{\tilde{\G}}^G(\cdot\,|\,\phi_{|K}=0)\text{ as }(\phi_x)_{x\in{\tilde{\G}_{K^c}}}\text{ under }\P_{\tilde{\G}_{K^c}}^G.
\end{equation}
One can identify trajectories in $W_{\tilde{\G}_{K^c}}^+$ which are not killed on $\hat{\partial} K,$ as defined in Remark~\ref{rqkilled},\ref{rk:killedonARI}), with trajectories in $W_{\tilde{\G}}^+$ which do not hit $\hat{\partial}K,$ which correspond $P_x^{\tilde{\G}}$-a.s.\ to trajectories in $W_{\tilde{\G}}^+$ which do not hit $K$ since $\h_K>0.$ One can thus show that for all compacts $K'$ of $\tilde{\G}_{K^c}$ that
\begin{equation}
\label{eK'Kvssurviving}
    e_{K',\tilde{\G}_{K^c}}(x)P^{K',\tilde{\G}_{K^c}}_x\big(\cdot,X\text{ is not killed on $\hat{\partial}K$}\big)=e_{K',\tilde{\G}}(x)P^{K',\tilde{\G}}_x(\cdot,H_K=\zeta)\text{ for all }x\in{\tilde{\G}_{K^c}}.
\end{equation} 
The proof of \eqref{eK'Kvssurviving} is easy when $\hat{\partial} K'\subset G_{K^c}$ and when considering events involving only the discrete process $Z$ on $G_{K^c}$ since for all $x\in{\hat{\partial}K'}$ by \eqref{defequicap} and \eqref{PxFforZ}
\begin{align*}
    e_{K',\tilde{\G}_{K^c}}(x)P^{K',\tilde{\G}_{K^c}}_x\big(Z\in{\cdot},Z&\text{ not killed on }\hat{\partial}K\big)\hspace{-1mm}
    \\&=\lambda'_xP_x^{\tilde{\G}_{K^c}}(Z\in{\cdot}, Z\text{ not killed on }\hat{\partial}K,\tilde{H}_{K'}=\infty)
    \\&=\lambda'_xP_x^{\tilde{\G}'}(Z\in{\cdot},H_K=\zeta,\tilde{H}_{K'}=\infty)
    \\&=e_{K',\tilde{\G}'}(x)P^{K',\tilde{\G}'}_x(Z\in{\cdot},H_K=\zeta),
\end{align*}
and $Z$ has the same law under $P^{K',\tilde{\G}'}_{\cdot}$ as the trace of $X$ on $G'=G_{K^c}$ under $P^{K,\tilde{\G}}_{\cdot}.$
In more generality, the equality \eqref{eK'Kvssurviving} can be justified for instance using the last exit decomposition \eqref{lastexitdec} from Appendix~\ref{app:inter}. Therefore, by \eqref{defQK} we obtain that
\begin{equation}
\label{survivingisneverhittingA}
\begin{gathered}
       \omega^{(K)}\text{ has the same law under }\P^I_{\tilde{\G}}\text{ as}\\\text{the surviving on $(\hat{\partial} K)^c$ random interlacement process under }\P^{I}_{\tilde{\G}_{K^c}}.
\end{gathered}
\end{equation}
Noting that since $\mathrm{cap}_{\tilde{\G}_{K^c}}(F)\geq\mathrm{cap}_{\tilde{\G}}(F)$ for all $F\subset G\cap K^c,$ we also have that if condition \eqref{capcondition} holds for $\G,$ then it also holds for $\G_{K^c}.$ Using Corollary~\ref{h0transformiso} for the graph $\G_{K^c}^{(\hat{\partial}K)_{<\infty}}$ and the identity 
\begin{equation*}
    \lambda_xP_x^{\tilde{\G}_{K}}(X\text{ is not killed on }\hat{\partial}K)^2P_x^{\tilde{\G}_K^{\hat{\partial} K'}}\big(\tilde{H}_{K'}=\infty\,|\,X\text{ is not killed on }\hat{\partial}K\big)=e_{K',\tilde{\G}}^{(K)}(x)
\end{equation*}
for all $K'\subset G\cap K^c$ and $x\in{\hat{\partial} K'},$ one can easily conclude as explained in Remark~\ref{rkisokilled},\ref{couplingkilledonA}).

\end{proof}

\begin{remark}
One could also find results similar to Corollary~\ref{couplingintergffK} for the Gaussian free field conditioned on being equal to $0$ on $K$ and the trajectories in the surviving random interlacement process not hitting $K,$ replacing $\h_K(x)$ by $P_x^{\tilde{\G}}(H_K=\zeta,W_{\tilde{\G}}^{\SV,+}),$ and adapting the definition of the capacity in \eqref{defcapK}. This can be proved using directly Corollary~\ref{h0transformiso} for surviving random interlacements on $\tilde{\G}_{K^c},$ as defined in the proof of Corollary~\ref{couplingintergffK}. Another possibility is to consider the trajectories in the killed random interlacement process not hitting $K,$ which can be proved using Corollary~\ref{h0transformiso} for killed on $K^c\cap G$ random interlacements on $\tilde{\G}_{K^c}.$
\end{remark}

\section{A non-trivial graph satisfying (\ref{0bounded}) but not (\ref{capcondition})}
\label{sec:Z20}

In this section, we give examples of graphs for which \eqref{capcondition} does not hold, but \eqref{0bounded} holds, thus showing that the implication \eqref{eq:capimplies0bounded} is not an equivalence. We first give trivial examples of such graphs, in the sense that \eqref{0bounded} can still be easily deduced from the property \eqref{eq:capfinite}, even if \eqref{capcondition} does not hold. We then present a condition \eqref{capGfinite}, stronger than the complement of \eqref{capcondition}, under which one cannot deduce \eqref{0bounded} from \eqref{eq:capfinite} anymore, and give an example of a graph satisfying \eqref{capGfinite} and \eqref{0bounded}, by considering the graph $\Z^{2,0}$ corresponding to the two-dimensional lattice $\Z^2$ with constant weights, infinite killing measure at the origin, and zero killing measure everywhere else, see Proposition~\ref{Z20counterexample}. Finally, we deduce Theorem~\ref{couplingintergffdim2} from Theorem~\ref{couplingintergffh}.

Let us first explain what we mean by a non-trivial example, and we first recall the trivial example given in \cite[Remark~\ref*{4R:mainresults1},\ref*{4signwithoutcap}]{DrePreRod3}. Let $0$ be the origin of the graph $\Z^3$ with unit weights and zero killing measure, $A\subset I_0$ be some infinite set without accumulation point, and let $\G^*=(\Z^3)^A,$ see \eqref{eq:enhancements}. The graph $\G^*$ simply corresponds to $\Z^3$ plus the set $A$ attached to $0,$ and $\tilde{\G}^*$ corresponds to $\tilde{\Z}^3,$ plus an additional infinite cable $I_x$ attached to each $x\in{A}(\subset I_0).$  Then by \cite[\eqref*{4capIx}]{DrePreRod3} we have $\mathrm{cap}_{\tilde{\G}^*}(A)=\mathrm{cap}_{\tilde{\Z}^3}(I_0)=\mathrm{cap}_{\tilde{\Z}^3}(\{0\})<\infty,$ and so the graph $\G^*$ does not satisfy condition \eqref{capcondition} in view of \cite[\eqref*{4capconditiondis}]{DrePreRod3}. But all the infinite connected components of ${\Z}^3$ have infinite capacity, and so by \eqref{eq:capfinite} the intersection with $\Z^3(\subset \tilde{\G}^*)$ of each connected component of $E^{\geq0}(\subset\tilde{\G}^*)$ is finite $\P^G_{\tilde{\G}^*}$-a.s. The intersection of each connected component of $E^{\geq 0}$ with $A$ is also finite, since the Gaussian free field on $I_0$ has the same law as a Brownian motion starting in $\phi_0$ with variance $2$ at time $1$ as explained in \cite[Section~\ref*{4subsec:GFF}]{DrePreRod3}, and so \eqref{0bounded} holds for $\G^*.$ Actually, $E^{\geq0}$ is not only bounded but also compact by \cite[Lemma~\ref*{41stpartofmaincor}]{DrePreRod3}.

Using a similar procedure, one could modify any graph $\G$ satisfying \eqref{capcondition} and \eqref{0bounded} such that $\kappa_x=0$ for some vertex $x\in{G}$ by adding infinitely many vertices on the cable $I_x$ to obtain a graph which does not satisfy \eqref{capcondition} anymore, since $I_x$ now correspond to an unbounded set with finite capacity, but for which \eqref{0bounded} can still be deduced from \eqref{eq:capfinite}, since every connected and unbounded components of $G$ have still infinite capacity. More generally, even when \eqref{capcondition} does not hold, one can often get some information from \eqref{eq:capfinite} about the shape of the connected components of $E^{\geq0},$ which may help in proving \eqref{0bounded}. However, if we assume that
\begin{equation}
\label{capGfinite}
    \mathrm{cap}(F)<\infty\text{ for some closed and connected set }F\subset\tilde{\G}\text{ with bounded complement,}
\end{equation}
then the result \eqref{eq:capfinite} is equivalent to $\P^G(\mathrm{cap}\big(E^{\geq0}(x_0)\cap F^c\big)<\infty)=1$ for all $x_0\in{\tilde{\G}},$ which does not provide us with any information on the boundedness of $E^{\geq 0}(x_0)$ since $F^c$ is bounded. Note that if $\G$ is unbounded, then \eqref{capGfinite} implies that \eqref{capcondition} does not hold. What we mean with "a non-trivial graph satisfying \eqref{0bounded} but not \eqref{capcondition}" is thus a graph satisfying \eqref{0bounded} and \eqref{capGfinite}. Such a graph does not percolate at level $0,$ but one cannot deduce this result from the information \eqref{capGfinite}, as opposed to graphs satisfying (or almost satisfying) \eqref{capcondition}, and we now present an example of such a graph.

Let $\Z^{2}$ be the graph with weights $\frac14$ between neighbours in $\Z^2$ and $0$ killing measure, and $\Z^{2,0}$ be the same graph as $\Z^2$ but with infinite killing at the origin. Identifying $\tilde{\Z}^{2,0}$ with $\tilde{\Z}^2\setminus I_0,$ let us denote for each $n\in\N$ by $B(n)$ the subset of $\tilde{\Z}^{2,0}$ identified with $\{x\in{\tilde{\Z}^2\setminus I_0}:\,0<d_{\tilde{\Z}^2}(0,x)\leq n\}.$ With a slight abuse of notation, let us also denote by $\textbf{a}$ the restriction to $\tilde{\Z}^{2,0}(\subset\tilde{\Z}^2)$ of the potential kernel $\textbf{a}$ defined page \pageref{pagedefa}. Since $\textbf{a}$ is linear on $I_e$ for each edge and vertex $e$ of $\Z^{2,0},$ it satisfies \eqref{hformulaonedges}, and since it also satisfies \eqref{hformulaoutsideofedge} by \cite[Proposition~4.4.2]{MR2677157}, we obtain that $\textbf{a}$ is harmonic on $\tilde{\Z}^{2,0},$ in the sense of Definition \ref{defharmonic}. One can relate the capacity of a set on the graph $\Z_\textbf{a}^{2,0},$ as defined in Definition \ref{harmo}, to the usual definition of the two-dimensional capacity. Indeed, let us define for all closed sets $A\subset\tilde{\Z}^2$ such that $0\in{A}$
\begin{equation}
\label{capdim2}
    \mathrm{cap}_{\tilde{\Z}^2}(A)\stackrel{\mathrm{def.}}{=}\mathrm{cap}_{\tilde{\Z}^{2,0}_\textbf{a}}\big(\psi_\textbf{a}(A\setminus I_0)\big).
\end{equation}
This definition is coherent with the usual definition of the two-dimensional capacity $\mathrm{cap}_{\Z^2}$ from \cite[Section~6.6]{MR2677157}, since by \cite[Proposition~2.2]{MR3475663}
\begin{equation}
    \label{capdim2rest}
    \mathrm{cap}_{\tilde{\Z}^2}(A)=\mathrm{cap}_{\Z^2}(A)\text{ for all finite sets }A\subset\Z^2\text{ with }0\in{A}.
\end{equation} 
Using Corollary~\ref{capimplieseverythingh}, we can now easily show that $\Z^{2,0}$ satisfies \eqref{0bounded}.

\begin{proposition}
\label{Z20counterexample}
\eqref{0bounded} and \eqref{capGfinite} are satisfied for the graph $\Z^{2,0}.$ 
\end{proposition}
\begin{proof}
It is clear that $n\mapsto\mathrm{cap}_{\Z^{2,0}}\big(\overline{B(n)\setminus B(1)}\big)$ is constant since a trajectory started on the boundary of $B(n)$ will come back in $B(n)$ with probability $1,$ and so $\overline{B(1)^c}$ has finite capacity, from which \eqref{capGfinite} readily follows.

It follows from \cite[Lemma~6.6.7,(b)]{MR2677157}, \eqref{capdim2} and \eqref{capdim2rest} that for all $n\in\N$ and connected sets $A_n\subset\Z^{2,0}_\textbf{a}$ with diameter at least $n$ containing $\{1,0\}$
\begin{equation*}
    \mathrm{cap}_{\tilde{\Z}_\textbf{a}^{2,0}}(A_n)\geq C\log(n+1),
\end{equation*}
for some constant $C>0.$ If $A\subset \Z^{2,0}_\textbf{a}$ is infinite and connected, let us denote by $A'\subset \Z^{2,0}_\textbf{a}$ a finite connected set connecting $A$ to $\{1,0\}.$ By subadditivity of the capacity, see \cite[Proposition~6.6.2]{MR2677157}, we have that 
\begin{equation*}
    \mathrm{cap}_{\tilde{\Z}_\textbf{a}^{2,0}}(A)\geq \mathrm{cap}_{\tilde{\Z}_\textbf{a}^{2,0}}(A\cup A')-\mathrm{cap}_{\tilde{\Z}_\textbf{a}^{2,0}}(A')\geq C\log(n+1)-\mathrm{cap}_{\tilde{\Z}_\textbf{a}^{2,0}}(A')
\end{equation*} 
for all $n\in\N$ since the diameter of $A\cup A'$ is at least $n.$ Using \cite[\eqref*{4capconditiondis}]{DrePreRod3} and letting $n\rightarrow\infty,$ we obtain that condition \eqref{capcondition} holds for the graph $\Z_\textbf{a}^{2,0},$ and so \eqref{0bounded} holds for $\Z^{2,0}$ by Corollary~\ref{capimplieseverythingh}. 
\end{proof}

\begin{remark}
\label{rk:Z20counterexample}
\begin{enumerate}[1)]
\item By Corollary~\ref{capimplieseverythingh}, if there exists an harmonic function $\h$ on $\tilde{\G}$ such that \eqref{capcondition} is satisfied for $\G_{\h},$ then \eqref{0bounded} holds, and it is an interesting open question to know whether a non-trivial graph satisfying \eqref{0bounded} but not \eqref{capcondition} on $\G_{\h}$ for all harmonic functions $\h$ on $\tilde{\G}$ exists. Note that $\Z^{2,0}$ is not such a graph, see the proof of Proposition~\ref{Z20counterexample}. 
\item \label{rk:h*-inftyZ20} Actually one can show that on the graph $\mathbb{Z}^{2,0},$ not only is \eqref{0bounded} satisfied, but actually $\tilde{h}_*=-\infty.$ Indeed, by \eqref{levelsetsvsIu} $E^{\geqslant -\sqrt{2u}}$ has unbounded connected components with positive probability if and only if the closure of the union of the clusters $\{x\in{\tilde{\G}}:\,|\phi_x|>0\}$ intersecting $\mathcal{I}^u$ has unbounded connected components with positive probability. But by \eqref{0bounded}, all the clusters of $\{x\in{\tilde{\G}}:\,|\phi_x|>0\}$ are bounded, and since $\Z^2$ is recurrent, the trajectories of $\mathcal{I}^u$ in $\mathbb{Z}^{2,0}$ are all finite (they correspond on $\Z^2$ to trajectories starting in $0$ and ending the first time they reach $0$), and we can conclude. There are also graphs for which $\tilde{h}_*=0,$ and so that \eqref{0bounded} cannot be directly deduced from \eqref{eq:capfinite}. For instance, instead of taking $\kappa_0=\infty$ on $\Z^2,$ one can consider the graph $\G$ obtained by adding to $\Z^2$ a distinct copy of $\Z^3,$ and identifying their respective origin so that they are equal. Indeed on this graph \eqref{0bounded} is satisfied by a similar reasoning as in Proposition~\ref{Z20counterexample} (extending the definition of $\mathbf{a}$ to be equal to $0$ on $\Z^3$); and $E^{\geqslant-\sqrt{2u}}$ is a.s.\ unbounded for all $u>0$ since $\mathcal{I}^u$ a.s.\ always contains unbounded connected components in $\Z^3,$ and thus $\tilde{h}_*=0.$ Moreover $\mathrm{cap}_{\tilde{\G}}(\Z^2)=\mathrm{cap}(\{0\})<\infty$ (but \eqref{capGfinite} is not satisfied since $\mathrm{cap}_{\tilde{\G}}(\Z^3)=\infty$); and so \eqref{eq:capfinite} only implies that $E^{\geqslant 0}\cap\Z^3$ is finite,  which is not enough to directly deduce \eqref{0bounded}.
\end{enumerate}
\end{remark}

One can also relate the Gaussian free field on $\tilde{\Z}^{2,0}$ with the pinned field $\phi^p$ defined in \eqref{def2dgff}. 

\begin{lemma}
\label{le:phiZ2otherdef}
The field $(\phi_x^p)_{x\in{\tilde{\Z}^2}\setminus I_0}$ has the same law under $\P^{G,p}$ as $(\phi_x)_{x\in\tilde{\Z}^{2,0}}$ under $\P^G_{\tilde{\Z}^{2,0}}.$
\end{lemma}
\begin{proof}
Let $\Z_{n}^{2,0}$ the same graph as $\Z^{2,0},$ but with infinite killing measure outside of $B(n).$ The Gaussian free field under $\P^G_{\tilde{\Z}_n^{2,0}}$ converges in law to the Gaussian free field under $\P^G_{\tilde{\Z}^{2,0}}$ since by the Markov property at the first time $H_{B(n)^c}$ that $X$ hits $B(n)^c$ we have for all $x,y\in{\tilde{\Z}^{2,0}}$
\begin{align*}
    g_{\tilde{\Z}^{2,0}}(x,y)-g_{\tilde{\Z}^{2,0}_n}(x,y)&=E_x^{\tilde{\Z}^{2,0}}\big[g_{\tilde{\Z}^{2,0}}(X_{H_{B(n)^c}},y)\1_{H_{B(n)^c}<\zeta}\big]
    \\&\leq g_{\tilde{\Z}^{2,0}}(y,y)P_x^{\tilde{\Z}^{2,0}}(H_{B(n)^c}<\zeta)\tend{n}{\infty}0.
\end{align*}
Moreover, it follows from \cite[(5.30)]{MR3936156}, whose proof can easily be extended to the cable system, that the Gaussian free field under $\P^G_{\tilde{\Z}_n^{2,0}}$ converges in law to $(\phi_x^p)_{x\in{\tilde{\Z}^2}\setminus I_0}$ under $\P^{G,p},$ and we can conclude.
\end{proof}

Combining Lemma~\ref{le:phiZ2otherdef} with Proposition~\ref{Z20counterexample} and Theorem~\ref{couplingintergffh} will let us easily deduce Theorem~\ref{couplingintergffdim2}. Recalling Definition \ref{defhtransforminter}, let us first define the the two-dimensional random interlacement process $\omega^{(2)}$ under a probability $\P^{I,2}$ as a point process of trajectories on $\tilde{\Z}^2$ such that
\begin{equation}
\label{defRIdim2}
\omega^{(2)}\text{ has the same law under }\P^{I,2}\text{ as }\omega^\textbf{a}\text{ under }\P_{\tilde{\Z}^{2,0}_\textbf{a}}^I,
\end{equation}
where we identified trajectories on $\tilde{\Z}^{2,0}$ with trajectories on $\tilde{\Z}^2$ avoiding $I_0.$ This definition is coherent with the previous definitions of two-dimensional random interlacements on the discrete graph $\Z^2,$ since the trace of $\omega^{(2)}$ on $\Z^2$ corresponds to the interlacement process defined above \cite[Corollary~4.3]{MR3936156}, and the associated discrete time skeleton to the interlacement process defined above \cite[Definition 2.1]{MR3475663}. The characterization \eqref{defIudim2} of the corresponding interlacement set $\I_2^u$ moreover directly follows from \eqref{defIu}, Definition \ref{defhtransforminter}, \eqref{capdim2} and \eqref{defRIdim2}.

\begin{proof}[Proof of Theorem~\ref{couplingintergffdim2}]
Since $\textbf{a}(0)=0,$ we have $\kappa^{(\textbf{a})}\equiv0$ on $\Z^{2,0},$ and so $E_\textbf{a}^{\geq h}$ contains an unbounded connected component with $\P^{G}_{\Z^{2,0}}$ positive probability for all $h<0$ by Corollary~\ref{capimplieseverythingh}. Since \eqref{0bounded} holds for $\Z^{2,0}$ by Proposition~\ref{Z20counterexample}, one can combine Lemma~\ref{le:phiZ2otherdef} with the results from Theorem~\ref{couplingintergffh} for $\G=\Z^{2,0}$ and $\h=\textbf{a}$ with \eqref{capdim2} and \eqref{defRIdim2} to obtain Theorem~\ref{couplingintergffdim2}, noting that \eqref{eq:laplacecapkilleddim2} and \eqref{eqcouplingintergffdim2} trivially extend to $I_0$ since $\mathrm{cap}_{\tilde{\Z}^2}(I_0)=0,$ $\I_2^u\cap I_0=\varnothing$ $\P^{I,2}$-a.s.\ and $\textbf{a}(x)=0$ for all $x\in{I_0}.$
\end{proof}

\begin{remark}
\label{endrkdim2}
\begin{enumerate}[1)]
    \item One could also try to prove Theorem~\ref{couplingintergffdim2} using Corollary~\ref{couplingintergffK}, similarly as the proof of Theorem~5.5 from Theorem~5.3 in \cite{MR3936156}. Indeed, one can easily deduce from \cite[Lemma~3.1]{MR3936156} that for all finite sets $K\subset{\Z^2\setminus\{0\}}$
    \begin{equation*}
        \lim\limits_{n\rightarrow\infty}\frac{4}{\pi^2}\log^2(n)\mathrm{cap}_{\tilde{\Z}^2_n}^{(\{0\})}(K)=\mathrm{cap}_{\Z^2}(K\cup\{0\}),
    \end{equation*}
    by \cite[(3.7)]{MR3936156} that $2\pi^{-1}\log(n)P_x^{\tilde{\Z}^{2}_n}(H_0=\zeta)\tend{n}{\infty}\textbf{a}(x)$ for all $x\in{\Z^2},$ and by \cite[(5.30)]{MR3936156} that $(\phi_x)_{x\in{\Z}^{2}_n}$ under $\P^G_{\tilde{\Z}^{2}_n}(\cdot\,|\,\phi_0=0)$ converges in law to $(\phi_x^p)_{x\in{{\Z}^2}}$ under $\P^{G,p}.$ Therefore if one could extend the previous results to the cable system, taking the limit in the version of \eqref{eq:laplacecaph} from Corollary~\ref{couplingintergffK} with $K=\{0\},$ $\G=\Z_n^2,$ $u=4\pi^{-2}\log^2(n)u'$ and $h=2\pi^{-1}\log(n)h'$ and extending the previous results to the cable system could give us \eqref{eq:laplacecapkilleddim2}. Moreover, by  \cite[Corollary~4.3]{MR3936156}, the law of the trace of $\omega_{4\pi^{-2}\log^2(n)u'}^{(0)}$ on $\Z_n^2$ under $\P^G_{\tilde{\Z}_n^2}$ converges to the law of the trace of $\omega^{(2)}_{u'}$ on $\Z^2$ under $\P^{I,2},$ and one could also try similarly to prove \eqref{eqcouplingintergffdim2} by taking the limit in the version of \eqref{eqcouplingintergffh} from Corollary~\ref{couplingintergffK}. 
    
    This strategy can be effectively applied to prove that the squares of the processes on both side of \eqref{eqcouplingintergffh} have the same law, which corresponds to proving \cite[Theorem~5.5]{MR3936156} but on the cable system. However, the weak convergence results for the Gaussian fields and random interlacements from \cite{MR3936156} seem not robust enough to prove rigorously that Theorem~\ref{couplingintergffdim2} follows from Corollary~\ref{couplingintergffK} using the previously explained strategy, see for instance the proof of \cite[Lemma~\ref*{4couplingisalwaystrue}]{DrePreRod3} which requires more robust convergence results, and using instead the $\textbf{a}$-transform directly on the graph $\Z^{2,0}$ solves this problem. 
    
    \item Some results similar to \eqref{eq:laplacecapkilleddim2} and \eqref{eqcouplingintergffdim2} also holds but for the usual level sets $E^{p,\geq h}=\{x\in{\tilde{\Z}^2}:\,\phi_x^{p}\geq h\}$ of the pinned free field. Indeed \eqref{0bounded} holds on $\Z^{2,0}$ by Proposition~\ref{Z20counterexample}, and using \cite[Theorem~\ref*{4mainresultcap}]{DrePreRod3} and Lemma~\ref{le:phiZ2otherdef}, one can thus obtain the law of the capacity (in terms of the graph $\Z^{2,0}$) of the usual level sets $E^{p,\geq h}$ of the pinned free field for all $h\geq0.$ Moreover, by \cite[Theorem~\ref*{4couplingintergff}]{DrePreRod3}, one can also obtain an isomorphism similar to \eqref{eqcouplingintergffdim2} but between the pinned free field on $\Z^{2}$ and random interlacements on $\Z^{2,0},$ and replacing $\sqrt{2u}\textbf{a}(x)$ by $\sqrt{2u}.$
    
    \item \label{normalpinnedlevelsets}One may also want to investigate percolation for the usual level sets $E^{p,\geq h}$ of the pinned free field. Since $E^{p,\geq0}=E^{p,\geq0}_\textbf{a},$ it follows from Theorem~\ref{couplingintergffdim2} and monotonicity that $E^{p,\geq h}$ contains $\P^{G,p}$-a.s.\ only bounded connected components for all $h\geq0.$ Moreover, using \eqref{levelsetsvsIu} and Lemma~\ref{le:phiZ2otherdef}, it is clear that for all $u\geq0,$ $E^{p,\geq-\sqrt{2u}}$ contains unbounded connected components with $\P^{G,p}$ positive probability if and only if $\mathcal{C}_u^p$ contains unbounded connected components with $\P^{G,p}\otimes\P^{I}_{\tilde{\Z}^{2,0}}$ positive probability, where $\mathcal{C}_u^p$ denotes the closure of the union of the connected components of $\{x\in{\tilde{\G}}:|\phi_x^{p}|>0\}$ intersecting $\I^u.$ Since $\h_{\text{kill}}\equiv1$ on $\tilde{\Z}^{2,0},$ the interlacements on $\tilde{\Z}^{2,0}$ consist only of trajectories whose forwards and backwards parts have been killed in $\{0\}$ and so $\I^u$ is $\P^I_{\tilde{\Z}^{2,0}}$-a.s.\ bounded. Therefore, since $\{x\in{\tilde{\G}}:|\phi_x^{p}|>0\}$ contains $\P^{G,p}$-a.s.\ only bounded connected components, we obtain that $E^{p,\geq h}$ contains $\P^{G,p}$-a.s.\ only bounded components for all $h\in\R,$ that is the associated critical parameter is equal to $-\infty.$
\end{enumerate}
\end{remark}

\section{A graph with infinite critical parameter}
\label{sec:example}
In this section, we give an example of a graph for which the critical parameter $\tilde{h}_*,$ see \eqref{defh*}, is strictly positive, and in fact infinite, thus proving that the implication \eqref{eq:capimplies0bounded} is not trivial. For any $\alpha\in{(0,1)}$ and $d\in\N,$ $d\geq2,$ we define $\mathbb{T}_d^{\alpha}$ as the rooted $(d+1)$-regular tree, such that, denoting by $T_n$ the set of vertices in $\mathbb{T}_d^{\alpha}$ at generation $n$ (seen from the root),
\begin{equation*}
    \lambda_{x,y}^{(\alpha)}=\alpha^n\text{ if }x\in{T_n}\text{ and }y\in{T_{n+1}},
\end{equation*}
and $0$ otherwise. Moreover, we take $\kappa^{(\alpha)}=0$ if $\alpha>\frac{1}{d}$ and $\kappa^{(\alpha)}=\1_{\varnothing}$ otherwise, where $\varnothing$ denotes the root of the tree. Since for $x\in{T_n}$, $n \geq 1$, and $\alpha \in (0,1),$
\begin{equation}
\label{eq:RWdrift}
    P_x^{\mathbb{T}_d^{\alpha}}(\hat{Z}_1\in{T_{n+1}})=d\frac{\alpha^n}{\alpha^{n-1}+d\alpha^n}=\frac{d\alpha}{1+d\alpha} \left({>}\frac{1}{2}\text{ if }\alpha>\frac{1}{d}\right),
\end{equation}
we have that $\mathbb{T}_d^{\alpha}$ is a transient graph for all $\alpha\in{(0,1)}$ and $d\in\N,$ $d\geq2.$ Taking $A$ to be an infinite connected line containing exactly one vertex per generation, and noting that the equilibrium measure of any sets at a point $x\in{T_n}$ is at most $\alpha^{n-1}+d\alpha^n,$ one can easily show that $\mathrm{cap}(A)<\infty,$ and so condition \eqref{capcondition} is not satisfied for $\T_d^{\alpha}.$
\begin{proposition}
\label{h*infinity}
There exists a constant $\overline{C}<\infty,$ such that for any $\alpha\in{(0,1)}$ and $d\in\N,$ $d\geq2,$ with 
\begin{equation}
\label{condondalpha}
    d\Big(1-\exp\Big(-\frac{\sqrt{\alpha}}{d\alpha+1}\Big)\Big)>\overline{C},
\end{equation} 
the set $E^{\geq h}$ contains $\P^G_{\tilde{\mathbb{T}}_d^{\alpha}}$-a.s.\ an unbounded connected component for all $h\in\R,$ and so $\tilde{h}_*=\infty.$ 
\end{proposition}
\begin{proof}
    Using the Markov property for the Gaussian free field, see \cite[(1.8)]{MR3492939} for instance, one can construct the Gaussian free field on $(\tilde{\mathbb{T}}_d^{\alpha})^{\Ed}$ (where $(\tilde{\mathbb{T}}_d^{\alpha})^{\Ed}$ is defined as in page \pageref{deftildeged} by removing the edges $I_x,$ $x\in{\mathbb{T}_d^\alpha},$ from $\tilde{\mathbb{T}}_d^{\alpha}$) recursively in the generation $n$ as follows. Let $Y_x,$ $x\in{\mathbb{T}_d^{\alpha}},$ be a family of i.i.d.\ $\mathcal{N}(0,1)$-distributed random variables under $\P,$ and let $\psi_0=Y_0(g_{\mathbb{T}_d^{\alpha}}(0,0))^{1/2}.$ Recursively in~$n \geq 0$, we then define 
    \begin{equation}
    \label{eq:tree_markov}
        \psi_x\stackrel{\mathrm{def.}}{=}\psi_{x^-}P_x(H_{\{x^-\}}<\infty)+Y_x\sqrt{g_{T_n^c}(x,x)}, \text{ for all }x\in{T_{n+1}},
    \end{equation}
    where $x^-$ is the first ancestor of $x,$ i.e., the neighbour of $x$ on a geodesic path from $x$ to $0,$ and $g_{T_n^c}(x,x)$ is the Green function defined as in \eqref{Greendef} but for diffusion killed on exiting $T_n^c.$ Using the Markov property for the Gaussian free field, see \cite[(1.8)]{MR3492939} for instance, one can then easily prove that $(\psi_x)_{x\in{\mathbb{T}_d^{\alpha}}}$ has the same law as $(\phi_x)_{x\in{\mathbb{T}_d^{\alpha}}}$ under $\P^G_{\mathbb{T}_d^{\alpha}}$. Moreover, let $B^e,$ $e\in{E},$ be a family of independent processes, such that for each edge $e=\{x,y\}\in{E}$ between $x\in{T_n}$ and $y\in{T_{n+1}},$ $B^e$ is a Brownian bridge of length $\frac{1}{2\alpha^n}$ between $0$ and $0$ of a Brownian motion with variance $2$ at time $1,$ and let
    \begin{equation*}
        \psi_{x+t\cdot I_e}=2\alpha^nt\psi_y+(1-2\alpha^nt)\psi_x+B_{t}^e\text{ for all }t\in{\big[0,1/2\alpha^n\big]}
    \end{equation*}
    (cf.\ the beginning of Section~\ref{sec:notation} for notation). Then $(\psi_x)_{x\in{(\tilde{\mathbb{T}}_d^{\alpha}})^{\Ed}}$ has the same law as $(\phi_x)_{x\in{(\tilde{\mathbb{T}}_d^{\alpha}})^{\Ed}}$ under $\P^G_{\mathbb{T}_d^{\alpha}};$ cf.\  \cite[Section~2]{DrePreRod} for the proof of an analogous construction on $\Z^d,$ $d\geq3$. Now for each $x\in{T_{n+1}},$ with $A_x=\{ \psi_{x}\geq(\lambda_{x}^{(\alpha)})^{-1/2} \}$, in view of \eqref{eq:tree_markov} we have that
    \begin{equation*}
        \P\big(A_x\,\big|\,\psi_{x^-}\big)\1_{ A_{x^-}}\geq\P\big(Y_x\geq(\lambda_{x}^{(\alpha)}g_{T_n^c}(x,x))^{-1/2}\big)\geq\P(Y_0\geq 1);
    \end{equation*}
    indeed, the inequality on the right-hand side follows since under $P_x^{\tilde{\T}_d^{\alpha}}$, $Z$ spends at least an exponential time with parameter $\lambda_x^{(\alpha)}$ in $x\in T_{n+1}$ before hitting $T_n$ and so $g_{T_n^c}(x,x)\lambda_x^{(\alpha)}\geq1.$ Moreover, using the exact
formula for the distribution of the maximum of a Brownian bridge, see e.g. \cite[Chapter~IV.26]{MR1912205}, we have for all $x\in{T_{n+1}}$ and $n$ large enough, writing $e=\{x,x^-\},$ on the event $A_x \cap A_{x^-}$,
    \begin{align*}
        \P\big(\psi_y\geq (\lambda_x^{(\alpha)})^{-\frac14} \, \forall\,y\in{I_{e}}\,|\,\psi_x,\psi_{x^-}\big)&=1-\exp\big(-2\alpha^n(\psi_x-(\lambda_x^{(\alpha)})^{-\frac14})(\psi_{x^-}-(\lambda_x^{(\alpha)})^{-\frac14})\big)
        \\&\geq1-\exp\big(-\alpha^n(\lambda_x^{(\alpha)})^{-\frac12}(\lambda_{x^-}^{(\alpha)})^{-\frac12}\big)
        \\&= 1-\exp\Big(-\frac{\sqrt{\alpha}}{d\alpha+1}\Big).
    \end{align*}
 Hence, for all $n$ large enough and any $y\in{T_{n+1}},$ the intersection of the cluster of $y$ in $\{x\in{(\tilde{\mathbb{T}}_d^{\alpha}})^{\Ed}:\,\psi_x\geq(\lambda_x^{(\alpha)})^{-\frac14}\}$ with $\mathbb{T}_d^{\alpha}$ stochastically dominates an independent Galton-Watson tree rooted at $y$, with average number of children equal to
    \begin{equation*}
        d\Big(1-\exp\Big(-\frac{\sqrt{\alpha}}{d\alpha+1}\Big)\Big)\P(Y_0\geq 1).
    \end{equation*}
    Choosing $\overline{C}=\P(Y_0\geq1)^{-1},$ we thus obtain that $\{x\in{(\tilde{\mathbb{T}}_d^{\alpha})^-}:\,\psi_x\geq(\lambda_x^{(\alpha)})^{-\frac14}\}$ contains {$\P$-a.s.}\ an unbounded connected component if $d\big(1-\exp\big(-\frac{\sqrt{\alpha}}{d\alpha+1}\big)\big)>\overline{C},$ and since $\lambda_x^{(\alpha)}\rightarrow{0}$ as $d(x,0)\rightarrow{\infty},$ we deduce that $\tilde{h}_*=\infty.$ 
    \end{proof}

\begin{remark}
\label{endremark}
\begin{enumerate}[1)]
    \item \label{alphadexist}In both cases, $\alpha>\frac1d$ as well as $\alpha\leq\frac1d,$ it is possible to find $\alpha\in{(0,1)}$ and $d\in\N,$ $d\geq2$ such that \eqref{condondalpha} holds. For instance, one can take $\alpha=\frac{a}{d}$ for arbitrary fixed $a>0,$ and choose $d$ large enough. In particular, in view of \eqref{eq:RWdrift}, this provides us with graphs such that $\kappa\equiv0$ and $\tilde{h}_*=\infty$ when $\alpha>\frac1d,$ or with $\h_{\text{kill}}\equiv1$ and $\tilde{h}_*=\infty$ when $\alpha\leq\frac1d.$
    
    \item By means of suitable enhancements, cf.\ \eqref{eq:enhancements}, one readily derives from Proposition~\ref{h*infinity} an example of a graph $\G$ with unit weights and zero killing, on which $\tilde{h}_*=\infty$. Indeed, fix some $d\in2\N$ such that $d\big(1-\exp\big(-\frac{\sqrt{2/d}}{3}\big)\big)>\overline{C}.$ Consider the set $A \subset (\tilde{\mathbb{T}}_d^{2/d})^{\Ed}$ (attached to the weights $\lambda^{(2/d)}, \kappa^{(2/d)}$ as above) with $A \cap I_{\{x^-,x\}}=\{ x^-+(k/2)\cdot I_{\{x^-,x\}}: 1\leq k \leq (d/2)^{n-1}-1\},$ for all $x\in{T_n}$, $n\geq 1$. Then \cite[\eqref*{4eq:GAsubsetG}]{DrePreRod3} yields that $ \tilde{\mathbb{T}}_d^{2/d} \subset \tilde{\G}$, where $\G\stackrel{\text{def.}}{=}(({\mathbb{T}}_d^{2/d})^A, \lambda^A,\kappa^A)$, with $\lambda^A \equiv 1$ and $\kappa^A\equiv 0$, whence $\tilde{h}_*(\tilde{\G})=\infty$ by Proposition~\ref{h*infinity}.
    \item \label{counterexampleh*cap<h*}By \cite[Theorem~\ref*{4mainresult}]{DrePreRod3}, we have that $\mathrm{cap}(E^{\geq 0}(x_0))<\infty$ $\P^G_{\tilde{\T}^{\alpha}_d}$-a.s.\ for all $x_0\in{\tilde{\G}}$, and so $\tilde{h}_*^{{\rm cap}}\leq0,$ where $\tilde{h}_*^{{\rm cap}}$ is the critical parameter associated to the percolation of $E^{\geq h}$ in terms of capacity, see \cite[\eqref*{4defh*cap}]{DrePreRod3}. In particular, when $\alpha>\frac1d$ and \eqref{condondalpha} holds, then $\T_d^{\alpha}$ is an example of a graph for which the inequality \cite[\eqref*{4capcomboukappa}]{DrePreRod3} is strict.
    \item \label{h*infiniteandlaw0}In \cite[Corollary~\ref*{4dichotomy}]{DrePreRod3}, it is proved that, under the condition that \eqref{eq:laplacecaph} holds for $h=0$ and $\h\equiv1,$ \eqref{0bounded} is not satisfied implies $\tilde{h}_*=\infty.$ If one could prove that there exist $\alpha\in{(0,1)}$ and $d\geq2$ fulfilling \eqref{condondalpha} such that $\mathbb{T}_d^{\alpha}$ satisfies \eqref{eq:laplacecaph} for $h=0$ and $\h\equiv1,$ this would show that this implication in \cite[Corollary~\ref*{4dichotomy}]{DrePreRod3} is not trivial. It is also an interesting question whether there exist graphs on which $\tilde{h}_*\in{\{0,\infty\}}$. We hope to come back to these questions soon.
\end{enumerate}
\end{remark}

\appendix
\renewcommand*{\thetheorem}{A.\arabic{theorem}}
\section{Appendix: Proof of Theorem~\ref{nuexists}}
\label{app:inter}
Before starting the proof of Theorem~\ref{nuexists}, let us recall some interesting facts from the theory of diffusions, which will lead to the important last exit decomposition \eqref{lastexitdec} of the diffusion $X$ on $\tilde{\G},$ that is a decomposition of the law of $X$ before and after the time $L_F$ at which $X$ leaves a closed set $F$ of $\tilde{\G}$. Using \cite[Theorems~4.1.2 and 4.2.4]{MR2778606}, one can associate to the diffusion $X$ a symmetric family of probability densities $(p_t(x,y))_{t>0},$ $x,y\in{\tilde{\G}},$ such that
\begin{equation}
\label{greenpt}
    P_x(X_t\in{\mathrm{d}}y)=p_t(x,y)m(\mathrm{d}y),\text{ and then } g(x,y)=\int_0^{\infty}p_t(x,y)\diff t.
\end{equation}
The fact that the formula \eqref{greenpt} for the Green function holds, recall the definition \eqref{Greendef}, can be for instance deduced from \cite[Theorem~3.6.5]{MR2250510}. Let us now recall some useful results from \cite[Section~2]{MR1278079} about the existence of Markovian bridges, that we apply to our $m$-symmetric diffusion $X.$ Under $P_x,$ the process $(p_{t-s}(X_s,y))_{s\in{[0,t)}}$ is a martingale, and thus we can define 
\begin{equation*}
    P_{x,y,t}(A)\stackrel{\mathrm{def.}}{=}\frac{E_{x}[p_{t-s}(X_s,y)\1_A]}{p_t(x,y)}\text{ for all }A\in{\mathcal{F}_s:=\sigma(X_u,u\leq s)}\text{ and }0\leq s<t,
\end{equation*}
and this definition is consistent. One can extend the definition of $P_{x,y,t}$ to a probability measure on $\F_t,$ which informally corresponds to the law of a bridge of length $t$ between $x$ and $y$ for $X.$ Applying the optional stopping theorem to the martingale $(p_{t-s}(X_s,y))_{s\in{[0,t)}},$ see for instance Theorem~3.2 in Chapter~II of \cite{MR1725357}, we have that for all $t>0$ and stopping times $T$
\begin{equation}
\label{stoppedbridge}
    E_x[p_{t-T}(X_T,y)\1_{A,T<t}]=P_{x,y,t}(A,T<t)p_t(x,y)\text{ for all }A\in{\F_T},
\end{equation}
where $\F_T=\{F\in{\mathcal{F}_t}:F\cap\{T\leq s\}\in\F_s\text{ for all }s<t\}$ is the filtration associated with $T.$ Moreover by $m$-symmetry of $X,$ we have for all $t>0$ and $x,y\in{\tilde{\G}}$ that
\begin{equation}
\label{symbridge}
   (X_{t-s})_{s\in{[0,t]}}\text{ has the same law under } P_{x,y,t}\text{ as }(X_s)_{s\in{[0,t]}}\text{ under }P_{y,x,t}.
\end{equation}
Using \eqref{symbridge}, one can derive a decomposition for stopping time on the reversed time scale: for all random times $\tau$ such that $\{\tau\geq t\}$ is in $\sigma(X_{t+u},u\geq0),$ we have that a.s.\
\begin{equation}
    \label{revmarkovbridge}
    (X_s)_{s\in{[0,\tau]}}\text{ has the same law under } P_{x}(\cdot\,|\,\mathcal{G}_{\tau})\text{ as }(X_s)_{s\in{[0,\tau]}}\text{ under }P_{x,X_{\tau},\tau},
\end{equation}
where $\mathcal{G}_{\tau}=\sigma(\tau,X_{\tau+u},u\geq0).$ Using results for general Hunt processes, see either \cite[Theorem~8]{MR336827}, \cite[Proposition~5.9]{MR334335} or \cite[Theorem~2.12]{MR521533}, under $P_x,$ if $L_F\in{(0,{\zeta})}$ then $(X_{s+L_F})_{s>0}$ is a Markov process depending on the past only through $X_{L_F},$ and so we have for all $x\in{\tilde{\G}},$ on the event $L_F\in{(0,{\zeta})},$ that
\begin{equation}
    \label{lawafterLK}
     (X_{s+L_F})_{s\geq0}\text{ has the same law under } P_{x}(\cdot\,|\,L_F,X_{L_F})\text{ as }(X_s)_{s\geq0}\text{ under }P_{X_{L_F}}^{F},
\end{equation}
where $P_{\cdot}^F$ is defined in \eqref{defPxF}. Combining \eqref{revmarkovbridge} and \eqref{lawafterLK}, one can thus describe the law of $(X_t)_{t\geq0}$ both before and after the last visit $L_F$ of $F.$ Let us now describe the law of $L_F$ and $X_{L_F}.$ Following the proof of \cite[(1.56)]{MR2932978} and \cite[\eqref*{4exitequi}]{DrePreRod3}, we moreover have that
\begin{equation}
    \label{exitequi}
    P_y(X_{L_F}=x,L_F\in{(0,\zeta)})={g(y,x)e_{F}(x)}\text{ for all }x,y\in{\tilde{\G}}.
\end{equation}
This leads to the following description of the law of $L_F$ and $X_{L_F}.$
\begin{lemma}
For all closed sets $F\subset\tilde{\G}$ and $x\in{\tilde{\G}}$ and $y\in{\hat{\partial}F}$ we have
\begin{equation}
    \label{disLK}
    P_x(L_F\in{\mathrm{d}t},L_F\in{(0,\zeta)},X_{L_F}=y)=p_t(x,y)e_{F}(y)\mathrm{d}t.
\end{equation}
\end{lemma}
\begin{proof}
For all $t>0,$ we have by the Markov property at time $t$ and \eqref{exitequi}
\begin{equation*}
    P_x(t< L_F<\zeta,X_{L_F}=y)=E_x\big[P_{X_t}(X_{L_F}=y,0<L_F<\zeta)\big]=E_x\big[g(X_t,y)\big]e_{F}(y).
\end{equation*}
Using \eqref{stoppedbridge}, we moreover have $E_x[p_{s-t}(X_t,y)]=p_s(x,y)$ for all $s>t,$ and so by \eqref{greenpt}
\begin{equation*}
    E_x\big[g(X_t,y)\big]=\int_{t}^{\infty}E_x\big[p_{s-t}(X_t,y)\big]\diff s=\int_{t}^{\infty}p_s(x,y)\diff s,
\end{equation*}
and we can conclude.
\end{proof}

We are now ready to give the last exit decomposition of $(X_t)_{t\geq0}$ before and after time $L_F.$ We denote by $W_{\tilde{\G}}^{+,f}$ the set of continuous trajectories in $\tilde{\G}$ with finite length, that is of continuous functions from $[0,t]$ to $\tilde{\G}$ for some $t>0.$ Let $\pi_t:\{w\in W_{\tilde{\G}}^+:\,t<\zeta\}\rightarrow W_{\tilde{\G}}^{+,f}$ be the application $w\mapsto w_{|[0,t]}$ and $\W_{\tilde{\G}}^{+,f}$ be the $\sigma$-algebra generated by $w\mapsto ((w(st))_{s\in{[0,1]}},t)$ when we endow $\{w:[0,1]\rightarrow\tilde{\G}:\,w\text{ is continuous}\}\times(0,\infty)$ with the product topology.  For each $A_1\in{\W_{\tilde{\G}}^{+,f}}$ and $A_2\in{\W_{\tilde{\G}}^+},$
using \eqref{revmarkovbridge} with $\tau=L_{F},$ we have that for all $x\in{\tilde{\G}}$ and $y\in{\hat{\partial} F}$
\begin{align*}
    &P_x\big((X_t)_{t\in{[0,L_F]}}\in{A_1},(X_{t+L_F})_{t\geq 0}\in{A_2},X_{L_F}=y,L_F\in{(0,\zeta)}\big)
    \\&=E_x\Big[\1_{(X_{t+L_{F}})_{t\geq 0}\in{A_2},X_{L_{F}}=y,L_F\in{(0,\zeta)}}(P_{x,y,L_{F}}\circ\pi_{L_F}^{-1})(A_1)\Big]
    \\&=P^{F}_y(A_2)E_x\Big[\1_{X_{L_{F}}=y,,L_F\in{(0,\zeta)}}(P_{x,y,L_{F}}\circ\pi_{L_F}^{-1})(A_1)\Big],
\end{align*}
where we used \eqref{lawafterLK} in the last equality. By \eqref{disLK}, we moreover have that
\begin{equation*}
    E_x\Big[\1_{X_{L_{F}}=y,L_F\in{(0,\zeta)}}(P_{x,y,L_{F}}\circ\pi_{L_F}^{-1})(A_1)\Big]=e_{F}(y)\int_{0}^\infty (P_{x,y,s}\circ\pi_{s}^{-1})(A_1)p_s(x,y)\diff s.
\end{equation*}
Summing over $y$ in $\hat{\partial} F,$ see \eqref{defpartialext}, we thus obtain the following last exit-decomposition for all closed sets $F\subset\tilde{\G},$ $x\in{\tilde{\G}},$ $A_1\in{\W_{\tilde{\G}}^{+,f}}$ and $A_2\in{\W^{+}}$ 
\begin{equation}
\label{lastexitdec}
\begin{split}
   &P_x\big((X_t)_{t\in{[0,L_F]}}\in{A_1},(X_{t+L_F})_{t\geq 0}\in{A_2},L_F\in{(0,\zeta)}\big)
   \\&=\sum_{y\in{\hat{\partial} F}}e_{F}(y)P^{F}_y(A_2)\int_{0}^\infty P_{x,y,s}(\pi_s^{-1}(A_1))p_s(x,y)\diff s.
   \end{split}
\end{equation}
Theorem~\ref{nuexists} follows classically from the last exit decomposition \eqref{lastexitdec}, as we now explain.

\begin{proof}[Proof of Theorem~\ref{nuexists}]
Let us fix a compact $K$ and a closed set $F$ with $K\subset F.$ For all $A\in{\mathcal{W}_{K,\tilde{\G}}^0},$ let $A'=\{(w(t+H_{F}))_{t\in\R}:\,w\in{A}\},$ with the convention that $\omega(t+H_{F})=\Delta$ for all $t\in\R$ if $H_F=\zeta^-.$ In order to prove \eqref{definter}, it is enough to prove that 
\begin{equation}
\label{QK=QK'}
    Q_{K,\tilde{\G}}(A)=Q_{F,\tilde{\G}}(A')\text{ for all $A\in{\W_{K,\tilde{\G}}^0}$ such that $A'\in{\W_{F,\tilde{\G}}^0}$}.
\end{equation}
Indeed one can then define $\1_{W_{K,\tilde{\G}}^*}\nu=Q_{K,\tilde{\G}}\circ(p^*_{\tilde{\G}})^{-1}$ for all compacts $K$ of $\tilde{\G},$ and this definition is consistent by \eqref{QK=QK'}, and we can conclude by taking a sequence of compacts increasing to $\tilde{\G}.$ The uniqueness of $\nu$ is clear since $W_{K,\tilde{\G}}^*$ increases to $W^*_{\tilde{\G}}$ as $K$ increases to $\tilde{\G},$ and \eqref{QK=QK'} directly implies \eqref{definter}.

Let us now prove \eqref{QK=QK'}. Using \eqref{defPxF}, \eqref{defQK} and \eqref{exitequi} we have
\begin{equation*}
     Q_{K,\tilde{\G}}(A)=\sum_{x\in{\hat{\partial} K}}\frac{1}{g(x,x)}P_x(A^+)P_x\big((X_{t+L_K})_{t\geq 0}\in{A^-},X_{L_K}=x,L_K\in{(0,\zeta)}\big),
\end{equation*}
where $A^+$ is the forwards part of $A$ and $A^-$ its backwards part, as defined above \eqref{defQK}. Since $A'\in{\W_{F,\tilde{\G}}^0},$ taking $A^{\pm}=\{(w(t))_{t\in{[0,H_{K}]}}:\,\omega\in{A'}\},$ one can easily check that $L_K\in{(0,\zeta)}$ and $(X_{t+L_K})_{t\geq 0}\in{A^-}$ if and only if $0<L_K\leq L_F<\zeta,$ $(X_{t+L_{F}})_{t\geq 0}\in{(A')^-}$ and $(X_{-t+L_{F}})_{t\in{[0,L_{F}-L_K]}}\in{A^\pm}.$ Therefore using \eqref{lastexitdec} for $F$ and \eqref{symbridge}, we obtain that for all $x\in{\hat{\partial} K}$
\begin{align*}
    &P_x\big((X_{t+L_K})_{t\geq 0}\in{A^-},X_{L_K}=x,L_K\in{(0,\zeta)}\big)
    \\&=\sum_{y\in{\hat{\partial} F}}\hspace{-1mm}e_{F}(y)P^{F}_y\big((A')^-\big)\int_0^{\infty}\hspace{-2mm}P_{x,y,s}\big((X_{s-t})_{t\in{[0,s-L_K]}}\in{A^{\pm}},X_{L_K}=x,L_K\in{(0,s]}\big)p_s(x,y)\diff s
    \\&=\sum_{y\in{\hat{\partial} F}}\hspace{-1mm}e_{F}(y)P^{F}_y\big((A')^-\big)\int_0^{\infty}\hspace{-2mm}P_{y,x,s}\big((X_{t})_{t\in{[0,H_K]}}\in{A^{\pm}},X_{H_K}=x,H_K\in{[0,s)}\big)p_s(y,x)\diff s.
\end{align*}
Moreover by \eqref{stoppedbridge}, we can write
\begin{align*}
    &\int_0^{\infty}P_{y,x,s}\big((X_{t})_{t\in{[0,H_K]}}\in{A^\pm},X_{H_K}=x,H_K\in{[0,s)}\big)p_s(y,x)\diff s
    \\&=\int_{0}^\infty E_y\big[p_{s-H_K}(x,x)\1_{(X_{t})_{t\in{[0,H_K]}}\in{A^\pm},X_{H_K}=x,H_K\in{[0,s\wedge\zeta)}}\big]\diff s
    \\&=E_y\Big[\1_{(X_{t})_{t\in{[0,H_K]}}\in{A^\pm},X_{H_K}=x,H_K<\zeta}\int_{H_K}^{\infty}p_{s-H_K}(x,x)\diff s\Big]
    \\&=g(x,x)P_y\big((X_{t})_{t\in{[0,H_K]}}\in{A^\pm},X_{H_K}=x,H_K<\zeta),
\end{align*}
where we used \eqref{greenpt} in the last equality. Combining the previous equations, we thus obtain by the strong Markov property at time $H_K$ that
\begin{align*}
     Q_{K,\tilde{\G}}(A)&=\sum_{x\in{\hat{\partial} K},y\in{\hat{\partial} F}}\hspace{-3mm}e_{F}(y)P_x(A^+)P_y\big((X_{t})_{t\in{[0,H_K]}}\in{A^\pm},X_{H_K}=x,H_K<\zeta\big) P^{F}_y\big((A')^-\big)
     \\&=\sum_{x\in{\hat{\partial} K},y\in{\hat{\partial} F}}\hspace{-3mm}e_{F}(y)P_y\big((A')^+,X_{H_K}=x,H_K<\zeta\big) P^{F}_y\big((A')^-\big)
     \\&=Q_{F,\tilde{\G}}(A'),
\end{align*}
where we used in the second equality the fact that $(X_t)_{t\geq0}\in{(A')^+}$ if and only if $H_K<\zeta,$ $(X_t)_{t\in{[0,H_K]}}\in{A^{\pm}}$ and $(X_{t+H_K})_{t\geq0}\in{A^+},$ and we can conclude.
\end{proof}

\renewcommand*{\thetheorem}{B.\arabic{theorem}}

\section{Appendix: Proof of Proposition~\ref{corh}}
\label{subsec:capandkappa=0}

In this Appendix, we prove the link between the graph $\G_{\h}$ and the notion of $\h$-transform, which is presented in Proposition~\ref{corh}. We recall the definition of the Dirichlet form $\mathcal{E}_{\tilde{\G}}$ and of the domain $D(\tilde{\G},\tilde{m})$ for any measures $\tilde{m}$ on $\tilde{\G}$ from \eqref{Dirichlet} and above, and for simplicity let us write $f'\in{L^2(\tilde{\G},\tilde{m})}$ if $f_{|I_e}\in{W^{1,2}(I_e,\tilde{m}_{|I_e}})$ for all $e\in{E\cup G}$ and $(f',f')_{\tilde{m}}<\infty.$ We begin with the following essential lemma.

\begin{lemma}
\label{htransformdirichlet}
If $\h$ is an harmonic function on $\tilde{\G},$ then $\xi_\h(X)$ is an $(\h^2\cdot m)$-symmetric diffusion on $\tilde{\G}$ under $P_{\psi_\h(x)}^{\tilde{\G}_\h},$ $x\in{\tilde{\G}},$ with associated Dirichlet form $\mathcal{E}_{\tilde{\G}}(f\h,g\h)$ on $L^2(\tilde{\G},\h^2\cdot m)$ with domain $D(\tilde{\G},\h^2\cdot m).$
\end{lemma}
\begin{proof}
Let $m_{\h}$ be the Lebesgue measure on $\tilde{\G}_{\h}.$ The process $X_{\theta_{\h}^{X}(t)}$ corresponds to a time changed process by a PCAF with Revuz measure $(\h\circ\psi_{\h}^{-1})^4,$ as follows from \cite[(5.1.13)]{MR2778606} for $f=1.$ Therefore by \cite[Theorem~6.2.1 and~(6.2.2)]{MR2778606}, the Dirichlet form associated with the $((\h\circ\psi_{\h}^{-1})^{4}\cdot m_{\h})$-symmetric diffusion $(X_{\theta_{\h}^X(t)})_{t<(\theta_{\h}^X)^{-1}(\zeta)}$ under $P_x^{\tilde{\G}_{\h}},$ $x\in{\tilde{\G}_{\h}}$ is $\mathcal{E}_{\tilde{\G}_{\h}}(f,g)$ on $L^2(\tilde{\G}_{\h},(\h\circ\psi_{\h}^{-1})^{4}\cdot m_{\h})$ with domain $\{f\in{L^2(\tilde{\G}_{\h},(\h\circ\psi_{\h}^{-1})^4\cdot m_{\h})\cap\mathcal{C}_0(\tilde{\G}_{\h})}:\,f'\in{L^2}(\tilde{\G}_{\h},m_{\h})\},$ where $m_{\h}$ is the Lebesgue measure on $\tilde{\G}_{\h}.$

Let $m'_{\h}=((\h\circ\psi_{\h}^{-1})^4\cdot m_{\h})\circ\psi_{\h}.$ Following \cite[Section~13]{MR958914}, if $\tilde{T}_t,$ $t\geq0,$ is the semigroup on $L^2(\tilde{\G}_{\h},(\h\circ\psi_{\h}^{-1})^{4}\cdot m_{\h})$ associated with $(X_{\theta_{\h}^X(t)})_{t<(\theta_{\h}^X)^{-1}(\zeta)},$ then the semigroup on $L^2(\tilde{\G},m'_{\h})$ associated with $\xi_{\h}(X)$ is $f\mapsto \tilde{T}_t(f\circ\psi_{\h}^{-1})\circ\psi_{\h}.$ Therefore using \cite[Lemma~1.3.4]{MR2778606}, one can easily prove that the Dirichlet form associated with the $m'_{\h}$-symmetric diffusion $\xi_{\h}(X)$ under $P_{\psi_\h(x)}^{\tilde{\G}_{\h}},$ $x\in{\tilde{\G}},$ is $\mathcal{E}_{\tilde{\G}_{\h}}(f\circ\psi_{\h}^{-1},g\circ\psi_{\h}^{-1})$ on $L^2(\tilde{\G},m'_{\h})$ with domain $\{f\in{L^2(\tilde{\G},m'_{\h})\cap\mathcal{C}_0(\tilde{\G})}:\,(f\circ\psi_{\h}^{-1})'\in{L^2({\tilde{\G}_{\h}},m_{\h})}\}.$ Let us fix some $e\in{E\cup{G}},$ and consider $\psi_{|I_e}$ and $h_{|I_e}$ as functions on $I_e,$ that we can identify with $(0,\rho_e)$ if $e\in{E},$ or $[0,\rho_e)$ if $e\in{x}.$ Then using \eqref{hformulaonedges} we have $(\psi_{\h})'(x)=\h(x)^{-2}$ for all $x\in{I_e},$ and $(\psi_{\h}^{-1})'(x)=(\h\circ\psi_{\h}^{-1}(x))^{2}$ for all $x\in{\psi_{\h}(I_{e})},$ and so we have by substitution for all $e\in{E\cup{G}}$ and for any Borel sets $A\subset I_e$
\begin{align*}
    m'_{\h}(A)=((\h\circ\psi_{\h}^{-1})^{4}\cdot m_{\h})(\psi_{\h}(A))&=\int_{\psi_{\h}(A)}(\h\circ\psi_{\h}^{-1})^4\,\mathrm{d}m_{\h}
    \\&=\int_A \h^2\,\mathrm{d}m=(\h^2\cdot m)(A),
\end{align*}
and so $m'_{\h}=\h^2\cdot m.$ Moreover for any functions $f,g$ with $(f\circ\psi_{\h}^{-1})',(g\circ\psi_{\h}^{-1})'\in{L^2({\tilde{\G}_{\h}},m_{\h})}$ we have
\begin{align*}
    \int_{\psi_{\h}(I_e)}(f\circ\psi_{\h}^{-1})'(g\circ\psi_{\h}^{-1})'\diff m_{\h}&=\int_{\psi_{\h}(I_e)}(f'\circ\psi_{\h}^{-1})(g'\circ\psi_{\h}^{-1})(\h\circ\psi_{\h}^{-1})^{4}\diff m_{\h}
    \\&=\int_{I_e}f'g'\h^2\diff m.
\end{align*}
Therefore, the domain of the Dirichlet form associated with $\xi_{\h}(X)$ is $D(\tilde{\G},\h^2\cdot m).$ Integrating by parts and noting that $\h'=\h'_e$ is constant on $I_e$ by Definition \ref{defharmonic},2), we have
\begin{align*}
    \int_{I_e}f'g'\h^2\diff m&=\int_{I_e}(f\h)'(g\h)'\diff m-\h'_e\int_{I_e}(fg)'\h\diff m-(\h'_e)^2\int_{I_e}fg\diff m
    \\&=\int_{I_e}(f\h)'(g\h)'\diff m-\h'_e[fg\h]_{I_e}.
\end{align*}
Moreover, 
\begin{align*}
    \sum_{e\in{E}}\h'_e[fg\h]_{I_e}&=\frac{1}2\sum_{\substack{x,y\in{G}\\x\sim y}}\frac{\mathrm{d}\h(x+t\cdot I_{\{x,y\}})}{\mathrm{d}t}\big(f(y)g(y)\h(y)-f(x)g(x)\h(x)\big)
    \\&=-\sum_{\substack{x,y\in{G}\\x\sim y}}\frac{\mathrm{d}\h(x+t\cdot I_{\{x,y\}})}{\mathrm{d}t}f(x)g(x)\h(x),
\end{align*}
and
\begin{equation*}
    \sum_{x\in{G}}\h'_x[fg\h]_{I_x}=-\sum_{x\in{G}}\frac{\mathrm{d}\h(x+t\cdot I_{x})}{\mathrm{d}t}f(x)g(x)\h(x).
\end{equation*}
Therefore we obtain by \eqref{Dirichlet} that the process $\xi_{\h}(X)$ under $P_{\psi_\h(x)}^{\tilde{\G}_{\h}}$ is a $(\h^2\cdot m)$-symmetric diffusion, and its associated Dirichlet form on $L^2(\tilde{\G},\h^2\cdot m)$ is 
\begin{align*}
    \mathcal{E}_{\tilde{\G}_{\h}}(f\circ\psi_{\h}^{-1},g\circ\psi_{\h}^{-1})&=\sum_{e\in{E\cup G}}\int_{\psi_{\h}(I_e)}(f\circ\psi_{\h}^{-1})'(g\circ\psi_{\h}^{-1})'\diff m_{\h}
    \\&=\sum_{e\in{E\cup G}}\int_{I_e}(f\h)'(g\h)'\diff m
    \\&\quad+\sum_{x\in{G}}f(x)g(x)\h(x)\Big(\frac{\mathrm{d}\h(x+t\cdot I_{x})}{\mathrm{d}t}+\sum_{y\sim x}\frac{\mathrm{d}\h(x+t\cdot I_{\{x,y\}})}{\mathrm{d}t}\Big)
    \\&=\mathcal{E}_{\tilde{\G}}(f\h,g\h),
\end{align*}
where we used \eqref{defharmonic} in the last equality.
\end{proof}

Lemma~\ref{htransformdirichlet} let us thus compare the law of $\xi_{\h}(X)$ under the $P_{\psi_{\h}(\cdot)}^{\tilde{\G}_{\h}}$ and the law of $X$ under $P_{\cdot}^{\tilde{\G}},$ from which Proposition~\ref{corh} follows easily.

\begin{proof}[Proof of Proposition~\ref{corh}]
It follows from \cite[Lemma~1.3.4]{MR2778606} that the Dirichlet form on $L^2(\tilde{\G},\h^2\cdot m)$ associated with the semigroup $\frac{1}{\h}T_t^{\tilde{\G}}(f\h)$ is
\begin{equation*}
    \lim\limits_{t\searrow0}\frac1t\big(f-\frac1{\h}T_t^{\tilde{\G}}(f\h),g\big)_{\h^2\cdot m}=\lim\limits_{t\searrow0}\frac1t(f\h-T_t^{\tilde{\G}}(f\h),g\h)_{m}=\mathcal{E}_{\tilde{\G}}(f\h,g\h),
\end{equation*}
with domain $D(\tilde{\G},\h^2\cdot m),$ and is thus the semigroup associated with $\xi_{\h}(X)$ by Lemma~\ref{htransformdirichlet}.

Let us now turn to the proof of \eqref{eq:localtimesh}. Following \cite[Section~2]{MR3502602}, we have for all $x\in{\tilde{\G}}$ and $t\geq0$ that $P^{\tilde{\G}_{\h}}_{\cdot}$-a.s.
    \begin{equation*}
        \ell_{\psi_{\h}(x)}(t)=\lim_{\eps\rightarrow0}\frac{1}{m_{\h}(\psi_{\h}(\tilde{B}(x,\eps)))}\int_{0}^t\1_{X_u\in{\psi_{\h}(\tilde{B}(x,\eps))}}\diff u,
    \end{equation*}
    where $\tilde{B}(x,\eps)=\{x+t\cdot I_e\in{\tilde{\G}}:t\in{[0,\eps]}\text{ and }e\in{E\cup G}\text{ with }x\in{\overline{I_e}}\}.$ Taking $u=\theta_{\h}^{X}(s),$ we have
    \begin{equation*}
        \int_{0}^t\1_{(\xi_{\h}(X))_s\in{\tilde{B}(x,\eps)}}\diff s=\int_{0}^{\theta_{\h}^{X}(t)}\1_{X_u\in{\psi_{\h}(\tilde{B}(x,\eps))}}\h\big(\psi_{\h}^{-1}(X_u)\big)^4\diff u,
    \end{equation*}
    and one can easily check that
    \begin{equation*}
        \frac{m_{{\h}}\big(\psi_{\h}(\tilde{B}(x,\eps))\big)}{m\big(\tilde{B}(x,\eps)\big)}\tend{\eps}{0}\frac{1}{\h(x)^2}.
    \end{equation*}
    We thus obtain that 
    \begin{align*}
        &\lim_{\eps\rightarrow0}\frac{1}{m(\tilde{B}(x,\eps))}\int_{0}^t\1_{(\xi_{\h}(X))_s\in{\tilde{B}(x,\eps)}}\diff s
        \\&=\lim_{\eps\rightarrow0} \frac{1}{\h(x)^2m_{{\h}}\big(\psi_{\h}(\tilde{B}(x,\eps))\big)}\int_{0}^{\theta_{\h}^{X}(t)}\1_{X_u\in{\psi_{\h}(\tilde{B}(x,\eps))}}\h\big(\psi_{\h}^{-1}(X_u)\big)^4\diff u
        \\&=\h(x)^2\ell_{\psi_{\h}(x)}(\theta_{\h}^{X}(t)).
    \end{align*}

Let us finally prove \eqref{eq:GFFh}. According to \eqref{eq:localtimesh}, with respect to the measure $(\h^2\cdot m),$ the field of local times associated with $\xi_{\h}(X)$ is $\big(\ell_{\psi_{\h}(x)}(\theta_{\h}^{X}(t))\big)_{t\geq0,x\in{\tilde{\G}}},$ and so the associated Green function is $g_{\tilde{\G}_{\h}}(\psi_{\h}(x),\psi_{\h}(y)),$ $x,y\in{\tilde{\G}}.$ Moreover, by \eqref{semigrouph} and \eqref{greenpt}, with respect to $(\h^2\cdot m),$ the family of probability densities associated with $\xi_{\h}(X)$ is $p_t(x,y)/(\h(x)\h(y)),$ $x,y\in{\tilde{\G}},$ and so by \eqref{greenpt} its Green function is $g_{\tilde{\G}}(x,y)/(\h(x)\h(y)),$ $x,y\in{\tilde{\G}}.$ We thus obtain that
\begin{equation}
\label{relationGreenfunctionh}
    \h(x)\h(y)g_{\tilde{\G}_{\h}}(\psi_{\h}(x),\psi_{\h}(y))=g_{\tilde{\G}}(x,y)\text{ for all }x,y\in{\tilde{\G}}.
\end{equation}
In particular, the two processes in \eqref{eq:GFFh} are centered Gaussian processes with the same covariance function, and we can conclude.
\end{proof}

\bibliography{bibliographie}
\bibliographystyle{abbrv}
\end{document}